\documentclass[preprint,12pt]{elsarticle}

\usepackage{lineno,hyperref}
\usepackage{amssymb, epsfig,amssymb, latexsym}
\usepackage{amsfonts,psfrag,amsmath,bbm,color,url}

\usepackage[ruled,vlined]{algorithm2e}

\def\PP{{{\rm l}\kern - .15em {\rm P} }}
\def\PN2{{\PP_{N}-\PP_{N-2}}}

\newcommand{\R}{\mathbbm{R}}

\newcommand{\cO}{\mathcal{O}}

\newcommand{\bphi}{\boldsymbol{\varphi}}

\newcommand{\btau}{\boldsymbol{\tau}}

\newcommand{\ba}{\boldsymbol{a}}

\newcommand{\bff}{\boldsymbol{f}}
\newcommand{\bFF}{{\boldsymbol F}}

\newcommand{\bH}{\boldsymbol{H}}

\newcommand{\bu}{\boldsymbol{u}}

\newcommand{\bur}{{\boldsymbol{u}}_r}

\newcommand{\bv}{\boldsymbol{v}}

\newcommand{\bx}{\boldsymbol{x}}
\newcommand{\bX}{\boldsymbol{X}}

\newcommand{\bXr}{{\bf X}^r}

\newcommand{\tA}{\tilde{A}}
\newcommand{\tB}{\widetilde{B}}

\newcommand{\deleted}[1]{{}}

\usepackage{tikz}
\usetikzlibrary{fit,positioning}
\usepackage{siunitx}
\usepackage{float}

\let\linenumbers\nolinenumbers\nolinenumbers

\usepackage{algorithmic}

\journal{arXiv}

\begin{document}

\begin{frontmatter}

\title{Data-Driven Variational Multiscale \\ Reduced Order Models}

 \author[xx]{Changhong Mou}
 \ead{cmou@vt.edu}
 \ead[url]{https://personal.math.vt.edu/people/cmou}
 \author[xx]{Birgul Koc}
 \ead{birgul@vt.edu}
 \ead[url]{https://personal.math.vt.edu/people/birgul}
 \author[os]{Omer San}
 \ead{osan@okstate.edu}
 \ead[url]{https://mae.okstate.edu/content/san-omer}
  \author[clemson]{Leo G. Rebholz}  
 \ead{rebholz@clemson.edu}
 \ead[url]{https://www.math.clemson.edu/~rebholz}
 \author[xx]{Traian Iliescu\corref{cor1}}
 \ead{iliescu@vt.edu}
 \ead[url]{https://personal.math.vt.edu/people/iliescu}
 \cortext[cor1]{corresponding author}
 \address[xx]{Department of Mathematics, Virginia Tech, Blacksburg, VA 24061, U.S.A.}
\address[os]{Department of Mechanical and Aerospace Engineering, Oklahoma State University, Stillwater, OK 74078, U.S.A.}
\address[clemson]{Department of Mathematical Sciences, Clemson University, Clemson, SC, 29634, U.S.A.}

\begin{abstract}
We propose a new data-driven reduced order model (ROM) framework that centers around the hierarchical structure of the variational multiscale (VMS) methodology and utilizes data to  increase the ROM accuracy at a modest computational cost.
The VMS methodology is a natural fit for the hierarchical structure of the ROM basis:
In the  first step, we use the ROM projection to separate the scales into three categories: (i) resolved large scales, (ii) resolved small scales, and (iii) unresolved scales. 
In the second step, we explicitly identify the VMS-ROM closure terms, i.e., the terms representing the interactions among the three types of scales.
In the third step, 
we use available data to model the VMS-ROM closure terms.
Thus, instead of phenomenological models used in VMS for standard numerical discretizations (e.g., eddy viscosity models), we utilize available data to construct new structural VMS-ROM closure models.
Specifically, we build ROM operators (vectors, matrices, and tensors) that are closest to the true ROM closure terms evaluated with the available data. 
We test the new data-driven VMS-ROM in the numerical simulation of
{four test cases:
(i) the 1D Burgers equation with  viscosity  coefficient $\nu = 10^{-3}$; 
(ii) a 2D flow past a circular cylinder at Reynolds numbers $Re=100$, $Re=500$, and $Re=1000$;
(iii) the quasi-geostrophic equations at Reynolds number $Re=450$ and Rossby number $Ro=0.0036$; and 
(iv) a 2D flow over a backward facing step at Reynolds number $Re=1000$.
}
The numerical results show that the data-driven VMS-ROM is significantly more accurate than standard ROMs.
\end{abstract}



\end{frontmatter}

\linenumbers

\section{Introduction}

For structure dominated systems, {\it reduced order models (ROMs)}~\cite{crommelin2004strategies,gunzburger2017ensemble,hesthaven2015certified,HLB96,noack2011reduced,perotto2017higamod,quarteroni2015reduced,sapsis2009dynamically,stefanescu2015pod,taira2019modal} can decrease the full order model (FOM) computational cost by orders of magnitude.
ROMs are low-dimensional models that are constructed from available data: 
In an offline stage, the FOM is run for a small set of parameters to construct a low-dimensional ROM basis $\{ \varphi_{1}, \ldots, \varphi_{r}\}$, which is used to build the ROM:
\vspace*{-0.05cm}
\begin{eqnarray}
	\overset{\bullet}{\ba}
	=
	\bFF(\ba),
	\label{eqn:eqn-rom} 
	\\[-0.55cm]
	\nonumber
\end{eqnarray}
where $\ba$ is the vector of coefficients in the ROM approximation $\sum_{i=1}^{r} a_{i}(t) \varphi_{i}(\bx)$ of the variable of interest and $\bFF$ comprises the ROM operators (e.g., vectors, matrices, and tensors) that are preassembled from the ROM basis in the offline stage.
In the online stage, the low-dimensional ROM~\eqref{eqn:eqn-rom} is then used in a regime that is different from the training regime.
Since the ROM~\eqref{eqn:eqn-rom} is low-dimensional, its computational cost is orders of magnitude lower than the FOM cost.

Unfortunately, current ROMs cannot be used in complex, realistic settings, since they require too many modes (degrees of freedom).
For example, to capture all the relevant scales in practical engineering flows, hundreds~\cite{osth2014need,shah2014very} and even thousands of ROM modes can be necessary~\cite{osth2014need,zhang2020characterizing}.
Thus, although ROMs decrease the FOM computational cost by orders of magnitude, they cannot be used in many important practical settings (e.g., digital twin applications, where a real-time control of physical assets may be required~\cite{hartmann2018model}). 

One of the main roadblocks in the development of ROMs for complex practical settings is their notorious inaccuracy. The drastic ROM truncation is one of the most important reasons for the ROMs' numerical inaccuracy:
Instead of using a sufficient number of ROM modes $\{ \varphi_{1}, \ldots, \varphi_{R}\}$ to  capture the dynamics of the  underlying system, current ROMs use only a handful of ROM modes $\{ \varphi_{1}, \ldots, \varphi_{r}\}$ to ensure a low computational cost.
This drastic truncation yields acceptable results in simple, academic test problems, but yields inaccurate results in many practical settings~\cite{osth2014need}, where the {\it ROM closure problem}~\cite{baiges2019finite,balajewicz2016minimal,benosman2017learning,chekroun2019variational,fick2018stabilized,harlim2019machine,HLB96,lu2017data,majda2018model,majda2012physics,pagani2019statistical,rebollo2017certified,stabile2019reduced,wan2018data,wang2012proper} needs to be solved:
One needs to model the effect of the discarded ROM modes $\{ \varphi_{r+1}, \ldots, \varphi_{R}\}$ on the ROM dynamics, i.e., on the time evolution of the resolved ROM modes $\{ \varphi_{1}, \ldots, \varphi_{r}\}$:
\begin{eqnarray}
	\overset{\bullet}{\ba} 
	=
	\bFF(\ba)
	+
	\text{Closure}(\ba),
	\label{eqn:eqn-rom-closure}
\end{eqnarray}
where $\text{Closure}(\ba)$ is a low-dimensional term that models the effect of the discarded ROM modes $\{ \varphi_{r+1}, \ldots, \varphi_{R}\}$ on $\{ \varphi_{1}, \ldots, \varphi_{r}\}$. 

The closure problem is ubiquitous in the numerical simulation of complex systems.
For example, classical numerical discretization of turbulent flows (e.g., finite element or finite volume methods), inevitably takes place in the {\it under-resolved regime} (e.g., on coarse meshes) and requires closure modeling (i.e., modeling the sub-grid scale effects).
In classical CFD, e.g., large eddy simulation (LES), there are hundreds (if not thousands) of closure models~\cite{sagaut2006large}.
This is in stark contrast with ROM, where only relatively few ROM closure models have been investigated.
The reason for the discrepancy between ROM closure and LES closure is that the latter has been entirely built around physical insight stemming from Kolmogorov's statistical theory of turbulence (e.g., the concept of eddy viscosity), which is generally posed in the Fourier setting~\cite{sagaut2006large}.
This physical insight is generally not available in a ROM setting.
Thus, current ROM closure models have generally been deprived of many tools of this powerful methodology that represents the core of most LES closure models.
Since physical insight cannot generally be used in the ROM setting, alternative ROM closure modeling strategies need to be developed.
Our vision is that {\it data} represents a natural solution for ROM closure modeling.

In this paper, we put forth a new ROM framework that centers around the hierarchical structure of {\it variational multiscale (VMS)} methodology~\cite{hughes1998variational,hughes2000large,hughes2001multiscale,hughes2001large}, which naturally separates the scales into (i) resolved large, (ii) resolved small, and (iii) unresolved.
We also construct new structural ROM closure models for the three scales by using available data.
We believe that the VMS methodology is a natural fit for the hierarchical structure of the ROM basis:
In the first step of the new VMS-ROM framework, we use the ROM projection to unambiguously separate the scales into three categories: (i) {\it resolved large} scales, (ii) {\it resolved small} scales, and (iii) {\it unresolved} scales. 
In the second step, we explicitly identify the {\it ROM closure} terms representing the interactions among the three types of scales by projecting the equations onto the corresponding resolved large, resolved small, and unresolved spaces.
In the third step, instead of phenomenological modeling techniques used in VMS for standard discretizations (e.g., finite element methods), we utilize {\it data-driven modeling}~\cite{brunton2019data,loiseau2018constrained,peherstorfer2016data} to construct novel, robust, {\it structural} ROM closure models.
Thus, instead of ad hoc, phenomenological models used in VMS for standard numerical discretizations (e.g., eddy viscosity models), we utilize available data to construct new structural models for the interaction among the three types of scales.
Specifically, we use FOM data to develop VMS-ROM closure terms that account for the {\it under-resolved} numerical regime.
We emphasize that, in the new {\it data-driven VMS-ROM (DD-VMS-ROM)} framework, we use data only to {\it  complement} classical physical modeling (i.e., only for closure modeling)~\cite{maulik2019time,xie2018data}, not to completely replace it~\cite{brunton2019data,rahman2019non}.
Thus, the resulting ROM framework combines the strengths of both physical and data-driven modeling.

\paragraph{Previous Relevant Work}

The VMS methodology has been used in ROM settings~\cite{bergmann2009enablers,eroglu2017modular,guler2019decoupled,iliescu2013variational,iliescu2014variational,roop2013proper,stabile2019reduced,wang2012proper}.
We emphasize, however, that the DD-VMS-ROM framework that we propose is different from the other VMS-ROMs.

The VMS-ROMs in~\cite{bergmann2009enablers,iliescu2013variational,iliescu2014variational,stabile2019reduced,wang2012proper} are {\it phenomenological} models in which the role of the VMS closure models is to {\it dissipate energy} from the ROM. 
In contrast, the new DD-VMS-ROM utilizes data to construct general structural VMS-ROM closure terms, which are {\it not required to be dissipative}.
(Of course, if deemed appropriate, we may impose additional constraints to mimic the physical properties of the underlying system~\cite{mohebujjaman2019physically}.)

The new DD-VMS-ROM is also different from the reduced-order subscales ROM proposed  in~\cite{baiges2015reduced} (see also~\cite{reyes2020projection,reyes2018reduced,tello2019fluid}):
The reduced-order subscales model in~\cite{baiges2015reduced} minimizes the difference between the {\it solutions} of the FOM and ROM (see equations (18)--(19) in [6]), whereas
the new DD-VMS-ROM minimizes the difference between the VMS-ROM {\it closure terms} and the ``true'' (i.e., high-resolution) closure terms.
Furthermore, the reduced-order subscales model in~\cite{baiges2015reduced} builds {\it linear} closure models (see also~\cite{oberai2016approximate}), whereas the new DD-VMS-ROM constructs {\it nonlinear} closure models.

Another ROM closure strategy that is related to the VMS-ROM framework is the adjoint Petrov-Galerkin method~\cite{parish2019adjoint} (see~\cite{carlberg2017galerkin,carlberg2011efficient,grimberg2020stability} for related work), which is based on the Mori-Zwanzig (MZ) formalism~\cite{gouasmi2017priori,lin2019data}.
In the MZ-ROM approach, the ROM closure model is represented by a memory term that depends on the temporal history of the resolved scales.
The memory term is approximated to construct effective ROM closure models and, therefore, practical ROMs.
The main difference between the adjoint Petrov-Galerkin method proposed in~\cite{parish2019adjoint} and the new DD-VMS-ROM is the tool used to define the ROM closure term:
The former uses a statistical tool (i.e., the MZ formalism), whereas the latter utilizes a spectral-like projection (i.e., the ROM projection).   

Finally, we note that the  VMS-ROM framework proposed herein belongs to the wider class of {\it  hybrid physical/data-driven ROMs}, in which data-driven modeling is used to model only the missing information (i.e., the ROM closure term) in ROMs constructed from first principles (i.e., from a  Galerkin projection of the underlying equations); see, e.g.,~\cite{baiges2019finite,chekroun2019variational,couplet2005calibrated,galletti2004low,hijazi2019data,lu2017data,noack2005need,parish2019adjoint,xie2018data}.

\bigskip

The rest of the paper is organized as follows:
In Section~\ref{sec:vms-rom}, we introduce the new DD-VMS-ROM.
In Section~\ref{sec:numerical-results}, we test the DD-VMS-ROM in the numerical simulation of  
{four test cases:
(i) the 1D Burgers equation with  viscosity  coefficient $\nu = 10^{-3}$; 
(ii) a 2D flow past a circular cylinder at Reynolds numbers $Re=100$, $Re=500$, and $Re=1000$;
(iii) the quasi-geostrophic equations at Reynolds number $Re=450$ and Rossby number $Ro=0.0036$; and 
(iv) a 2D flow over a backward facing step at Reynolds number $Re=1000$.
}
Finally, in Section~\ref{sec:conclusions}, we draw conclusions and outline future research directions.

\section{Data-Driven Variational Multiscale Reduced Order Models (DD-VMS-ROMs)}
    \label{sec:vms-rom}

In this section, we construct the new {\it data-driven VMS-ROM (DD-VMS-ROM)} framework, which can significantly increase the accuracy of under-resolved ROMs, i.e., ROMs whose dimension is too low to capture the complex dynamics of realistic applications. 
In Section~\ref{sec:vms-rom-vms} we briefly sketch the VMS methodology for general numerical discretizations (see, e.g.,~\cite{ahmed2017review,john2016finite} for more details), and in Section~\ref{sec:vms-rom-g-rom} we outline the standard Galerkin ROMs.

We construct the new DD-VMS-ROM in two stages:
In Section~\ref{sec:vms-rom-two-scales}, we construct the two-scale DD-VMS-ROM, which is the simplest DD-VMS-ROM. 
We note that the two-scale data-driven VMS-ROM was investigated in~\cite{xie2018data} under the name ``data-driven filtered ROM'' and in~\cite{mohebujjaman2019physically,mou2019data} under the name ``data-driven correction ROM.''
However, we decided to outline the construction of the two-scale data-driven VMS-ROM since it is the most straightforward illustration of the DD-VMS-ROM framework. 

In Section~\ref{sec:vms-rom-three-scales}, we construct the novel three-scale DD-VMS-ROM. This new model separates the scales into three categories (instead of two, as in the two-scale DD-VMS-ROM), which allows more flexibility in constructing the ROM closure models and could lead to more accurate ROMs.. 

\subsection{Classical VMS}
    \label{sec:vms-rom-vms}

The VMS methods are general numerical discretizations that increase the {\it accuracy} of classical Galerkin approximations in {\it under-resolved} simulations, e.g., on coarse meshes or when not enough basis functions are available.
The VMS framework, which was proposed by Hughes and coworkers~\cite{hughes1998variational,hughes2000large,hughes2001multiscale,hughes2001large}, has made a profound impact in several areas of computational mathematics (see, e.g.,~\cite{ahmed2017review,codina2018variational,john2016finite,rasthofer2018recent} for surveys). 
To illustrate the standard VMS methodology, we consider a general nonlinear system/PDE
\vspace*{-0.1cm}
\begin{eqnarray}
	\overset{\bullet}{\bu}
	= 
	\bff(\bu) \, ,
	\label{eqn:pde-strong}
\vspace*{-0.1cm}
\end{eqnarray}
whose weak (variational) form is  
\begin{eqnarray}
	( \overset{\bullet}{\bu} , \bv ) 
	= 
	( \bff(\bu) , \bv ), \  
	\qquad
	\forall \, \bv \in \bX ,
 	\label{eqn:pde-weak}
 \end{eqnarray}
where $\bff$ is a general nonlinear function and $\bX$ is an appropriate infinite dimensional space.
To build the VMS framework, we start with a sequence of {\it hierarchical spaces} of increasing resolutions:
$\bX_1, \bX_1 \oplus\bX_2, \bX_1 \oplus\bX_2 \oplus \bX_3, \ldots$.
Next, we project system~\eqref{eqn:pde-strong} onto {\it each} of the spaces $\bX_1, \bX_2, \bX_3, \ldots$, which yields a separate equation for each space.
The goal is, of course, to solve for the $\bu$ component that lives in the coarsest space (i.e., $\bX_1$), since this yields the lowest-dimensional system:
\begin{eqnarray}
	( \overset{\bullet}{\bu} , \bv_1 ) 
	= 
	\bigl( \bff(\bu) , \bv_1 \bigr) 
	\qquad
	\forall \, \bv_1 \in \bX_1 \, .
	\label{eqn:vms-1}
\end{eqnarray}
System~\eqref{eqn:vms-1}, however, is {\it not} closed, since its right-hand side 
\begin{eqnarray}
	\bigl( \bff(\bu) , \bv_1 \bigr) 
	=
	\bigl( \bff(\bu_1 + \bu_2 + \bu_3 + \ldots) , \bv_1 \bigr)
	\qquad
	\forall \, \bv_1 \in \bX_1 \, ,
	\label{eqn:vms-2}
\end{eqnarray}
involves $\bu$ components that  do {\it not} live in $\bX_1$ (i.e., $\bu_2 \in \bX_2, \bu_3  \in \bX_3, \ldots$).
This coupling is mainly due to the {\it nonlinearity} of $\bff$.
Thus, the {\it VMS closure problem} needs to be solved, i.e.,~\eqref{eqn:vms-2} needs to be approximated in $\bX_1$. 
The VMS~\eqref{eqn:vms-1} equipped with an appropriate closure model yields an accurate approximation of the large scale $\bX_1$ component of $\bu$. \\[-0.1cm]

The main reasons for the VMS framework's impressive success are its {\it  utter simplicity} and its {\it generality} (it can be applied to {\it any} Galerkin based numerical discretization).
The classical VMS methodology, however, is facing several major challenges:
(i) The {\it hierarchical spaces} can be difficult to construct in classical Galerkin methods (e.g., finite elements); and 
(ii) Developing {\it VMS closure} models for the coupling terms (i.e., the terms that model the interactions among scales) can be challenging. 

In this  paper, we propose a new data-driven VMS-ROM framework that overcomes these major challenges of standard VMS methodology:
(i) The ROM setting allows a natural, straightforward construction of ROM hierarchical spaces.
(ii) We use available data to construct data-driven VMS-ROM closure models.
Thus, we avoid the ad hoc assumptions and phenomenological arguments that are often used in traditional VMS closures.

\subsection{Galerkin ROM (G-ROM)}
    \label{sec:vms-rom-g-rom}

Before building the new VMS-ROM framework, we sketch the standard Galerkin ROM derivation:
(i) Use available data (snapshots) for few parameter values to construct orthonormal modes $\{ \bphi_1, \ldots, \bphi_R \}$, $R = \cO(10^{3})$, which represent the recurrent spatial structures;
(ii) Choose the dominant modes $\{ \bphi_1, \ldots, \bphi_r \}$, $r= \cO(10)$, as  basis functions for the ROM;
(iii) Use a Galerkin truncation $\bur(\bx,t) = \sum_{j=1}^{r} a_j(t) \, \bphi_j(\bx)$;
(iv) Replace $\bu$ with $\bur$ in~\eqref{eqn:pde-strong};
(v) Use a Galerkin projection of the PDE obtained in step (iv) onto the ROM space $\bXr := \text{span} \{ \bphi_1, \ldots, \bphi_r \}$ to obtain an $r$-dimensional  system, which is the {\it Galerkin ROM (G-ROM)}: 
\vspace*{-0.2cm}
\begin{eqnarray}
	\bigl( \overset{\bullet}{\bu_{r}} \, , \bphi_{i} \bigr) 
	= \bigl( \bff(\bu_{r}) \, , \bphi_{i} \bigr) \, ,
	\qquad 
	i = 1, \ldots, r; 	
	\label{eqn:g-rom-general}
	\\[-0.6cm]
	\nonumber
\end{eqnarray}
(vi) In an offline stage, compute the ROM operators; 
(vii) In an online stage, repeatedly use the G-ROM~\eqref{eqn:g-rom-general} (for parameters  different from the training parameters  and/or longer time intervals).

We illustrate the G-ROM for the Navier-Stokes equations (NSE):
\vspace*{-0.1cm}
 \begin{eqnarray}
     && \frac{\partial \bu}{\partial t}
     - Re^{-1} \Delta \bu
     + \bu \cdot \nabla \bu
     + \nabla p
     = {\bf 0} \, ,
     \label{eqn:nse-1} \\
     && \nabla \cdot \bu
     = 0 \, ,
     \label{eqn:nse-2}
	\\[-0.6cm]
	\nonumber
 \end{eqnarray}
 where $\bu$ is the velocity, $p$ the pressure, and $Re$ the Reynolds number. 
For clarity of presentation, we use homogeneous Dirichlet boundary conditions.
The NSE~\eqref{eqn:nse-1}--\eqref{eqn:nse-2} can be cast in the general form~\eqref{eqn:pde-strong} by choosing $\bff = Re^{-1} \Delta \bu - \bu \cdot \nabla \bu$ and $\bX$ the space of weakly divergence-free functions in $\bH_{0}^{1}$.
For the NSE, the G-ROM reads
\vspace*{-0.1cm}
\begin{eqnarray}
	\overset{\bullet}{\ba} 
	= A \, \ba + \ba^{\top} \, B \, \ba ,
\label{eqn:g-rom}
	\\[-0.6cm]
	\nonumber
\end{eqnarray}
where 
$\ba(t)$ is the vector of unknown coefficients $a_j(t), 1 \leq j \leq r$, 
$A$ is an $r \times r$ matrix with entries 
$
	A_{im}
	= - Re^{-1} \, \left( \nabla \bphi_m , \nabla \bphi_i \right) ,
$
and $B$ is an $r \times r \times r$ tensor with entries
$
	B_{imn}
	= - \bigl( \bphi_m \cdot \nabla \bphi_n , \bphi_i \bigr) \, ,
	\ 
	1 \leq i, m, n \leq r \, .
$
The G-ROM~\eqref{eqn:g-rom} does not include a pressure approximation, since we assumed that the ROM modes are discretely divergence-free (which is the case if, e.g., the snapshots are discretely divergence-free).
ROMs that provide a pressure approximation are discussed in, e.g.,~\cite{decaria2020artificial,hesthaven2015certified,quarteroni2015reduced}.
Once the matrix $A$ and tensor $B$  are assembled in the offline stage, the G-ROM~\eqref{eqn:g-rom} is a low-dimensional, efficient dynamical system that can be used in the online stage for numerous parameter values.
We emphasize, however, that the G-ROM generally yields inaccurate results when used in {\it under-resolved}, realistic, complex flows~\cite{HLB96,noack2011reduced,osth2014need,wang2012proper}.

\subsection{Two-Scale Data-Driven Variational Multiscale ROMs (2S-DD-VMS-ROM)}
    \label{sec:vms-rom-two-scales}

The first DD-VMS-ROM that we outline is the {\it two-scale data-driven VMS-ROM (2S-DD-VMS-ROM)}, which utilizes two orthogonal spaces, $\bX_1$ and $\bX_2$.
Since the ROM basis is orthonormal by construction,  we can build the two orthogonal spaces in a natural way:
$\bX_1 := \text{span}  \{ \bphi_1, \ldots, \bphi_{r} \}$, which represents the resolved ROM scales, 
and 
$\bX_2 := \text{span}  \{ \bphi_{r+1}, \ldots, \bphi_{R} \}$, which represents the unresolved ROM scales.  
We note that, in practical settings, we are forced to use {\it under-resolved} ROMs, i.e., ROMs whose dimension $r$ is much lower than the dimension of the snapshot data set (i.e., $R$).
Next, we use the best ROM approximation of $\bu$ in the space $\bX_1 \oplus \bX_2$,  i.e.,   $\bu_{R} \in \bX_1 \oplus \bX_2$  defined as
\vspace*{-0.1cm}
\begin{eqnarray}
	\bu_{R}
	= \sum_{j=1}^{R} a_{j} \, \bphi_{j}
	= \sum_{j=1}^{r} a_{j} \, \bphi_{j}
	+ \sum_{j=r+1}^{R} a_{j} \, \bphi_{j} = \bu_{r}
	+ \bu' \, ,
	\label{eqn:vms-rom-1}
\end{eqnarray}
where $\bu_{r} \in \bX_1$ represents the resolved ROM component of $\bu$, and $\bu' \in \bX_2$ represents the unresolved ROM component of $\bu$.
Plugging $\bu_R$ in~\eqref{eqn:pde-strong}, projecting the resulting equation onto $\bX_1$, and using the ROM basis orthogonality  to show that  
$\bigl( \overset{\bullet}{\bu_{R}} \, , \bphi_{i} \bigr) = \bigl( \overset{\bullet}{\bu_{r}} \, , \bphi_{i} \bigr), \  \forall \, i = 1, \ldots, r$, we obtain
\begin{eqnarray}
		\bigl( \overset{\bullet}{\bu_{r}} \, , \bphi_{i} \bigr) 
		= 
		\bigl( \bff(\bu_{r}) \, , \bphi_{i} \bigr) 
		+
		\underbrace{
		\boxed{
		\bigl[ 
			\bigl( \bff(\bu_{R}) \, , \bphi_{i} \bigr)
			-
			\bigl( \bff(\bu_{r}) \, , \bphi_{i} \bigr)
		\bigr] 
		} 
		}_{\text{VMS-ROM closure term}}
		\ , 
		\  
		\forall \, i = 1, \ldots, r .
	\label{eqn:vms-rom-two-scales-1}
\end{eqnarray}
The boxed term in~\eqref{eqn:vms-rom-two-scales-1} is the {\it VMS-ROM closure term}, which models the interaction between the ROM modes $\{ \bphi_1, \ldots, \bphi_r \}$ and the discarded ROM modes $\{ \bphi_{r+1}, \ldots, \bphi_{R} \}$. 
The VMS-ROM closure term is essential for the accuracy of~\eqref{eqn:vms-rom-two-scales-1}:
If we drop the VMS-ROM closure term, we are left with the  G-ROM~\eqref{eqn:g-rom-general}, which yields inaccurate results in the under-resolved regime.
The VMS-ROM closure term is a {\it correction term} that ensures an accurate approximation of $\bu_r \in \bX_1$ in the higher-dimensional space $\bX_1 \oplus \bX_2$.

Next, we approximate the VMS-ROM closure term with ${\boldsymbol{g}}(\bu_r)$, where ${\boldsymbol{g}}$ is a {\it generic} function whose coefficients/parameters still need to be determined:
\begin{eqnarray}
	\text{VMS-ROM closure term} 
	= \bigl[ 
			\bigl( \bff(\bu_{R}) \, , \bphi_{i} \bigr)
			-
			\bigl( \bff(\bu_{r}) \, , \bphi_{i} \bigr)
		\bigr] 
		\approx
		\bigl( {\boldsymbol{g}}(\bu_{r}) \, , \bphi_{i} \bigr) \, .
	\label{eqn:g}
\end{eqnarray}
To determine the coefficients/parameters in ${\boldsymbol{g}}$ used in~\eqref{eqn:g}, in the offline  stage, we solve the following low-dimensional {\it least squares problem}:
\begin{eqnarray}
	\boxed{
	 \begin{aligned}
	\min_{ {\boldsymbol{g}}\text{ parameters}} \ \sum_{j=1}^{M} \biggl\| 
	\bigl[ 
			\bigl( \bff(\bu_{R}^{FOM}(t_j)) \, , \bphi_{i} \bigr)
			&-
			\bigl( \bff(\bu_{r}^{FOM}(t_j)) \, , \bphi_{i} \bigr)
		\bigr] 
	 \\&-
	\bigl( {\boldsymbol{g}}(\bu_{r}^{FOM}(t_{j})) \, , \bphi_{i} \bigr) \biggr\|^2  ,
	 \end{aligned}
	}
	\label{eqn:least-squares-general}
\end{eqnarray}
where $\bu_{R}^{FOM}$ and $\bu_{r}^{FOM}$ are obtained from the FOM data and $M$ is the number of snapshots.
Once ${\boldsymbol{g}}$ is determined, the model~\eqref{eqn:vms-rom-two-scales-1} with the VMS-ROM closure term replaced by ${\boldsymbol{g}}$ yields the {\it two-scale data-driven VMS-ROM (2S-DD-VMS-ROM)}:
\begin{eqnarray}
	\boxed{
	\bigl( \overset{\bullet}{\bu_{r}} \, , \bphi_{i} \bigr) 
	= \bigl( \bff(\bu_{r}) \, , \bphi_{i} \bigr) 
	+ \bigl( {\boldsymbol{g}}(\bu_{r}) \, , \bphi_{i} \bigr) ,
	}
	\qquad 
	i = 1, \ldots, r.
	\label{eqn:vms-rom-two-scales-2}
\end{eqnarray}

We emphasize that, in contrast to the traditional VMS methodology, the 2S-DD-VMS-ROM framework allows {\it great flexibility} in choosing the {\it structure} of the closure term.
For example, for the NSE, the approximation~\eqref{eqn:g} becomes:
$\forall \, i = 1, \ldots, r,$
\vspace*{-0.1cm}
\begin{eqnarray}
	\text{VMS-ROM closure term} 
	&=& - \bigl[ 
			\bigl( \bigl( {\bu_R} \cdot \nabla \bigr) \, {\bu_R} \, , \bphi_{i} \bigr)	
			- \bigl( \bigl( {\bu_r} \cdot \nabla \bigr) \, {\bu_r} \, , \bphi_{i} \bigr)
		\bigr] 
	\nonumber \\[0.3cm]
	&\approx& \bigl( {\boldsymbol{g}}(\bu_{r}) \, , \bphi_{i} \bigr) 
	\nonumber \\[0.3cm]
	&=& \bigl( \tilde{A} \, \ba + \ba^\top \tilde{B} \, \ba \bigr)_{i} \ ,
	\\[-0.6cm]
	\nonumber
	\label{eqn:ansatz-nse}
\end{eqnarray}
where, for computational efficiency, we assume that the structures of ${\boldsymbol{g}}$  and $\bff$ are similar.
Thus, in the least squares  problem~\eqref{eqn:least-squares-general}, we solve for {\it  all} the entries in the $r \times r$ matrix $\tilde{A}$ and the $r \times r \times r$ tensor $\tilde{B}$:
\begin{eqnarray}
    &&
	\min_{\tA, \tB} \ \sum_{j=1}^{M} 
	\biggl\| 
	- \biggl[ 
			\bigl( \bigl( {\bu_R}^{FOM}(t_{j}) \cdot \nabla \bigr) \, {\bu_R}^{FOM}(t_{j}) \, , \bphi_{i} \bigr)	
	\nonumber \\
	&& \hspace*{2.5cm}
			- \bigl( \bigl( {\bu_r}^{FOM}(t_{j}) \cdot \nabla \bigr) \, {\bu_r}^{FOM}(t_{j}) \, , \bphi_{i} \bigr)
	\biggr]
	\nonumber \\[0.2cm]
	&& \hspace*{1.8cm}
 	- \bigl( \tilde{A} \, \ba^{FOM}(t_{j}) + \ba^{FOM}(t_{j})^\top \tilde{B} \, \ba^{FOM}(t_{j}) \bigr) 
	\biggr\|^2 \, ,
	\label{eqn:least-squares-nse}
	\\[-0.6cm]
	\nonumber
\end{eqnarray}
where $\bu_{R}^{FOM}, \bu_{r}^{FOM}$, and $\ba^{FOM}$ are obtained from the available FOM data.
Specifically, the values $\ba^{FOM}(t_j)$,
computed at snapshot time instances t$_j, j = 1,\cdots,M$, are obtained by projecting the
corresponding snapshots $\bu(t_j)$ 
onto the ROM basis functions $\bphi_i$
and using the orthogonality of the ROM basis functions: $\forall i=1,\cdots,R,\, \forall j=1,\cdots, M$, 
\begin{eqnarray}
    a_i^{FOM}(t_j) = 
    \Bigl(\bu(t_j),\bphi_i\Bigr).
    \label{eqn:a-fom}
\end{eqnarray}
In addition,
\begin{eqnarray}
    \bu_{R}^{FOM}(t_j) 
    = \sum_{k=1}^Ra_k^{FOM}(t_j)\bphi_k,
    \qquad 
    \bu_{r}^{FOM}(t_j) 
    = \sum_{k=1}^ra_k^{FOM}(t_j)\bphi_k.
    \label{eqn:u-r-fom}
\end{eqnarray}
The least squares problem~\eqref{eqn:least-squares-nse} is {\it low-dimensional} since, for a  small $r$ value, seeks the optimal $(r^2 + r^3)$ entries in $\tilde{A}$ and $\tilde{B}$, respectively.
Thus,~\eqref{eqn:least-squares-nse} can be efficiently solved in the offline stage.
For the NSE, the 2S-DD-VMS-ROM~\eqref{eqn:vms-rom-two-scales-2} takes the form
\begin{eqnarray}
	\boxed{
	\overset{\bullet}{\ba} 
	= ( A + \tilde{A}) \ba 
	+ \ba^\top (B + \tilde{B}) \ba \, ,
	}
	\label{eqn:vms-rom-two-scales-3}
\end{eqnarray}
where $A$ and $B$ are the G-ROM  operators in~\eqref{eqn:g-rom}, and $\tilde{A}$ and $\tilde{B}$ are the VMS-ROM closure operators constructed in~\eqref{eqn:least-squares-nse}.

\subsection{Three-Scale Data-Driven Variational Multiscale ROMs (3S-DD-VMS-ROM)}
    \label{sec:vms-rom-three-scales}

The 2S-DD-VMS-ROM~\eqref{eqn:vms-rom-two-scales-2} is based on the two-scale decomposition of $\bu_{R} \in \bX_1 \oplus \bX_2$ into resolved and unresolved scales:   $\bu_{R} = \bu_{r}+ \bu'$.
The flexibility of the hierarchical structure of the ROM space allows a three-scale decomposition of $\bu_{R}$, which yields a {\it three-scale data-driven VMS-ROM (3S-DD-VMS-ROM)} that is more accurate than the 2S-DD-VMS-ROM~\eqref{eqn:vms-rom-two-scales-2}.
To construct the new 3S-DD-VMS-ROM, we first build three orthogonal spaces, $\bX_1, \bX_2$, and $\bX_3$:
$\bX_1 := \text{span}  \{ \bphi_1, \ldots, \bphi_{r_1} \}$, which represents the {\it large resolved} ROM scales, $\bX_2 := \text{span}  \{ \bphi_{r_1+1}, \ldots, \bphi_{r} \}$, which represents the {\it small resolved} ROM scales, and $\bX_3 := \text{span}  \{ \bphi_{r+1}, \ldots, \bphi_{R} \}$, which represents the {\it unresolved} ROM scales.  
Next, we consider the best ROM approximation of $\bu$ in the space $\bX_1 \oplus \bX_2 \oplus \bX_3$,  i.e.,   $\bu_{R} \in \bX_1 \oplus \bX_2\oplus \bX_3$  defined as
\vspace*{-0.3cm}
\begin{eqnarray}
	\bu_{R}
	&=& \sum_{j=1}^{R} a_{j} \, \bphi_{j}
	\nonumber \\[0.3cm]
	&=& \sum_{j=1}^{r_{1}} a_{j} \, \bphi_{j}
	+ \sum_{j=r_{1}+1}^{r} a_{j} \, \bphi_{j}
	+ \sum_{j=r+1}^{R} a_{j} \, \bphi_{j} 
	\nonumber \\[0.3cm]
	&=& \bu_{L}
	+ \bu_{S}
	+ \bu' \, ,
	\label{eqn:vms-ddc-rom-1}
	\\[-0.5cm]
	\nonumber
\end{eqnarray}
where $\bu_{L} \in \bX_1$ represents the large resolved ROM component of $\bu_R$, $\bu_{S} \in \bX_2$ represents the small resolved ROM component of $\bu_R$, and $\bu' \in \bX_3$ represents the unresolved ROM component of $\bu_R$.
Thus, with the notation from Section~\ref{sec:vms-rom-two-scales}, $\bu_{r} = \bu_{L} + \bu_{S}$.
We plug $\bu_R$ in~\eqref{eqn:pde-strong}, and project the resulting equation onto both $\bX_1$ and $\bX_2$:
\vspace*{-0.1cm}
\begin{eqnarray}
\hspace*{-0.7cm}
		\biggl( \overset{\bullet}{\bu_{L}} \, , \bphi_{i} \biggr) 
		&=& 
		\biggl( \bff(\bu_{L} + \bu_{S}) \, , \bphi_{i} \biggr) 
		+
		\boxed{
		\biggl[ 
			\biggl( \bff(\bu_{R}) \, , \bphi_{i} \biggr)
			-
			\biggl( \bff(\bu_{L} + \bu_{S}) \, , \bphi_{i} \biggr)
		\biggr] 
		} \ , 
		\nonumber \\[0.5cm] 
		&& \hspace*{7.0cm}
		\forall \, i = 1, \ldots, r_{1} ,
		\label{eqn:vms-rom-three-scales-1} \\[0.5cm]
		\biggl( \overset{\bullet}{\bu_{S}} \, , \bphi_{i} \biggr) 
		&=& 
		\biggl( \bff(\bu_{L} + \bu_{S}) \, , \bphi_{i} \biggr) 
		+
		\boxed{
		\biggl[ 
			\biggl( \bff(\bu_{R}) \, , \bphi_{i} \biggr)
			-
			\biggl( \bff(\bu_{L} + \bu_{S}) \, , \bphi_{i} \biggr)
		\biggr] 
		} \ , 
		\nonumber \\[0.5cm] 
		&& \hspace*{7.0cm}
		\forall \, i = r_{1}+1, \ldots, r .
		\label{eqn:vms-rom-three-scales-2}
\end{eqnarray}
The two boxed terms in~\eqref{eqn:vms-rom-three-scales-1}--\eqref{eqn:vms-rom-three-scales-2} are the {\it VMS-ROM closure terms},  which have {\it fundamentally different roles}: 
The VMS-ROM closure term in~\eqref{eqn:vms-rom-three-scales-1} models the interaction between the large resolved ROM modes and the small resolved ROM modes;
the VMS-ROM closure term in~\eqref{eqn:vms-rom-three-scales-2} models the interaction between the small resolved ROM modes and the unresolved ROM modes.
The new 3S-DD-VMS-ROM framework allows {\it great flexibility} in choosing the {\it structure} of the two VMS-ROM closure terms.
For the NSE, we can use the following approximations:
\vspace*{-0.1cm}
\begin{eqnarray}
    (\btau_L)_i
    &:=&
	- \bigl[ 
			\bigl( \bigl( {\bu_R} \cdot \nabla \bigr) \, {\bu_R} \, , \bphi_{i} \bigr)	
			- \bigl( \bigl( (\bu_{L} + \bu_{S}) \cdot \nabla \bigr) \, (\bu_{L} + \bu_{S}) \, , \bphi_{i} \bigr)
		\bigr] 
    \nonumber \\[0.2cm]
	&=& \bigl( \tilde{A}_L \, \ba + \ba^\top \tilde{B}_L \, \ba \bigr)_{i} , 
	\qquad
	\forall \, i = 1, \ldots, r_1 \ ,
	\label{eqn:vms-rom-three-scales-3}  \\[0.3cm]
    (\btau_S)_i
    &:=&
	- \bigl[ 
			\bigl( \bigl( {\bu_R} \cdot \nabla \bigr) \, {\bu_R} \, , \bphi_{i} \bigr)	
			- \bigl( \bigl( (\bu_{L} + \bu_{S}) \cdot \nabla \bigr) \, (\bu_{L} + \bu_{S}) \, , \bphi_{i} \bigr)
		\bigr] 
	    \nonumber \\[0.2cm]
	&=& \bigl( \tilde{A}_S \, \ba + \ba^\top \tilde{B}_S \, \ba \bigr)_{i} 
	\qquad
	\forall \, i = r_1+1, \ldots, r  ,
	\label{eqn:vms-rom-three-scales-4}  
\end{eqnarray}
where $\tilde{A}_L \in \R^{r_1 \times r}, \tilde{A}_S \in \R^{(r-r_1) \times r}, \tilde{B}_L \in \R^{r_1 \times r \times r},$ and $\tilde{B}_S \in \R^{(r-r_1) \times r \times r}$.
To determine the entries in $\tilde{A}_L,  \tilde{A}_S, \tilde{B}_L,$ and $\tilde{B}_S$, we solve two least squares problems:
\vspace*{-0.2cm}
\begin{eqnarray}
    &&
	\hspace*{-0.7cm}
	\min_{\tA_L, \tB_L} \ \sum_{j=1}^{M} 
	\bigl\| 
	\btau_L^{FOM}
 	- \bigl( \tilde{A}_L \, \ba^{FOM}(t_{j}) + \ba^{FOM}(t_{j})^\top \tilde{B}_L \, \ba^{FOM}(t_{j}) \bigr) 
	\bigr\|^2 \, ,
	\label{eqn:vms-rom-three-scales-5}  \\
   &&
	\hspace*{-0.7cm}
	\min_{\tA_S, \tB_S} \ \sum_{j=1}^{M} 
	\bigl\| 
	\btau_S^{FOM}
 	- \bigl( \tilde{A}_S \, \ba^{FOM}(t_{j}) + \ba^{FOM}(t_{j})^\top \tilde{B}_S \, \ba^{FOM}(t_{j}) \bigr) 
	\bigr\|^2 \, ,
	\label{eqn:vms-rom-three-scales-6}
	\\[-0.4cm]
	\nonumber
\end{eqnarray}
where $\btau_L^{FOM}, \btau_S^{FOM}$, and $\ba^{FOM}$ are obtained from the available FOM data.

For the NSE, the {\it three-scale data-driven VMS-ROM (3S-DD-VMS-ROM)} is
\vspace*{-0.1cm}
\begin{eqnarray}
	\left[
		\begin{array}{c}
			\overset{\bullet}{\ba_{L}} \\[0.3cm]
			\overset{\bullet}{\ba_{S}} 
		\end{array}
	\right]
	&=&
	A \, \ba
	+ 
	\ba^{\top} \, B \, \ba
	+ 
	\left[
		\begin{array}{c}
			\tA_{L} \, \ba
			+
			\ba^{\top} \, \tB_{L} \, \ba \\[0.3cm]
			\tA_{S} \, \ba
			+
			\ba^{\top} \, \tB_{S} \, \ba 
		\end{array}
	\right] \, ,
	\qquad
	\label{eqn:vms-rom-three-scales-7} 	
	\\[-0.4cm]
	\nonumber
\end{eqnarray}
where
$\ba^{\top} = [\ba_{L} , \ba_{S}]^{\top}$, $A$ and $B$ are the G-ROM  operators in~\eqref{eqn:g-rom}, and $\tilde{A}_L,  \tilde{A}_S, \tilde{B}_L,$ and $\tilde{B}_S$ are the VMS-ROM closure operators constructed in~\eqref{eqn:vms-rom-three-scales-5}--\eqref{eqn:vms-rom-three-scales-6}. 
Compared to the 2S-DD-VMS-ROM, in the 3S-DD-VMS-ROM we have {\it more flexibility} in choosing the VMS-ROM closure operators $\tilde{A}_L,  \tilde{A}_S, \tilde{B}_L,$ and $\tilde{B}_S$ in the least squares problems~\eqref{eqn:vms-rom-three-scales-5}--\eqref{eqn:vms-rom-three-scales-6}.  
For example,  for $\tilde{A}_L, \tilde{B}_L$ we can specify physical constraints, sparsity patterns, or regularization parameters, that are different from those for $\tilde{A}_S, \tilde{B}_S$.  
Because of this increased flexibility, we expect that the 3S-DD-VMS-ROM~\eqref{eqn:vms-rom-three-scales-7} is {\it more accurate} than the 2S-DD-VMS-ROM~\eqref{eqn:vms-rom-two-scales-3}.

\section{Numerical Results}
	\label{sec:numerical-results}


In this section, we perform a numerical investigation of the new DD-VMS-ROM framework.
As noted in Section~\ref{sec:vms-rom}, the 2S-DD-VMS-ROM~\eqref{eqn:vms-rom-two-scales-3} was investigated in~\cite{xie2018data} under the name ``data-driven filtered ROM'' and in~\cite{mohebujjaman2019physically,mou2019data} under the name ``data-driven correction ROM.''
In~\cite{xie2018data}, it was shown that the 2S-DD-VMS-ROM is more accurate than the standard G-ROM in the numerical simulation of 2D flow past a circular cylinder at Reynolds numbers $Re=100,  Re=500$, and $Re=1000$.
Furthermore, the 2S-DD-VMS-ROM was more accurate and more efficient than other modern ROM closure models.
In~\cite{mou2019data}, it was shown that the 2S-DD-VMS-ROM is more accurate than the standard G-ROM in the numerical simulation of the quasi-geostrophic equations modeling the large scale ocean circulation.

Since the 2S-DD-VMS-ROM has already been shown to perform well, the focus of the current numerical investigation is on the new 3S-DD-VMS-ROM~\eqref{eqn:vms-rom-three-scales-7}.
Specifically, we investigate whether the 3S-DD-VMS-ROM is more accurate than the 2S-DDC-ROM.
To this end, we consider {four} test cases:
(i) the 1D viscous Burgers equation with viscosity coefficient $\nu=10^{-3}$ (Section~\ref{sec:numerical-results-burgers});  
(ii) a 2D flow past a circular cylinder at Reynolds numbers $Re=100$, $Re=500$, and $Re=1000$ (Section~\ref{sec:numerical-results-nse});
{
(iii) the quasi-geostrophic equations at Reynolds number $Re=450$ and Rossby number $Ro=0.0036$ (Section~\ref{sec:qge}); and 
(iv) a 2D flow over a backward facing step at Reynolds number $Re=1000$ (Section~\ref{sec:bfs}).
}
For each test case, we investigate three ROMs:
the 2S-DD-VMS-ROM~\eqref{eqn:vms-rom-two-scales-3}, the new 3S-DD-VMS-ROM~\eqref{eqn:vms-rom-three-scales-7}, and (for comparison purposes) the standard G-ROM~\eqref{eqn:g-rom}.
As a benchmark, we use the FOM results.

We test 
the ROMs in three different regimes:

(i) {\it Reconstructive} regime: 
The ROM basis and ROM operators $A$ and $B$ are constructed from FOM data obtained on the time interval $[0,T_{1}]$, and then the resulting ROMs are tested on the {\it  same} time interval $[0,T_{1}]$.
To construct the DD-VMS-ROM operators $\tA$  and $\tB$ (for the 2S-DD-VMS-ROM)  and  $\tA_{L},  \tA_{S}, \tB_{L},$ and $\tB_{S}$ (for the 3S-DD-VMS-ROM), we use different approaches for the four test cases: 
For the Burgers equation, {quasi-geostrophic equations, and backward facing step test cases}, we construct the DD-VMS-ROM operators by using FOM data from the entire time interval $[0,T_{1}]$. 
For the flow past a circular cylinder test case, for computational efficiency, we construct the DD-VMS-ROM operators from FOM data obtained on a  shorter time interval, which does not significantly decrease the accuracy of the resulting DD-VMS-ROM.
Specifically, we use FOM data for one period~\cite{xie2018data,mohebujjaman2019physically}, i.e., (i) from $t = 7$ to $t = 7.332$ for $Re = 100$, (ii) from $t = 7$ to $t = 7.442$ for $Re = 500$, and (iii) from $t = 13$ to $t = 13.268$ for $Re = 1000$.

(ii) {\it Cross-validation} regime: 
The ROM  basis and ROM operators $A$ and $B$ are constructed from FOM data obtained on the time interval $[0, T_{2}]$, and then the resulting ROMs are tested on the time interval $[0, T_{3}]$, where $T_{3} > T_{2}$.
We note that the two time intervals are different, but they do overlap over $[0, T_{2}]$.
To construct the DD-VMS-ROM operators, we use different approaches for the two test cases: 
For the Burgers equation test case, we construct the DD-VMS-ROM operators by using FOM data from the entire time interval $[0, T_{2}]$. 
For the flow past a circular cylinder test case, for computational efficiency, we construct the DD-VMS-ROM operators from FOM data for one period~\cite{xie2018data,mohebujjaman2019physically}.




(iii) {\it Predictive} regime: 
The ROM  basis and ROM operators $A$ and $B$  are constructed from FOM data obtained on the time interval $[0, T_{2}]$, and then the resulting ROMs are tested on the time interval $[T_{2}, T_{3}]$, where $T_{3} > T_{2}$.
We emphasize that the two time intervals are completely different, without any overlap.
To construct the DD-VMS-ROM operators, we use different approaches for the two test cases: 
For the Burgers equation test case, we construct the DD-VMS-ROM operators by using FOM data from the entire time interval $[0, T_{2}]$. 
For the flow past a circular cylinder test case, for computational efficiency, we construct the DD-VMS-ROM operators from FOM data for half a  period~\cite{xie2018data,mohebujjaman2019physically}.

\subsection{Computational Setting}
	\label{sec:numerical-results-setting}

In this section, we present the computational setting used in the numerical investigation.

First, as explained in detail on page B843 of~\cite{xie2018data}, we rewrite the optimization problem~\eqref{eqn:least-squares-nse} as the least squares problem
\begin{equation}
		\min_{\substack{\bx \in \R^{(r^2 + r^3) \times 1}}} \, 
		\left\| \, \bff - E \, \bx \, \right\|^2 \, ,
	\label{eqn:least-squares}
\end{equation}
where $\bx \in \R^{(r^2 + r^3) \times 1}$ contains all the entries of $\tA$ and $\tB$, and the vector $\bff \in \R^{(M \, r) \times 1}$ and matrix $E \in \R^{(M \, r) \times (r^2 + r^3)}$ are computed from $\bu_{R}^{FOM}, \bu_{r}^{FOM}$, and $\ba^{FOM}$ (see (4.8) in~\cite{xie2018data}).
The optimal $\tA$ and $\tB$ (i.e., the entries in $\bx$ that solves the linear least squares problem~\eqref{eqn:least-squares}) are used to build the 2S-DD-VMS-ROM~\eqref{eqn:vms-rom-two-scales-3}.

Furthermore, as explained on page B843 of~\cite{xie2018data}, the least squares problem~\eqref{eqn:least-squares} is ill-conditioned.
This ill-conditioning is common in data-driven least squares problems (see, e.g.,~\cite{peherstorfer2016data}).
To alleviate this ill-conditioning, we use the truncated singular value decomposition (SVD)~\cite{mohebujjaman2019physically,xie2018data}.

The algorithm for the 2S-DD-VMS-ROM~\eqref{eqn:vms-rom-two-scales-3} is presented in Algorithm~\ref{alg:2s-dd-vms-rom}.
In most of our numerical experiments, we choose the optimal tolerance $tol$ in the truncated SVD step of   Algorithm~\ref{alg:2s-dd-vms-rom}.
Specifically, for each value $1 \leq  m \leq  R$ (where $R$ is the dimension of the snapshot matrix), we consider the truncated SVD approximation of dimension $m$, construct the operators $\tA_{m}$ and $\tB_{m}$, integrate the resulting 2S-DD-VMS-ROM in~\eqref{eqn:alg-ddf-rom-5}, and choose the $\widetilde{m}$ value yielding the lowest $L^{2}$ error.
The only exception is in some of the numerical experiments for the Burgers equation (Section~\ref{sec:burgers}), where we fix $tol=tol_L$ or $tol = tol_S$ (see  Tables~\ref{table:burgers-reconstructive-1a}--\ref{table:burgers-reconstructive-3b}).

\begin{algorithm}[H]
	\caption{2S-DD-VMS-ROM}
	\label{alg:2s-dd-vms-rom}
	\begin{algorithmic}[1]
		\STATE{
					Use all the entries of $\tA$ and $\tB$ in~\eqref{eqn:vms-rom-two-scales-3} to define vector of unknowns, $\bx$.
					}
		\STATE{
					Use $\bu_{R}^{FOM}, \bu_{r}^{FOM}$, and $\ba^{FOM}$ to assemble the vector $\bff$ and matrix $E$ in~\eqref{eqn:least-squares}.
					}
		\STATE{
					Use the truncated SVD algorithm to solve the linear least squares problem~\eqref{eqn:least-squares}.
					}
  			\begin{enumerate}\itemsep0em 
				\item[(i)] Calculate the SVD of $E$: 
					\begin{eqnarray}
						E = U \, \Sigma V^{\top} \, ,
						\label{eqn:svd}
					\end{eqnarray}
					where the rank of matrix $E$($\Sigma$) is $\mathcal{M}$.
				\item[(ii)] Specify tolerance $tol = \sigma_i\, , i=1,\cdots,\mathcal{M}$.
				\item[(iii)] Construct matrix $\widehat{\Sigma}^m$ from $\Sigma$ as follows: 
							$\widehat{\sigma}_m = \sigma_m$ if $\sigma_m \ge tol$, $m=1,\cdots,\mathcal{M}$.
				\item[(iv)] Construct $\widehat{E}^m$, the truncated SVD of $E$:
					\begin{eqnarray}
						\widehat{E}^m = \widehat{U}^m \, \widehat{\Sigma}^m \, \bigl(\widehat{V}^m\bigr)^{\top} \, ,
						\label{eqn:truncated-svd}
					\end{eqnarray}
					where $\widehat{U}^m$ and $\widehat{V}^m$ are the entries of $U$ and $V$ in~\eqref{eqn:svd} that correspond to $\widehat{\Sigma}^m$.
				\item[(v)] The solution of the least squares problem~\eqref{eqn:least-squares} is 
					\begin{eqnarray}
						\bx
						= \left( \widehat{V}^m \, \bigl(\widehat{\Sigma}^m\bigr)^{-1} \, \bigl(\widehat{U}^m\bigr)^{\top} \right) \, \bff \, .
					\label{eqn:least-squares-solution}
				\end{eqnarray}
			\end{enumerate}
			\vspace{-1pc}
		\STATE{
			The 2S-DD-VMS-ROM~\eqref{eqn:vms-rom-two-scales-3} has the following form:
								\begin{eqnarray}
									\boxed{
									\overset{\bullet}{\ba} 
									= \left( A + \tA^m \right) \, \ba
									+ \ba^{\top} \, \left( B + \tB^m \right) \, \ba
									}\,,  
									\label{eqn:alg-ddf-rom-5}
								\end{eqnarray}
								where $\tA^m$ and $\tB^m$ are the appropriate entries of $\bx$ found in~\eqref{eqn:least-squares-solution} with $tol=\sigma_m$.
					}
			\STATE{
			Integrate the resulting 2S-DD-VMS-ROM in~\eqref{eqn:vms-rom-two-scales-3} over the given time domain and calculate the average $L^{2}$ error $\mathcal{E}^{m}(L^2)$ by using formula~\eqref{eqn:l2-error}.
           The optimal $\widetilde{m}$ value (the optimal operators $\tA$ and $\tB$) is found by solving the following minimization problem:
\begin{eqnarray}
\mathcal{E}^{\widetilde{m}}(L^2) 
=
\min_{1 \le m \le R} \mathcal{E}^m(L^2)\,.
	\label{numerical-results-nse-1}
\end{eqnarray}
}					
\end{algorithmic}
\end{algorithm}

\bigskip

The algorithm for the 3S-DD-VMS-ROM~\eqref{eqn:vms-rom-three-scales-7} is the same as Algorithm~\ref{alg:2s-dd-vms-rom}, except that we are using two different truncated SVDs to solve two different linear least squares problems, which correspond to the large and small resolved scales. 
Thus, we have two different control parameters, $tol_L$ and $tol_S$. 
Similar to the 2S-DD-VMS-ROM , we rewrite the optimization problems~\eqref{eqn:vms-rom-three-scales-5} and~\eqref{eqn:vms-rom-three-scales-6} as the least squares problems
\begin{equation}
		\min_{\substack{\bx_L \in \R^{\left[r_1(r+r^2) \right] \times 1}}\quad} \, 
		\left\| \, \bff_L - E_L \, \bx_L \, \right\|^2 \, ,
	\label{eqn:least-squares-3S-1}
\end{equation}
\begin{equation}
		\min_{\substack{\bx \in \R^{\left[(r-r_1)(r+r^2) \right] \times 1}}} \, 
		\left\| \, \bff_S - E_S \, \bx_S \, \right\|^2 \, ,
	\label{eqn:least-squares-3S-2}
\end{equation}
where $\bx_L \in \R^{\left[r_1(r+r^2) \right] \times 1}$ contains all the entries of the operators $\tA_L$ and $\tB_L$, $\bx_S \in \R^{\left[(r-r_1)(r+r^2) \right] \times 1}$ contains all the entries of the operators $\tA_S$ and $\tB_S$, and the vectors $\bff_L \in \R^{(M \, r_1) \times 1}$, $\bff_S \in \R^{(M \,(r- r_1)) \times 1}$ and the matrices $E_L \in \R^{(M \, r_1) \times (r_1(r+r^2))}$, $E_S \in \R^{(M \, (r-r_1)) \times ((r-r_1)(r+r^2))}$ are computed from $\bu_{R}^{FOM}, \bu_{r}^{FOM}$, and $\ba^{FOM}$ (see (4.8) in~\cite{xie2018data}).
The optimal $\tA_L$, $\tB_L$ and $\tA_S$, $\tB_S$ (i.e., the entries in $\bx_L$ and $\bx_S$ that solve the linear least squares problems~\eqref{eqn:least-squares-3S-1} and ~\eqref{eqn:least-squares-3S-2}) are used to build the 3S-DD-VMS-ROM~\eqref{eqn:vms-rom-three-scales-7}.
Again, to address the  ill-conditioning of the least squares problems~\eqref{eqn:least-squares-3S-1}--\eqref{eqn:least-squares-3S-2}, we use the truncated SVD algorithm. 
The algorithm for the 3S-DD-VMS-ROM~\eqref{eqn:vms-rom-three-scales-7} is presented in  Algorithm~\ref{alg:3s-dd-vms-rom}.
We note that, if $tol_L$ = $tol_S$ = $tol$, the 2S-DD-VMS-ROM and 3S-DD-VMS-ROM yield the same results, since we solve the same minimization problem. Thus, the interesting case is when $tol_L$ and/or $tol_S$ are different from $tol$.
\vspace{0.05in}

\begin{algorithm}[H]
	\caption{3S-DD-VMS-ROM
	}
	\label{alg:3s-dd-vms-rom}
	\begin{algorithmic}[1]
		\STATE{
					Choose $r_1\, ,1\le r_1<r$, and use all the entries of $\tA_L$ and $\tB_L$ as well as $\tA_S$ and $\tB_S$ in~\eqref{eqn:vms-rom-three-scales-7} to define vectors of unknowns, $\bx_L$ and $\bx_S$, respectively.
					}
		\STATE{
					Use $\bu_{R}^{FOM}, \bu_{r}^{FOM}$, and $\ba^{FOM}$ to assemble the vectors $\bff_L$ and $\bff_S$, and the  matrices $E_L$ and $E_S$ in~\eqref{eqn:least-squares-3S-1} and~\eqref{eqn:least-squares-3S-2}. 
					}
		\STATE{
					Use the truncated SVD algorithm to solve the linear least squares problems~\eqref{eqn:least-squares-3S-1} and~\eqref{eqn:least-squares-3S-2}.
					}
  			\begin{enumerate}
				\item[(i)] Calculate the SVD of $E_L$ and  $E_S$: 
					\begin{eqnarray}
						E_L = U_L \, \Sigma_L V_L^{\top} \,,\quad E_S = U_S \, \Sigma_S V_S^{\top} \, .
						\label{eqn:svd-2}
					\end{eqnarray}
				\item[(ii)] Specify tolerances $tol_L = \sigma_{L,i}\, , i=1,\cdots,\mathcal{M}_L$,  and $tol_S = \sigma_{S,j}\, , j=1,\cdots,\mathcal{M}_S$, where $\mathcal{M}_L$ is the rank of $\Sigma_L$ ($E_L$) and $\mathcal{M}_S$ is the rank of $\Sigma_S$ ($E_S$).
				\item[(iii)] Construct matrix $\widehat{\Sigma}_L^m$ from $\Sigma_L$ as follows: 
							$\widehat{\sigma}_{L,m_L} = \sigma_{m_L}$ if $\sigma_{m_L}\ge tol_L$, $m_L=1,\cdots,\mathcal{M}_L$; construct matrix $\widehat{\Sigma}_L^m$ from $\Sigma_L$ as follows: 
							$\widehat{\sigma}_{S,m_L} = \sigma_{m_S}$ if $\sigma_{m_S}\ge tol_S$, $m_S=1,\cdots,\mathcal{M}_S$.
				\item[(iv)] Construct $\widehat{E}_L^{m_L}$ and $\widehat{E}_S^{m_S}$ with the truncated SVD of $E_L$ and $E_S$.
				\item[(v)] Construct the operators $\tA_{L}^{m_L}$ and $\tB_{L}^{m_L}$ as well as $\tA_{S}^{m_S}$ and $\tB_{S}^{m_S}$.
				\item[(vi)] Integrate the resulting 3S-DD-VMS-ROM in~\eqref{eqn:vms-rom-three-scales-7} over the given time domain and calculate the average $L^{2}$ error $\mathcal{E}^{r_1,m_L,m_S}(L^2)$ by using formula~\eqref{eqn:l2-error}.
			\end{enumerate}
		\STATE{Find the optimal $\widetilde{r_1}$, $\widetilde{m}_L$ and $\widetilde{m}_S$ values (i.e., the optimal operators $\tA_L,\tB_L$ and $\tA_S,\tB_S$ corresponding to the optimal $r_1$) by solving the following minimization problem:
		\begin{align}
					\mathcal{E}^{\widetilde{r_1},\widetilde{m}_L,\widetilde{m}_S}(L^2) 
=
\min_{
		\substack{ 1\le r_1 < r\\[0.1cm] 
		                 1\le m_L\le \mathcal{M}_L \\[0.1cm] 
				1\le m_S\le \mathcal{M}_S
				}
		} \,  \mathcal{E}^{r_1,m_L,m_S}(L^2)\,.
	\label{numerical-results-nse-2}
\end{align}
}
\end{algorithmic}
\end{algorithm}


\bigskip

The focus of the current numerical investigation is on the numerical accuracy of the new DD-VMS-ROMs.
Thus, we use all the available data to build the DD-VMS-ROM operators.
We emphasize, however, that the computational cost of the construction of the DD-VMS-ROM operators can be significantly decreased by using the approach proposed on page B848 in~\cite{xie2018data}.

To compare the ROMs' performance, in the {Burgers equation, flow past a circular cylinder, and backward facing step test cases}, we use the  error metric
\begin{align}
&\text{average $L^2$ norm:} \  \ \mathcal{E} (L^2)= \frac{1}{M}\sum_{j=1}^M\, \left\| \bu_r(t_j)-\sum_{i=1}^r\left(\bu^{\text{FOM}}(t_j),\varphi_i\right)\varphi_i \right\|_{L^2} \, ,
\label{eqn:l2-error}
\end{align}
{ whereas in the quasi-geostrophic equations test case we use the error metric~\eqref{Eq_err_streamfunc}.}
In the flow past a circular cylinder, {quasi-geostrophic equations, and flow over a backward facing step test cases}, we  plot the time evolution of the ROM kinetic energy.
{
Furthermore, in the quasi-geostrophic equations test case, we use the $L^2$ error of the time-averaged streamfunction, and plot the time-averaged streamfunction.
Finally, in the backward facing step test case, we  plot the time evolution of the $y$-component of the velocity, and the spectrum of the $y$-component of the velocity at a control point. 
}


\subsection{Burgers Equation}
	\label{sec:numerical-results-burgers}

In this section, we investigate the 2S-DD-VMS-ROM~\eqref{eqn:vms-rom-two-scales-3} and the new 3S-DD-VMS-ROM~\eqref{eqn:vms-rom-three-scales-7} in the numerical simulation of the one-dimensional viscous Burgers equation: 
\begin{equation}
\begin{cases}
\displaystyle~~	u_t -\nu u_{xx} + u u_x = 0~,~~~x \in [0,1],~t\in[0,1], \\
\displaystyle~~	u(0,t)= u(1,t) = 0~,~~~t \in [0,1],  
\end{cases}
\label{eqn:burgers}
\end{equation}
with the initial condition 
\begin{equation}
u_0(x)=\begin{cases}
\displaystyle~ 1, & x \in (0,1/2],\\
~\displaystyle 0, & x \in (1/2,1], 
\end{cases}
	\label{eqn:initial-condition}
\end{equation}
and $\nu=10^{-3}$.
This test problem has been used in~\cite{ahmed2018stabilized,iliescu2014are,KV01,xie2018data}.

\paragraph{Snapshot Generation}
We generate the FOM results by using a linear FE spatial discretization with mesh size $h=1/2048$ and a Crank-Nicolson time discretization with timestep size $\Delta t=10^{-3}$. 

\paragraph{ROM Construction}
We run the FOM from $t = 0$ to $t = 1$. 
To generate the ROM basis functions, we collect a total of 1000 snapshots for the reconstructive regime, and 700 snapshots for the cross-validation and predictive regimes. 
These snapshots are the  solutions from $t = 0$ to $t = 1$ for the reconstructive regime, and $t = 0$ to $t = 0.7$ for the cross-validation and predictive regimes. 
To train $\tilde{A}$, $\tilde{B}$ (for the 2S-DD-VMS-ROM) and $\tilde{A_L}$, $\tilde{B_L}$ and $\tilde{A_S}$, $\tilde{B_S}$ (for
the 3S-DD-VMS-ROM), we use FOM data on the time interval $[0, 1]$ for the reconstructive regime, and FOM data on the time interval $[0, 0.7]$ for the cross-validation and predictive regimes. 
We test all the ROMs on the time interval $[0, 1]$ for the  reconstructive and cross-validation regimes, and $[0.7, 1]$ for the predictive regime.


\color{black}

\paragraph{Implementation Details}

To implement the 2S-DD-VMS-ROM~\eqref{eqn:alg-ddf-rom-5}, we use Algorithm~\ref{alg:2s-dd-vms-rom}.
To implement the 3S-DD-VMS-ROM~\eqref{eqn:vms-rom-three-scales-7}, we use Algorithm~\ref{alg:3s-dd-vms-rom}.
For a fair comparison of the 2S-DD-VMS-ROM with the 3S-DD-VMS-ROM, we choose optimal tolerances in the two algorithms, i.e., optimal $tol$ in Algorithm~\ref{alg:2s-dd-vms-rom} and optimal $tol_L$ and $tol_s$ in Algorithm~\ref{alg:3s-dd-vms-rom}.
We also investigate whether there is any relationship between the 2S-DD-VMS-ROM tolerance and the 3S-DD-VMS-ROM tolerances.
To this end, we perform two sets of numerical  experiments:
(a) In the first set of experiments, we fix $tol$ in the 2S-DD-VMS-ROM, choose $tol_L = tol$ in the 3S-DD-VMS-ROM, and search the optimal $tol_S$ in the 3S-DD-VMS-ROM. 
(b) In the second set of experiments, we fix $tol$ in the 2S-DD-VMS-ROM, choose $tol_S = tol$ in the 3S-DD-VMS-ROM, and search the optimal $tol_L$ in the 3S-DD-VMS-ROM.



\subsubsection{Numerical Results} 
	\label{sec:burgers}
	
In this section, we present numerical results for the Burgers equation~\eqref{eqn:burgers} with $\nu=10^{-3}$ in the reconstructive, cross-validation, and predictive regimes. 
In all the tables, we list the average $L^{2}$ error~\eqref{eqn:l2-error} for the G-ROM, the 2S-DD-VMS-ROM, and the new 3S-DD-VMS-ROM.
We also list the tolerances used in the truncated SVD algorithm for the 2S-DD-VMS-ROM  and 3S-DD-VMS-ROM, as well as the $r_1$ values for the 3S-DD-VMS-ROM.

In Table~\ref{table:burgers-reconstructive-1c}, we  list the ROMs errors for the reconstructive regime with optimal $tol$ in Algorithm~\ref{alg:2s-dd-vms-rom} and optimal $tol_L$ and $tol_s$ in Algorithm~\ref{alg:3s-dd-vms-rom}.
These results show that, for all $r$ values, the 2S-DD-VMS-ROM  and 3S-DD-VMS-ROM are several times (sometimes one or even two orders of magnitude) more accurate than the standard G-ROM. 
Overall, the 3S-DD-VMS-ROM is more accurate than the 2S-DD-VMS-ROM.
For example, for $r=7$, the 3S-DD-VMS-ROM is more than twice more accurate than the 2S-DD-VMS-ROM.
We also note that, for low $r$  values,  the ROM errors do not seem to converge monotonically. 
We emphasize, however, that for large $r$ values, we recover the expected asymptotic convergence.

\begin{figure}[h!]
	\begin{center}
		\includegraphics[width=0.4\textwidth,height=0.3\textwidth]{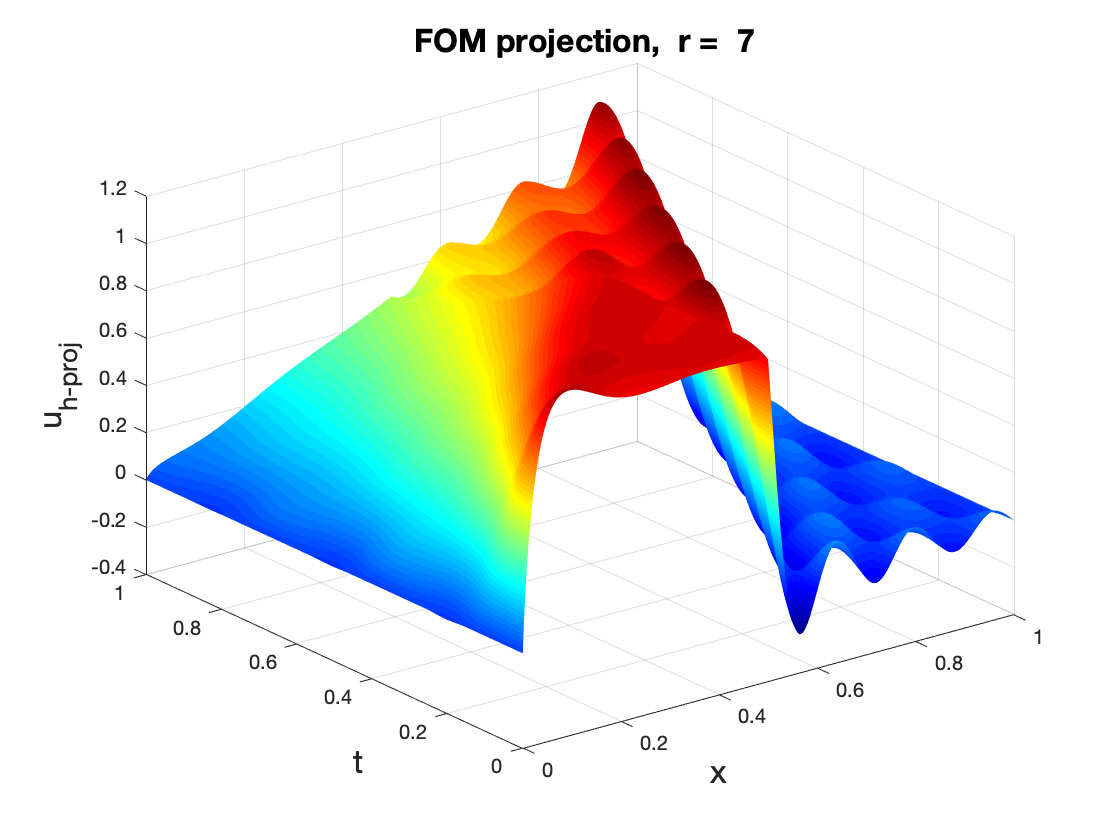}		
		\includegraphics[width=0.4\textwidth,height=0.3\textwidth]{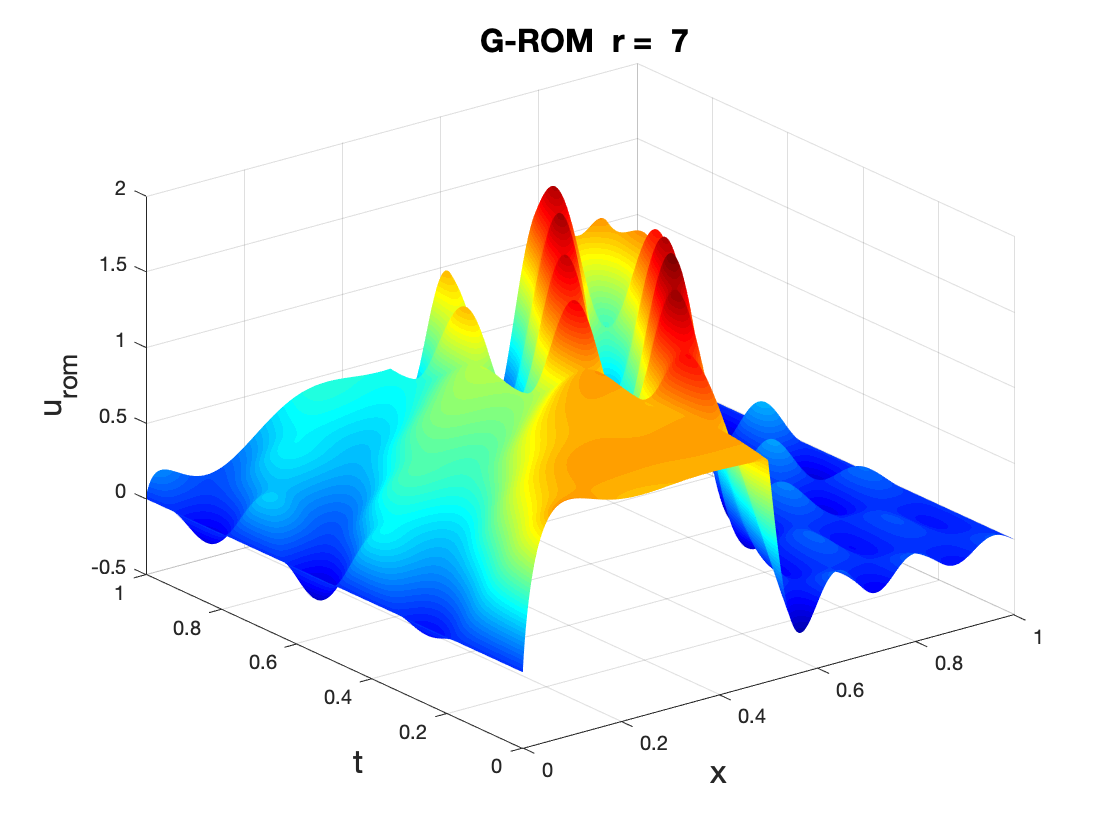}
		\includegraphics[width=0.4\textwidth,height=0.3\textwidth]{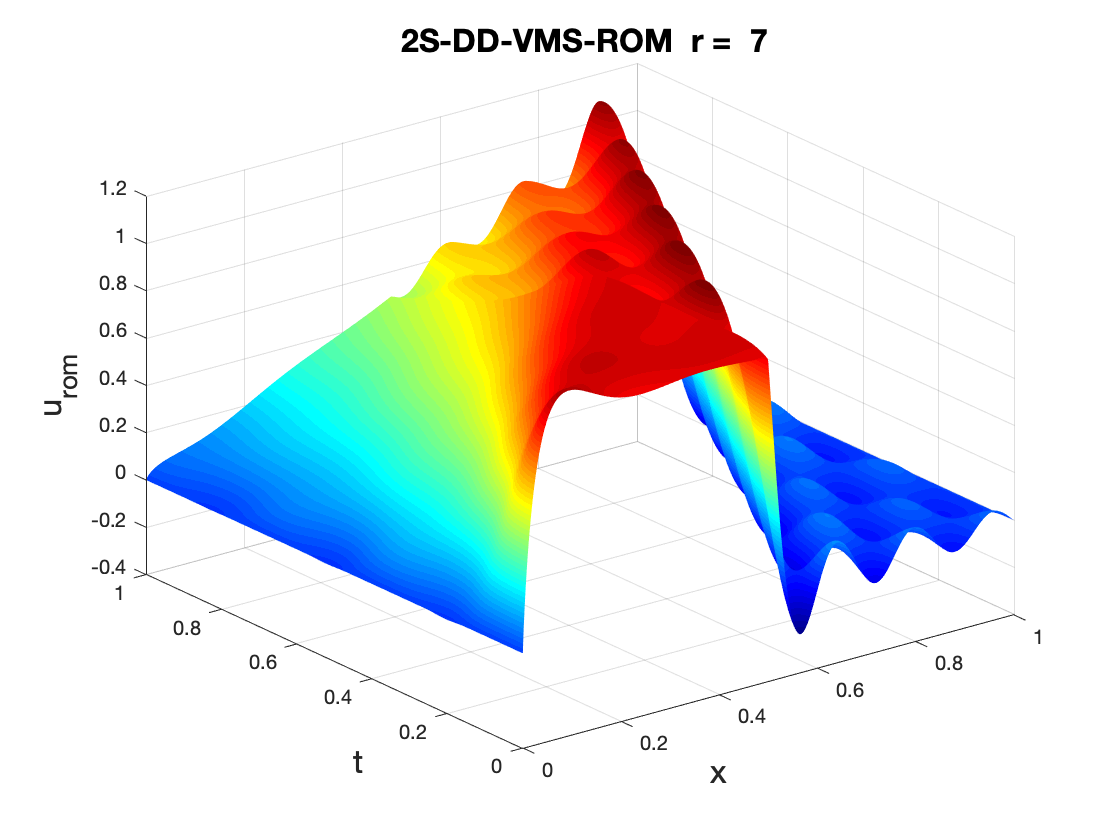}
		\includegraphics[width=0.4\textwidth,height=0.3\textwidth]{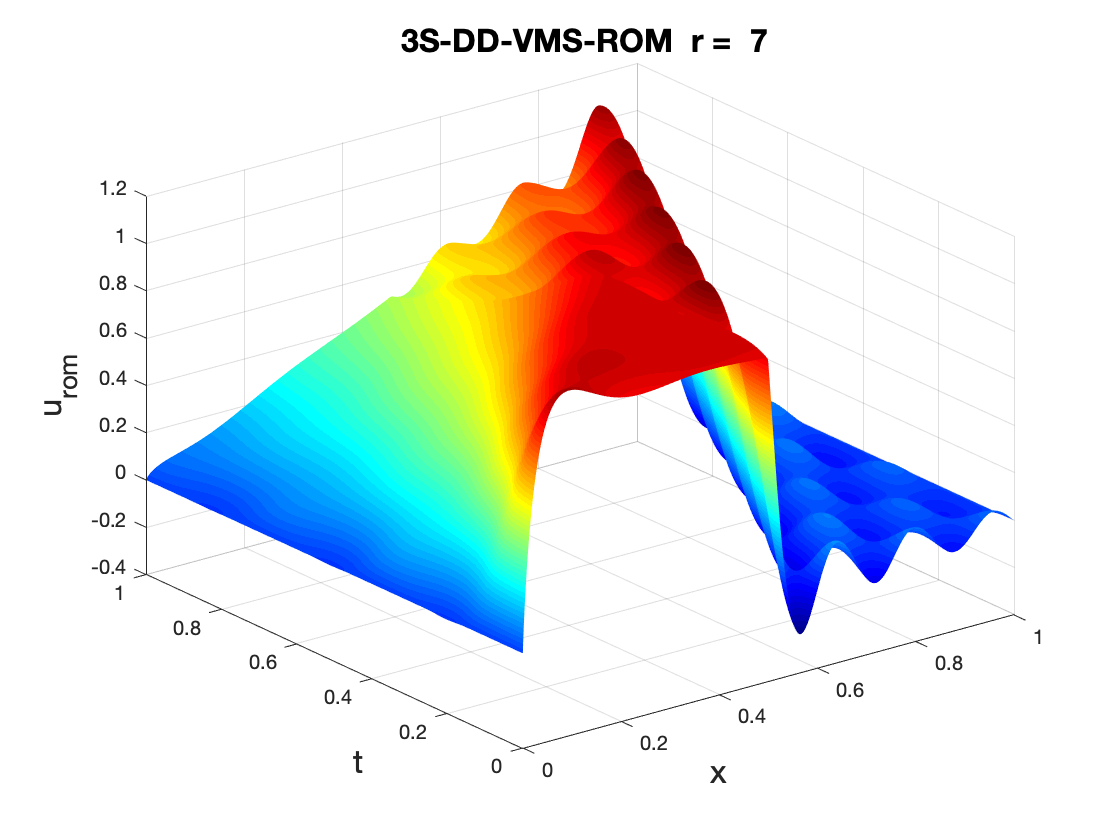}		
	\end{center}
	\caption {
			Burgers equation, $\nu = 10^{-3}$, reconstructive regime.
			FOM projection, G-ROM, 2S-DD-VMS-DDC-ROM, and 3S-DD-VMS-DDC-ROM plots for $r = 7$. 
	\label{fig:burgers-reconstructive}
	} 
\end{figure}

\begin{table}[h!]
\small 
	\begin{center}
		\smallskip
		\begin{tabular}{|c|c|c|c|c|c|c|c|c|}
			\hline
			\multicolumn{1}{|c|}{ $r$} &  
			\multicolumn{1}{c|}{ G-ROM } &  
			\multicolumn{2}{c|}{ 2S-DD-VMS-ROM } &  
			\multicolumn{4}{c|}{ 3S-DD-VMS-ROM}
						 \\  \cline{2-8}&{$\mathcal{E}(L^2)$ }   & {$tol$}& {$\mathcal{E}(L^2)$ }  &{$r_1$}& {$tol_S$}& {$tol_L$}&{ $\mathcal{E}(L^2)$ } 
									 \\
			 \hline  
$3$
 &  1.181e-01 
& 1e-02 & 1.548e-03
& 1 &  1e-02 & 1e-02 &  1.548e-03
 	 \\
			 \hline    
$7$
 &  1.828e-01
& 1e-04 & 3.542e-03
& 4 & 1e-02 & 1e-04 & 1.688e-03
 	 \\
			 \hline    
$11$
 &  1.258e-01 
& 1e-02 & 2.213e-03 
& 4 & 1e-02 & 1e-04 & 1.675e-03  
 	 \\
			 \hline    
$17$
 &  6.551e-02 
& 1e-02 & 2.312e-03
& 4 & 1e-02 & 1e-04 & 1.971e-03  
	 \\
	 
			 \hline   
		\end{tabular}
	\end{center}
	\caption{Burgers equation, $\nu = 10^{-3}$, reconstructive regime, optimal $tol$, $tol_S$, and $tol_L$. 
	Average $L^2$ error for G-ROM, 2S-DD-VMS-ROM, and 3S-DD-VMS-ROM for different r values. \label{table:burgers-reconstructive-1c}
    }
		\end{table} 
In Figure \ref{fig:burgers-reconstructive}, we plot the time evolution of the solutions for the FOM  projection, G-ROM, 2S-DD-VMS-ROM, and 3S-DD-VMS-ROM for the reconstructive regime.
These plots show that both the 2S-DD-VMS-ROM and the 3S-DD-VMS-ROM are significantly more accurate than the standard G-ROM, as indicated by the results in Table \ref{table:burgers-reconstructive-1c}. 

\clearpage
		
In Tables~\ref{table:burgers-reconstructive-1a}, \ref{table:burgers-reconstructive-2a}, \ref{table:burgers-reconstructive-3a}, and \ref{table:burgers-reconstructive-4a}, we list the ROMs errors for the reconstructive regime with fixed $tol$ in the 2S-DD-VMS-ROM, and $tol_L = tol$ and optimal $tol_S$ in the 3S-DD-VMS-ROM.
We also list the optimal value 
of $tol_S$.
We consider the following values for $tol_L = tol$: $10^2$ (Table~\ref{table:burgers-reconstructive-1a}), $10^1$ (Table~\ref{table:burgers-reconstructive-2a}), $10^0$ (Table~\ref{table:burgers-reconstructive-3a}), and $10^{-1}$ (Table~\ref{table:burgers-reconstructive-4a}).  
These results yield the following conclusions:
For large $tol_L = tol$ values (i.e., $10^2$ and $10^1$), the 2S-DD-VMS-ROM is slightly more or as accurate as the G-ROM, whereas the 3S-DD-VMS-ROM is several times (and sometimes more than one order of magnitude) more accurate than the G-ROM and 2S-DD-VMS-ROM.
For small $tol_L = tol$ values (i.e., $10^0$ and $10^{-1}$), the 2S-DD-VMS-ROM is several times (and sometimes more than one order of magnitude) more accurate than the G-ROM.
Even in these cases, however, the 3S-DD-VMS-ROM is several times (and sometimes more than one order of magnitude) more accurate than the 2S-DD-VMS-ROM.
Overall, the 3S-DD-VMS-ROM is by far the most accurate ROM.

\begin{table}[h!]
\small 
	\begin{center}
		\smallskip
		\begin{tabular}{|c|c|c|c|c|c|c|}
			\hline
			\multicolumn{1}{|c|}{ $r$} &  
			\multicolumn{1}{c|}{ G-ROM } &  
			\multicolumn{1}{c|}{ 2S-DD-VMS-ROM } &  
			\multicolumn{3}{c|}{ 3S-DD-VMS-ROM}
						 \\  \cline{2-6}&{$\mathcal{E}(L^2)$ }   &{$\mathcal{E}(L^2)$ }  &{$r_1$}& {$tol_S$}&{ $\mathcal{E}(L^2)$ } 
									 \\
			 \hline  
$3$
 &  1.181e-01 
& 1.181e-01
& 1 &  1e+00 &  1.609e-02
 	 \\
			 \hline    
$7$
 &  1.828e-01
& 1.828e-01
& 1 & 1e-01 & 6.241e-03
 	 \\
			 \hline    
$11$
 &  1.258e-01 
& 1.258e-01 
& 1 & 1e-01 & 4.955e-03  
 	 \\
			 \hline    
$17$
 &  6.551e-02 
& 6.551e-02
& 1 & 1e-02 & 2.826e-03  
	 \\
	 
			 \hline   
		\end{tabular}
	\end{center}
	\caption{ Burgers equation, $\nu = 10^{-3}$, reconstructive regime, $tol = tol_L = 10^{2}$, and optimal $tol_S$.
	Average $L^2$ error for G-ROM, 2S-DD-VMS-ROM, and 3S-DD-VMS-ROM for
different r values.} \label{table:burgers-reconstructive-1a}
		\end{table}

\begin{table}[h!]
\small 
	\begin{center}
		\smallskip
		\begin{tabular}{|c|c|c|c|c|c|c|}
			\hline
			\multicolumn{1}{|c|}{ $r$} &  
			\multicolumn{1}{c|}{ G-ROM } &  
			\multicolumn{1}{c|}{ 2S-DD-VMS-ROM } &  
			\multicolumn{3}{c|}{ 3S-DD-VMS-ROM}
						 \\  \cline{2-6}&{$\mathcal{E}(L^2)$ }   &{$\mathcal{E}(L^2)$ }  &{$r_1$}& {$tol_S$}&{ $\mathcal{E}(L^2)$ } 
									 \\
			 \hline  
$3$
 &  1.181e-01 
& 7.278e-02
& 1 &  1e+00 &  1.322e-02
 	 \\
			 \hline    
$7$
 &  1.828e-01
& 1.755e-01
& 2 & 1e-03 & 3.915e-03
 	 \\
			 \hline    
$11$
 &  1.258e-01 
& 1.229e-01 
& 1 & 1e-03 & 1.787e-03  
 	 \\
			 \hline    
$17$
 &  6.551e-02 
& 6.456e-02
& 1 & 1e-02 & 2.310e-03  
	 \\
	 
			 \hline   
		\end{tabular}
	\end{center}
	\caption{ Burgers equation, $\nu = 10^{-3}$, reconstructive regime, $tol = tol_L = 10^{1}$, and optimal $tol_S$. Average $L^2$ error for G-ROM, 2S-DD-VMS-ROM, and 3S-DD-VMS-ROM for
different r values.} \label{table:burgers-reconstructive-2a}
		\end{table}

\begin{table}[h!]
\small 
	\begin{center}
		\smallskip
		\begin{tabular}{|c|c|c|c|c|c|c|}
			\hline
			\multicolumn{1}{|c|}{ $r$} &  
			\multicolumn{1}{c|}{ G-ROM } &  
			\multicolumn{1}{c|}{ 2S-DD-VMS-ROM } &  
			\multicolumn{3}{c|}{ 3S-DD-VMS-ROM}
						 \\  \cline{2-6}&{$\mathcal{E}(L^2)$ }   &{$\mathcal{E}(L^2)$ }  &{$r_1$}& {$tol_S$}&{ $\mathcal{E}(L^2)$ } 
									 \\
			 \hline  
$3$
 &  1.181e-01 
& 1.333e-01
& 1 &  1e-02 &  5.292e-03
 	 \\
			 \hline    
$7$
 &  1.828e-01
& 2.590e-02
& 2 & 1e-03 & 3.549e-03
 	 \\
			 \hline    
$11$
 &  1.258e-01 
& 3.607e-02 
& 2 & 1e-02 & 2.045e-03  
 	 \\
			 \hline    
$17$
 &  6.551e-02 
& 5.029e-02
& 5 & 1e-02 & 2.237e-03  
	 \\
	 
			 \hline   
		\end{tabular}
	\end{center}
	\caption{ Burgers equation, $\nu = 10^{-3}$, reconstructive regime, $tol = tol_L = 10^{0}$, and optimal $tol_S$. Average $L^2$ error for G-ROM, 2S-DD-VMS-ROM, and 3S-DD-VMS-ROM for
different r values.} \label{table:burgers-reconstructive-3a}
		\end{table}

\begin{table}[h!]
\small 
	\begin{center}
		\smallskip
		\begin{tabular}{|c|c|c|c|c|c|c|}
			\hline
			\multicolumn{1}{|c|}{ $r$} &  
			\multicolumn{1}{c|}{ G-ROM } &  
			\multicolumn{1}{c|}{ 2S-DD-VMS-ROM } &  
			\multicolumn{3}{c|}{ 3S-DD-VMS-ROM}
						 \\  \cline{2-6}&{$\mathcal{E}(L^2)$ }   &{$\mathcal{E}(L^2)$ }  &{$r_1$}& {$tol_S$}&{ $\mathcal{E}(L^2)$ } 
									 \\
			 \hline  
$3$
 &  1.181e-01 
& 3.729e-03
& 1 &  1e-02 &  2.061e-03
 	 \\
			 \hline    
$7$
 &  1.828e-01
& 4.232e-03
& 4 & 1e-03 & 2.557e-03
 	 \\
			 \hline    
$11$
 &  1.258e-01 
& 4.556e-03 
& 2 & 1e-02 & 2.086e-03  
 	 \\
			 \hline    
$17$
 &  6.551e-02 
& 5.962e-03
& 5 & 1e-02 & 2.255e-03  
	 \\
	 
			 \hline   
		\end{tabular}
	\end{center}
	\caption{ Burgers equation, $\nu = 10^{-3}$, reconstructive regime $tol = tol_L = 10^{-1}$, and optimal $tol_S$. Average $L^2$ error for G-ROM, 2S-DD-VMS-ROM, and 3S-DD-VMS-ROM for different r values.} \label{table:burgers-reconstructive-4a}
		\end{table} 
		\color{black}


In Tables~\ref{table:burgers-reconstructive-1b}, \ref{table:burgers-reconstructive-2b}, and  \ref{table:burgers-reconstructive-3b}, we list the ROMs errors for the reconstructive regime with fixed $tol$ in the 2S-DD-VMS-ROM, and $tol_S = tol$ and optimal $tol_L$ in the 3S-DD-VMS-ROM.
We also list the optimal value 
of $tol_L$.
We consider the following values for $tol_S = tol$: $10^0$ (Table~\ref{table:burgers-reconstructive-1b}), $10^{-1}$ (Table~\ref{table:burgers-reconstructive-2b}), and $10^{-2}$ (Table~\ref{table:burgers-reconstructive-3b}).  
These results yield the following conclusions:
For all $tol_L = tol$ values and all $r$ values,
the 2S-DD-VMS-ROM is several times (and sometimes more than one order of magnitude) more accurate than the G-ROM.
Furthermore, the 3S-DD-VMS-ROM is significantly (and sometimes several times) more accurate than the 2S-DD-VMS-ROM.
Overall, the 3S-DD-VMS-ROM is the most accurate ROM.



\begin{table}[h!]
\small 
	\begin{center}
		\smallskip
		\begin{tabular}{|c|c|c|c|c|c|c|}
			\hline
			\multicolumn{1}{|c|}{ $r$} &  
			\multicolumn{1}{c|}{ G-ROM } &  
			\multicolumn{1}{c|}{ 2S-DD-VMS-ROM } &  
			\multicolumn{3}{c|}{ 3S-DD-VMS-ROM}
						 \\  \cline{2-6}&{$\mathcal{E}(L^2)$ }   &{$\mathcal{E}(L^2)$ }  &{$r_1$}& {$tol_L$}&{ $\mathcal{E}(L^2)$ } 
									 \\
			 \hline  
$3$
 &  1.181e-01 
& 1.133e-02
& 2 &  1e-01 &  8.085e-03
 	 \\
			 \hline    
$7$
 &  1.828e-01
&2.590e-02
& 6 & 1e-03 & 1.762e-02
 	 \\
			 \hline    
$11$
 &  1.258e-01 
& 3.607e-02 
& 10 & 1e-02 & 2.390e-02  
 	 \\
			 \hline    
$17$
 &  6.551e-02 
& 5.029e-02
& 16 & 1e-02 & 1.486e-02  
	 \\
	 
			 \hline   
		\end{tabular}
	\end{center}
	\caption{ Burgers equation, $\nu = 10^{-3}$, reconstructive regime: $tol = tol_S = 10^{0}$ and optimal $tol_L$. Average $L^2$ error for G-ROM, 2S-DD-VMS-ROM, and 3S-DD-VMS-ROM for
different r values.} \label{table:burgers-reconstructive-1b}
		\end{table}

\begin{table}[h!]
\small 
	\begin{center}
		\smallskip
		\begin{tabular}{|c|c|c|c|c|c|c|}
			\hline
			\multicolumn{1}{|c|}{ $r$} &  
			\multicolumn{1}{c|}{ G-ROM } &  
			\multicolumn{1}{c|}{ 2S-DD-VMS-ROM } &  
			\multicolumn{3}{c|}{ 3S-DD-VMS-ROM}
						 \\  \cline{2-6}&{$\mathcal{E}(L^2)$ }   &{$\mathcal{E}(L^2)$ }  &{$r_1$}& {$tol_L$}&{ $\mathcal{E}(L^2)$ } 
									 \\
			 \hline  
$3$
 &  1.181e-01 
& 3.729e-03
& 2 &  1e-02 &  2.568e-03
 	 \\
			 \hline    
$7$
 &  1.828e-01
& 4.232e-03
& 4 & 1e-03 & 3.678e-03
 	 \\
			 \hline    
$11$
 &  1.258e-01 
& 4.556e-03 
& 10 & 1e-02 & 3.995e-03  
 	 \\
			 \hline    
$17$
 &  6.551e-02 
& 5.962e-03
& 16 & 1e-02 & 2.995e-03  
	 \\
	 
			 \hline   
		\end{tabular}
	\end{center}
	\caption{ Burgers equation, $\nu = 10^{-3}$, reconstructive regime: $tol = tol_S = 10^{-1}$ and optimal $tol_L$. Average $L^2$ error for G-ROM, 2S-DD-VMS-ROM, and 3S-DD-VMS-ROM for
different r values.} \label{table:burgers-reconstructive-2b}
		\end{table}

\begin{table}[h!]
\small 
	\begin{center}
		\smallskip
		\begin{tabular}{|c|c|c|c|c|c|c|}
			\hline
			\multicolumn{1}{|c|}{ $r$} &  
			\multicolumn{1}{c|}{ G-ROM } &  
			\multicolumn{1}{c|}{ 2S-DD-VMS-ROM } &  
			\multicolumn{3}{c|}{ 3S-DD-VMS-ROM}
						 \\  \cline{2-6}&{$\mathcal{E}(L^2)$ }   &{$\mathcal{E}(L^2)$ }  &{$r_1$}& {$tol_L$}&{ $\mathcal{E}(L^2)$ } 
									 \\
			 \hline  
$3$
 &  1.181e-01 
& 1.548e-03
& 1 &  1e-02 &  1.548e-03
 	 \\
			 \hline    
$7$
 &  1.828e-01
& 1.062e-02
& 4 & 1e-04 & 1.688e-03
 	 \\
			 \hline    
$11$
 &  1.258e-01 
& 2.213e-03 
& 4 & 1e-04 & 1.675e-03  
 	 \\
			 \hline    
$17$
 &  6.551e-02 
& 2.312e-03
& 4 & 1e-04 & 1.974e-03  
	 \\
	 
			 \hline   
		\end{tabular}
	\end{center}
	\caption{ Burgers equation, $\nu = 10^{-3}$, reconstructive regime: $tol = tol_S = 10^{-2}$ and optimal $tol_L$. Average $L^2$ error for G-ROM, 2S-DD-VMS-ROM, and 3S-DD-VMS-ROM for
different r values.} \label{table:burgers-reconstructive-3b}
		\end{table}

The results in Tables~\ref{table:burgers-reconstructive-1a}--\ref{table:burgers-reconstructive-3b} suggest that there is no apparent relationship between the 2S-DD-VMS-ROM tolerance $tol$ and the 3S-DD-VMS-ROM tolerances $tol_L$ and $tol_S$.
We intend to perform a more thorough investigation of potential relationships among these tolerances in a future study.


In Table~\ref{table:burgers-cross-validation}, we  list the ROMs errors for the cross-validation regime with optimal $tol$ in Algorithm~\ref{alg:2s-dd-vms-rom} and optimal $tol_L$ and $tol_s$ in Algorithm~\ref{alg:3s-dd-vms-rom}.
These results show that, for all $r$ values, the 2S-DD-VMS-ROM  and 3S-DD-VMS-ROM are several times (sometimes even one order of magnitude) more accurate than the standard G-ROM. 
Overall, the 3S-DD-VMS-ROM is more accurate than the 2S-DD-VMS-ROM.



\begin{table}[h!]
\small 
	\begin{center}
		\smallskip
		\begin{tabular}{|c|c|c|c|c|c|c|c|c|}
			\hline
			\multicolumn{1}{|c|}{ $r$} &  
			\multicolumn{1}{c|}{ G-ROM } &  
			\multicolumn{2}{c|}{ 2S-DD-VMS-ROM } &  
			\multicolumn{4}{c|}{ 3S-DD-VMS-ROM}
						 \\  \cline{2-8}&{$\mathcal{E}(L^2)$ }   & {$tol$}& {$\mathcal{E}(L^2)$ }  &{$r_1$}& {$tol_S$}& {$tol_L$}&{ $\mathcal{E}(L^2)$ } 
									 \\
			 \hline  
$3$
 &  2.015e-01 
& 1e-01 & 2.028e-02
& 2 &  1e-01 & 1e+00 &  1.863e-02
 	 \\
			 \hline    
$7$
 &  1.796e-01
& 5e-02 & 1.400e-02
& 3  & 5e-02 & 1e+00 & 1.188e-02
 	 \\
			 \hline    
$11$
 &  1.163e-01 
& 3e-02& 8.981e-03
& 6  & 3e-02 &  1e+00 & 8.383e-03
 	 \\
			 \hline    
$17$
 &  6.897e-02
&  1e-02 & 8.542e-03
&  6 &  1e-02 &  1e+00 & 8.452e-03
	 \\
	 
			 \hline   
		\end{tabular}
	\end{center}
	\caption{ Burgers equation, $\nu = 10^{-3}$, cross-validation regime, optimal $tol$, $tol_S$, and $tol_L$. 
	Average $L^2$ error for G-ROM, 2S-DD-VMS-ROM, and 3S-DD-VMS-ROM for different r values.} \label{table:burgers-cross-validation}
		\end{table} 


In Table~\ref{table:burgers-predictive}, we  list the ROMs errors for the predictive regime with optimal $tol$ in Algorithm~\ref{alg:2s-dd-vms-rom} and optimal $tol_L$ and $tol_s$ in Algorithm~\ref{alg:3s-dd-vms-rom}.
These results show that, for all $r$ values, the 2S-DD-VMS-ROM  and 3S-DD-VMS-ROM are several times (sometimes even one order of magnitude) more accurate than the standard G-ROM. 
Overall, the 3S-DD-VMS-ROM is more accurate than the 2S-DD-VMS-ROM.


\begin{table}[h!]
\small 
	\begin{center}
		\smallskip
		\begin{tabular}{|c|c|c|c|c|c|c|c|c|}
			\hline
			\multicolumn{1}{|c|}{ $r$} &  
			\multicolumn{1}{c|}{ G-ROM } &  
			\multicolumn{2}{c|}{ 2S-DD-VMS-ROM } &  
			\multicolumn{4}{c|}{ 3S-DD-VMS-ROM}
						 \\  \cline{2-8}&{$\mathcal{E}(L^2)$ }   & {$tol$}& {$\mathcal{E}(L^2)$ }  &{$r_1$}& {$tol_S$}& {$tol_L$}&{ $\mathcal{E}(L^2)$ } 
									 \\
			 \hline  
$3$
 &  2.185e-01 
& 1e-01 & 3.623e-02
& 2 &  1e-01 & 1e+00 &  3.029e-02
 	 \\
			 \hline    
$7$
 &  2.054e-01
& 3e-02 & 2.004e-02
& 6  & 5e-02 & 3e-02 & 1.428e-02
 	 \\
			 \hline    
$11$
 &  1.620e-01 
& 3e-02& 1.608e-02
& 10  & 5e-02 &  3e-02 & 1.418e-02
 	 \\
			 \hline    
$17$
 &  1.103e-01
&  1e-02 & 1.524e-02
&  6 &  1e-02 &  1e-01 & 1.506e-02
	 \\
	 
			 \hline   
		\end{tabular}
	\end{center}
	\caption{ Burgers equation, $\nu = 10^{-3}$, predictive regime, optimal $tol$, $tol_S$, and $tol_L$. 
	Average $L^2$ error for G-ROM, 2S-DD-VMS-ROM, and 3S-DD-VMS-ROM for different r values.} \label{table:burgers-predictive}
		\end{table}

\subsection{Flow Past A Cylinder}
	\label{sec:numerical-results-nse}

In this section, we investigate the 2S-DD-VMS-ROM~\eqref{eqn:vms-rom-two-scales-3} and the new 3S-DD-VMS-ROM~\eqref{eqn:vms-rom-three-scales-7} in the numerical simulation of a 2D channel flow past a circular cylinder at Reynolds numbers $Re=100, Re=500$, and $Re=1000$. 

\paragraph{Computational Setting}
As a mathematical model, we use the NSE~\eqref{eqn:nse-1}--\eqref{eqn:nse-2}.
The computational domain is a $2.2\times 0.41$ rectangular channel with a radius $=0.05$ cylinder, centered at $(0.2,0.2)$, see Figure~\ref{cyldomain}.  
{
\color{black}
\begin{figure}[h!]
\centering
\begin{tikzpicture}
  \coordinate (a) at (0,0);
  \coordinate (b) at (11,0);
  \coordinate (c) at (11,2.05);
  \coordinate (d) at (0,2.05);
  \coordinate(bb) at  (11.4,0);
  \coordinate(cc1) at  (11.4,2.05);
  \coordinate(cc2) at  (11,2.45);
  \coordinate(dd) at  (0,2.45);
  \coordinate(cir) at  (1,1);
  \coordinate(cir1) at  (1,0);
  \coordinate(cir2) at  (0,1);

  \draw[line width=1.2pt, mark=none] (a) -- (b) -- (c) -- (d)  -- cycle;
  \draw[line width=1pt, mark=none] (b) -- (bb) ;
  \draw[line width=1pt, mark=none] (c) -- (cc1) ;
  \draw[line width=1pt, mark=none] (c) -- (cc2) ;
  \draw[line width=1pt, mark=none] (d) -- (dd) ;
  \draw[line width=1pt, mark=none] (a) -- (-0.5,0) ;
  \draw[line width=1pt, mark=none] (cir2) -- (-0.5,1) ;

  \draw[latex-latex,line width=0.8pt, mark=none] ([xshift=-0.3cm]cir2) -- ([xshift=-0.3cm]a) node[midway,right,xshift=-0.8cm]{\small\SI{0.2}{}};
  \draw[latex-latex,line width=0.8pt, mark=none] (cir) -- (cir2) node[midway,above,yshift=0.2cm]{\small\SI{0.2}{}};
  \draw[latex-latex,line width=0.2pt, mark=none] (cir) -- (1.2121,1.2121) node[above,yshift=-0.03cm,xshift=0.42cm]{\small\SI{0.05}{}};
  
    \node at (cir) {
    \begin{tikzpicture}
      \draw[line width=1.2pt, mark=none] circle (0.3cm);
    \end{tikzpicture}
      };
      
  \draw[latex-latex,line width=0.8pt, mark=none] ([xshift=+0.25cm]b) -- ([xshift=+0.25cm]c) node[midway, right]{\SI{0.41}{}};
 \draw[latex-latex,line width=0.8pt, mark=none] ([yshift=+0.25cm]c) -- ([yshift=+0.25cm]d) node[midway, above]{\SI{2.2}{}};
\end{tikzpicture}
	\caption{\label{cyldomain} 
		Geometry of the flow past a circular cylinder numerical experiment.
		}
\end{figure}
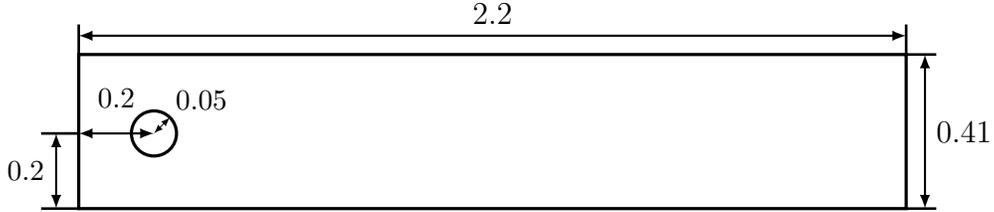
}

We  prescribe no-slip boundary conditions on the walls and cylinder, and the following inflow and outflow profiles~\cite{john2004reference,mohebujjaman2019physically,mohebujjaman2017energy,rebholz2017improved}:
\begin{align}
u_{1}(0,y,t)&=u_{1}(2.2,y,t)=\frac{6}{0.41^{2}}y(0.41-y), \\ u_{2}(0,y,t)&=u_{2}(2.2,y,t)=0,
\end{align} 
where $\bu=\langle u_1, u_2 \rangle$.  
There is no forcing and the flow starts from rest.

\paragraph{Snapshot Generation}
For the spatial discretization, we use the pointwise divergence-free, LBB stable $(P_2, P_1^{disc})$ Scott-Vogelius finite element pair on a barycenter refined regular triangular mesh~\cite{john2016divergence}. 
The mesh provides $103K$ ($102962$) velocity and $76K$ ($76725$) pressure degrees of freedom. 
We utilize the commonly used linearized BDF2 temporal discretization and a time step size $\Delta t=0.002$ for both FOM and ROM time discretizations. 
On the first time step, we use a backward Euler scheme so that we have two initial time step solutions required for the BDF2 scheme.

\paragraph{ROM Construction}

The FOM simulations achieve the statistically steady state at different time instances for the three Reynolds numbers used in the numerical investigation:
For $Re = 100$, after $t = 5s$; for $Re = 500$, after $t=7s$; and for $Re=1000$, after $t=13s$.
To build the ROM basis functions, we decided to use $10 s$ of FOM data.
Thus, to ensure a fair comparison of the numerical results at different Reynolds numbers, we collect FOM snapshots on the following time intervals:
For $Re=100$, from $t=7$ to $t=17$; for $Re=500$, from $t=7$ to $t=17$; and for $Re=1000$, from $t=13$ to $t=23$.

To train $\tA,\tB$ (for the 2S-DD-VMS-ROM) and $\tA_L,\tB_L$ and $\tA_S,\tB_S$ (for the 3S-DD-VMS-ROM), we use FOM  data for one period in the reconstructive and cross-validation regimes, and FOM data for half a period in the predictive regime.
We note that the period length of the statistically steady state is different for the three different Reynolds numbers:
From $t = 7$ to $t = 7.332$ for $Re = 100$; from $t = 7$ to $t = 7.442$ for $Re = 500$; and from $t = 13$ to $t = 13.268$ for $Re = 1000$. 
Thus, the reconstructive and cross-validation regimes, we collect $167$ snapshots for $Re = 100$; $222$ snapshots for $Re = 500$; and $135$ snapshots for $Re = 1000$.
For the predictive regime, we collect $84$ snapshots for $Re = 100$; $111$ snapshots for $Re = 500$; and $68$ snapshots for $Re = 1000$.

\bigskip



\subsubsection{Numerical Results for  $Re=100$} 
	\label{sec:numerical-results-nse-re100}


In this section, we present numerical results for the flow past a cylinder at $Re=100$.

In  Table~\ref{table:nse-reconstructive-re100}, for different $r$ values, we list the average $L^{2}$ error~\eqref{eqn:l2-error}
for the G-ROM, the 2S-DD-VMS-ROM, and the new 3S-DD-VMS-ROM in the reconstructive regime.
We also list 
the $r_1$ values for the 3S-DD-VMS-ROM.
These results show that, for all $r$ values, the 2S-DD-VMS-ROM  and 3S-DD-VMS-ROM are several times (sometimes one and even two orders of magnitude) more accurate than the standard G-ROM.
Furthermore, the 3S-DD-VMS-ROM is generally more accurate than the 2S-DD-VMS-ROM, especially for large $r$ values: 
For example, for $r=8$, the 3S-DD-VMS-ROM is {\it more than twice more accurate} than the 2S-DD-VMS-ROM.
We also note that the ROM errors in Table~\ref{table:nse-reconstructive-re100} converge to $0$ according to an even/odd pattern:
The ROM errors for even $r$ values converge to $0$ and the ROM errors for odd $r$ values also converge to $0$.
This behavior is related to the flow past a cylinder configuration, in which the ROM modes appear in pairs.
We emphasize, however, that for large $r$ values, we recover the asymptotic convergence that does not depend on the odd/even $r$ values, just as in the Burgers equation test case (Section~\ref{sec:burgers}). 

\begin{table}[h!]
\small
	\begin{center}
		\smallskip
		\begin{tabular}{|c|c|c|c|c|c|c|c|}
			\hline
			\multicolumn{1}{|c|}{ $r$} &  
			\multicolumn{1}{c|}{ G-ROM } &  
			\multicolumn{1}{c|}{ 2S-DD-VMS-ROM } &  
			\multicolumn{2}{c|}{3S-DD-VMS-ROM}
						 \\  \cline{2-5}&{$\mathcal{E}(L^2)$ }  
						 &{$\mathcal{E}(L^2)$ }  
						 &{$r_1$} & { $\mathcal{E}(L^2)$ } 
									 \\
			 \hline  
$2$
&9.902e-02 
&5.118e-04 
&1&5.088e-04
 	 \\
			 \hline    
$3$
&1.029e-01
&3.208e-02
&2& 3.018e-02 
 	 \\
			 \hline    
$4$
&5.840e-02  
&1.553e-03   
&2&1.479e-03
 	 \\
			 \hline    
$5$
&6.492e-02 
&2.270e-02  
&4&2.191e-02   
	 \\
			 \hline  
$6$
&1.370e-02   
&5.336e-04   
&1&4.804e-04 
	 \\
			 \hline   
$7$
&1.403e-02 
&6.038e-03  
&6& 5.817e-03 
	 \\
			 \hline  
$8$
&1.214e-02  
& 9.302e-04  
&6&4.415e-04  
\\	\hline	 
		\end{tabular}
	\end{center}
	\caption {
			Flow past a cylinder, $Re=100$, reconstructive regime.
			Average $L^{2}$ errors for G-ROM, 2S-DD-VMS-ROM, and 3S-DD-VMS-ROM for different $r$ values. 
	\label{table:nse-reconstructive-re100}
	} 
\end{table}

In  Table~\ref{table:nse-cross-re100}, for different $r$ values, we list the average $L^{2}$ error~\eqref{eqn:l2-error} for the G-ROM, the 2S-DD-VMS-ROM, and the new 3S-DD-VMS-ROM in the cross-validation regime. 
These results show that, for all $r$ values, the 2S-DD-VMS-ROM  and 3S-DD-VMS-ROM are several times (sometimes even two orders of magnitude) more accurate than the standard G-ROM.
Furthermore, the 3S-DD-VMS-ROM is generally more accurate than the 2S-DD-VMS-ROM, especially for large $r$ values: 
For example, for $r=8$, the 3S-DD-VMS-ROM is almost {\it three times more accurate} than the 2S-DD-VMS-ROM.

\begin{table}[h!]
\small
	\begin{center}
		\smallskip
		\begin{tabular}{|c|c|c|c|c|c|c|c|}
			\hline
			\multicolumn{1}{|c|}{ $r$} &  
			\multicolumn{1}{c|}{ G-ROM } &  
			\multicolumn{1}{c|}{ 2S-DD-VMS-ROM } &  
			\multicolumn{2}{c|}{3S-DD-VMS-ROM}
						 \\  \cline{2-5}&{$\mathcal{E}(L^2)$ }  
						 &{$\mathcal{E}(L^2)$ }  
						 &{$r_1$} & { $\mathcal{E}(L^2)$ } 
									 \\
			 \hline  
$2$
&4.891e-01 
&1.536e-03  
&1&   1.458e-03
 	 \\
			 \hline    
$3$
&4.088e-01 
&   3.514e-02
&2&   3.106e-02
 	 \\
			 \hline    
$4$
&9.291e-02
&   2.187e-03  
&2&   2.015e-03 
 	 \\
			 \hline    
$5$
&1.013e-01 
&   2.279e-02
&4&   2.220e-02
	 \\
			 \hline  
$6$
& 3.270e-02  
&   5.113e-04
&2&   4.921e-04 
	 \\
			 \hline   
$7$
&3.059e-02 
&   7.476e-03  
&1&   7.260e-03 
	 \\
			 \hline  
$8$
&3.600e-02  
&  1.221e-03
&6&   4.385e-04  
\\	\hline	 
		\end{tabular}
	\end{center}
	\caption{
	Flow past a cylinder, $Re=100$, cross-validation regime.
			Average $L^{2}$ errors for G-ROM, 2S-DD-VMS-ROM, and 3S-DD-VMS-ROM for different $r$ values. 
			 \label{table:nse-cross-re100}
	}
	\end{table}

In  Table~\ref{table:nse-predictive-re100}, for different $r$ values, we list the average $L^{2}$ error~\eqref{eqn:l2-error} for the G-ROM, the 2S-DD-VMS-ROM, and the new 3S-DD-VMS-ROM in the predictive regime. 
These results show that, for all $r$ values, the 2S-DD-VMS-ROM  and 3S-DD-VMS-ROM are several times (sometimes even one order of magnitude) more accurate than the standard G-ROM. 
Furthermore, the 3S-DD-VMS-ROM is generally more accurate than the 2S-DD-VMS-ROM: Specifically, for $r\ge 4$, the 3S-DD-VMS-ROM is at least {\it twice more accurate} than the  2S-DD-VMS-ROM.

\begin{table}[h!]
\small
	\begin{center}
		\smallskip
		\begin{tabular}{|c|c|c|c|c|c|c|c|}
			\hline
			\multicolumn{1}{|c|}{ $r$} &  
			\multicolumn{1}{c|}{ G-ROM } &  
			\multicolumn{1}{c|}{ 2S-DD-VMS-ROM } &  
			\multicolumn{2}{c|}{3S-DD-VMS-ROM}
						 \\  \cline{2-5}&{$\mathcal{E}(L^2)$ }  
						 &{$\mathcal{E}(L^2)$ }  
						 &{$r_1$} & { $\mathcal{E}(L^2)$ } 
									 \\
			 \hline  
$2$
& 3.883e-01 
&8.172e-02 
&1& 7.388e-02 
 	 \\
			 \hline    
$3$
&3.616e-01 
&3.374e-02  
&1&   3.141e-02 
 	 \\
			 \hline    
$4$
&1.366e-01  
&8.127e-03 
&2& 4.115e-03  
 	 \\
			 \hline    
$5$
&1.464e-01  
&4.248e-02 
&3& 2.602e-02  
	 \\
			 \hline  
$6$
&1.348e-01  
&5.946e-03  
&4&1.051e-03  
	 \\
			 \hline   
$7$
& 1.291e-01  
& 1.529e-02  
&3& 6.613e-03 
	 \\
			 \hline  
$8$
&9.638e-02  
&6.798e-03  
&6&3.170e-03 
\\	\hline	 
		\end{tabular}
	\end{center}
	\caption {
			Flow past a cylinder, $Re=100$, predictive regime.
			Average $L^{2}$ errors for G-ROM, 2S-DD-VMS-ROM, and 3S-DD-VMS-ROM for different $r$ values. 
	\label{table:nse-predictive-re100}
	} 
	\end{table}

\begin{figure}[H]
	\begin{center}
		\includegraphics[width=0.9\textwidth,height=0.22\textwidth]{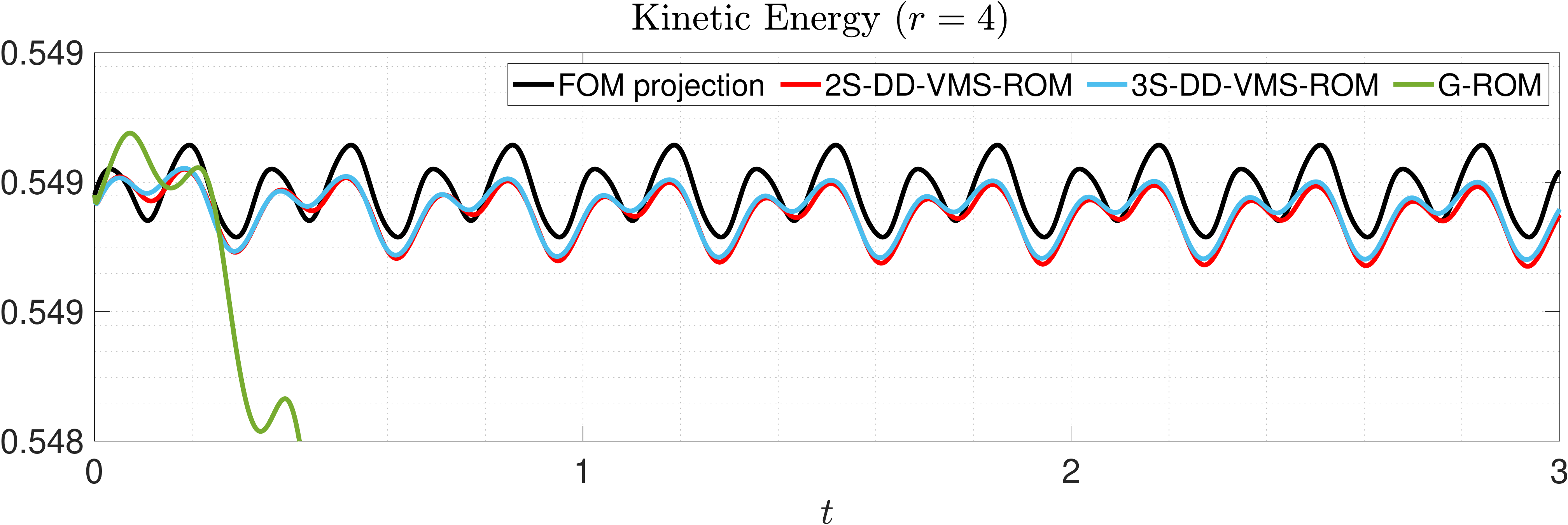}
		\includegraphics[width=0.9\textwidth,height=0.22\textwidth]{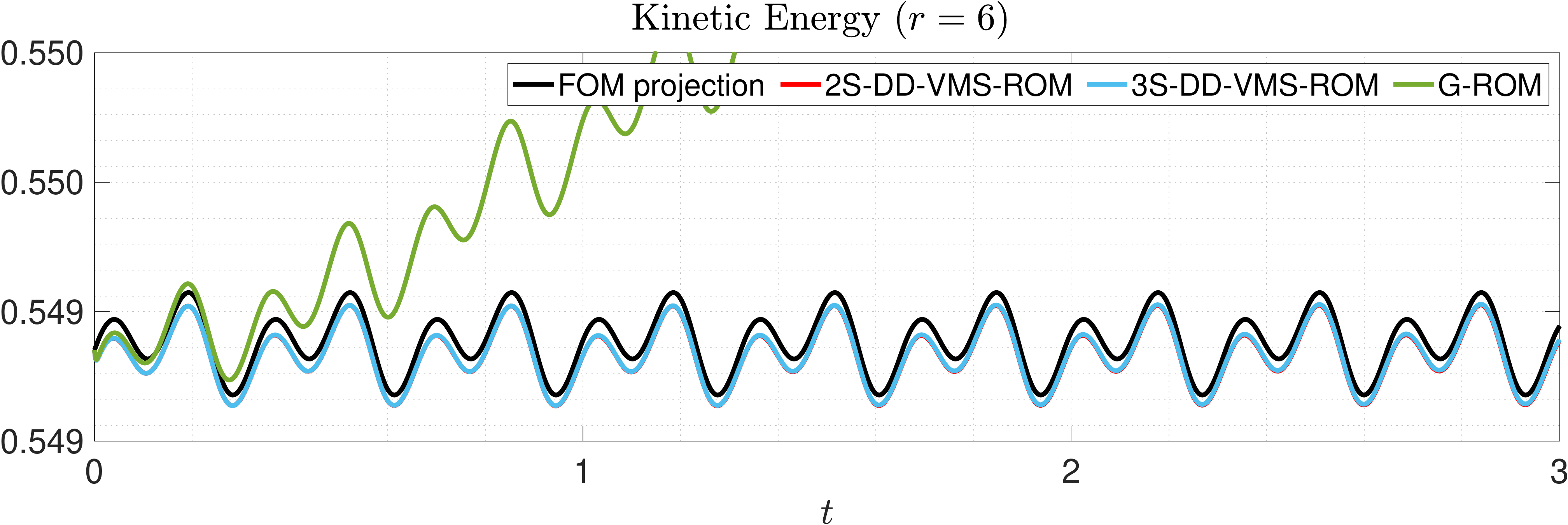}
		\includegraphics[width=0.9\textwidth,height=0.22\textwidth]{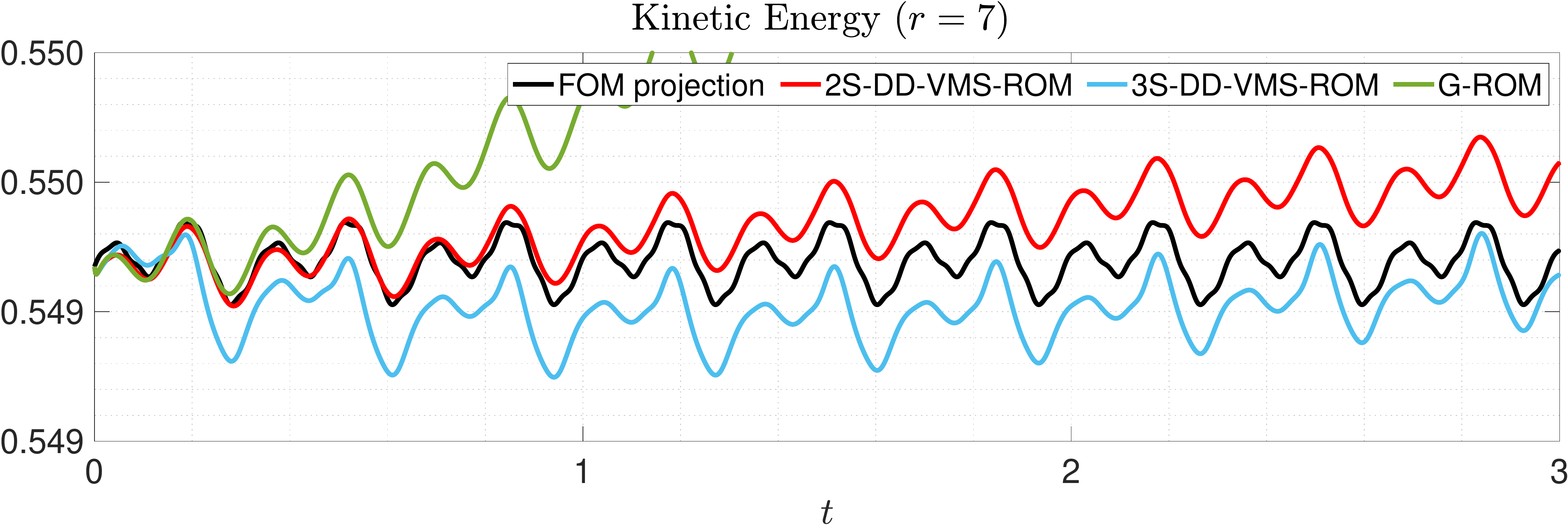}
	\end{center}
	\caption {
			Flow past a cylinder, $Re=100$, reconstructive regime.
			Time evolution of the kinetic energy for FOM projection, G-ROM, 2S-DD-VMS-ROM, and 3S-DD-VMS-ROM for different $r$ values. 
	\label{fig:nse-reconstructive-re100-ke}
	} 
\end{figure}

In Figure~\ref{fig:nse-reconstructive-re100-ke}, for $r=4,6,7$, we plot the time evolution of the kinetic energy of the 
FOM,
the G-ROM, the 2S-DD-VMS-ROM, and the new 3S-DD-VMS-ROM in the reconstructive regime.
These plots support the conclusions in Table~\ref{table:nse-reconstructive-re100}: Both the 3S-DD-VMS-ROM and the 2S-DD-VMS-ROM accurately approximate the 
FOM kinetic energy and are significantly more accurate than the standard G-ROM. 
Furthermore, 3S-DD-VMS-ROM is slightly more accurate than the 2S-DD-VMS-ROM, especially for $r=7$.

\begin{figure}[H]
	\begin{center}
		\includegraphics[width=0.9\textwidth,height=0.22\textwidth]{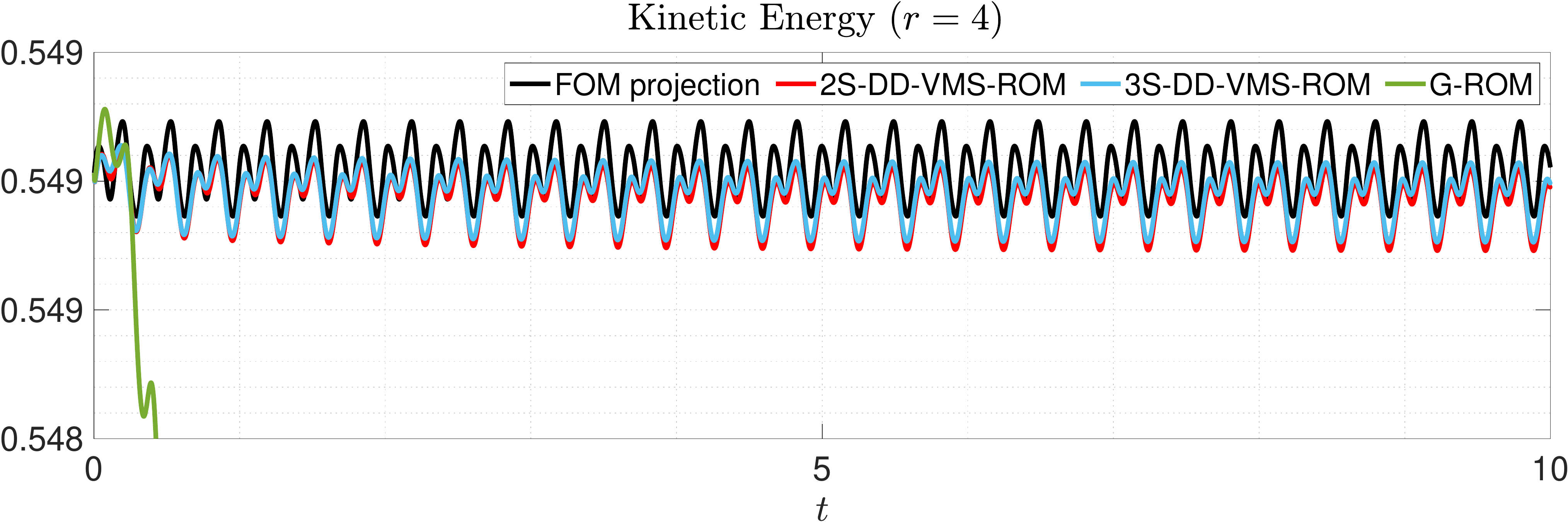}
		\includegraphics[width=0.9\textwidth,height=0.22\textwidth]{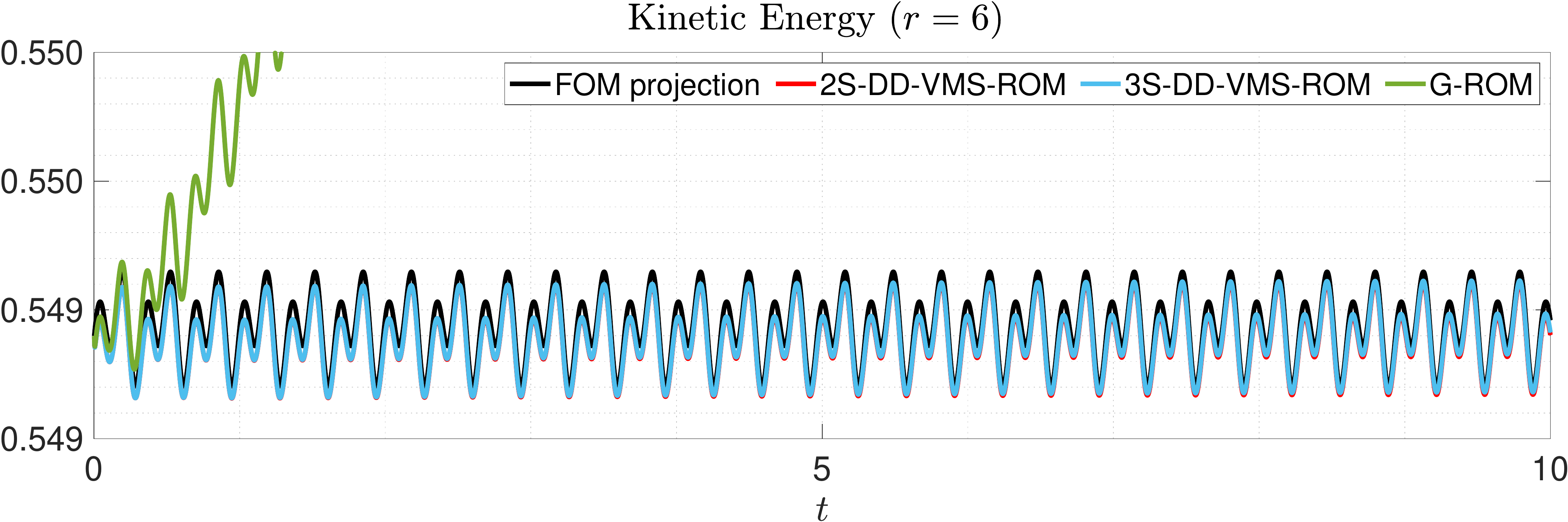}
		\includegraphics[width=0.9\textwidth,height=0.22\textwidth]{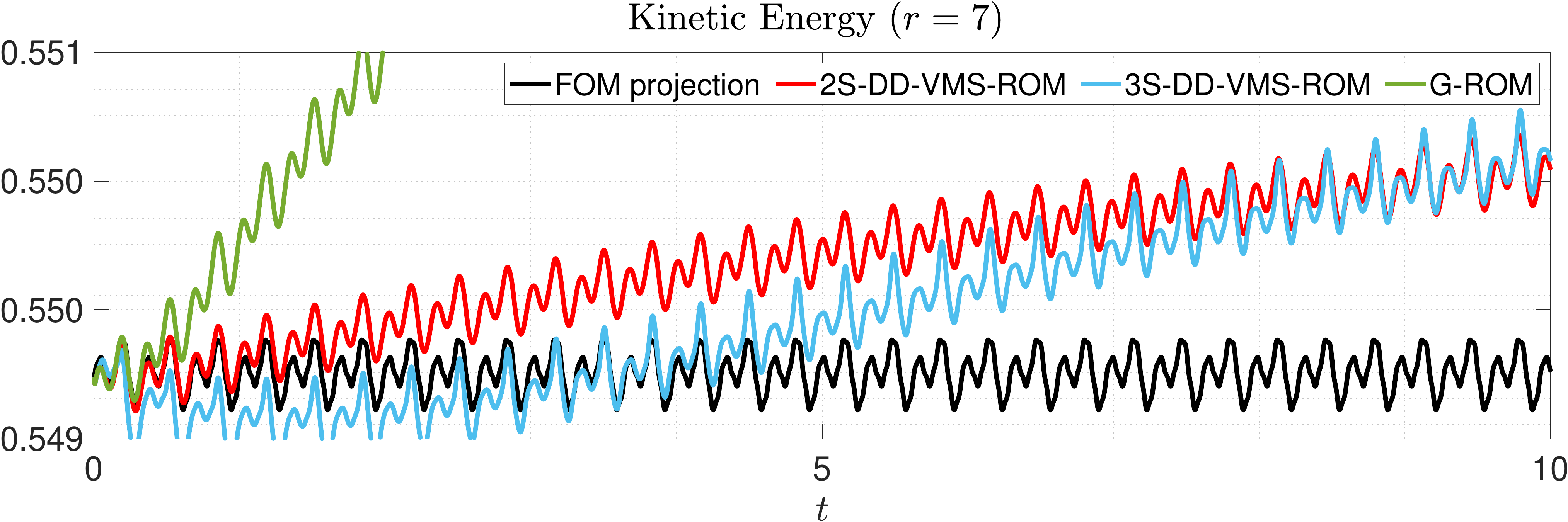}
	\end{center}
	\caption {
			Flow past a cylinder, $Re=100$, cross-validation regime.
			Time evolution of the kinetic energy for FOM projection, G-ROM, 2S-DD-VMS-ROM, and 3S-DD-VMS-ROM for different $r$ values. 
	\label{fig:nse-cross-re100-ke}
	} 
\end{figure}

 In Figure~\ref{fig:nse-cross-re100-ke}, for $r=4, 6, 7$, we plot the time evolution of the kinetic energy of the 
FOM
 the G-ROM, the 2S-DD-VMS-ROM, and the new 3S-DD-VMS-ROM in the cross-validation regime. 
 For all cases, the evolution of the G-ROM kinetic energy is very inaccurate. 
 In contrast, for $r=4$ and $r=6$, both the 3S-DD-VMS-ROM and the 2S-DD-VMS-ROM successfully reproduce the 
FOM
 kinetic energy. For $r=7$, the 3S-DD-VMS-ROM accurately approximates the 
FOM
 kinetic energy between $t=0$ and $t=4$. 
 For $t \ge 4$, although the 3S-DD-VMS-ROM and 2S-DD-VMS-ROM kinetic energy approximations are not as accurate, they are still much more accurate than the G-ROM kinetic energy approximation.  

\begin{figure}[H]
	\begin{center}
		\includegraphics[width=0.9\textwidth,height=0.22\textwidth]{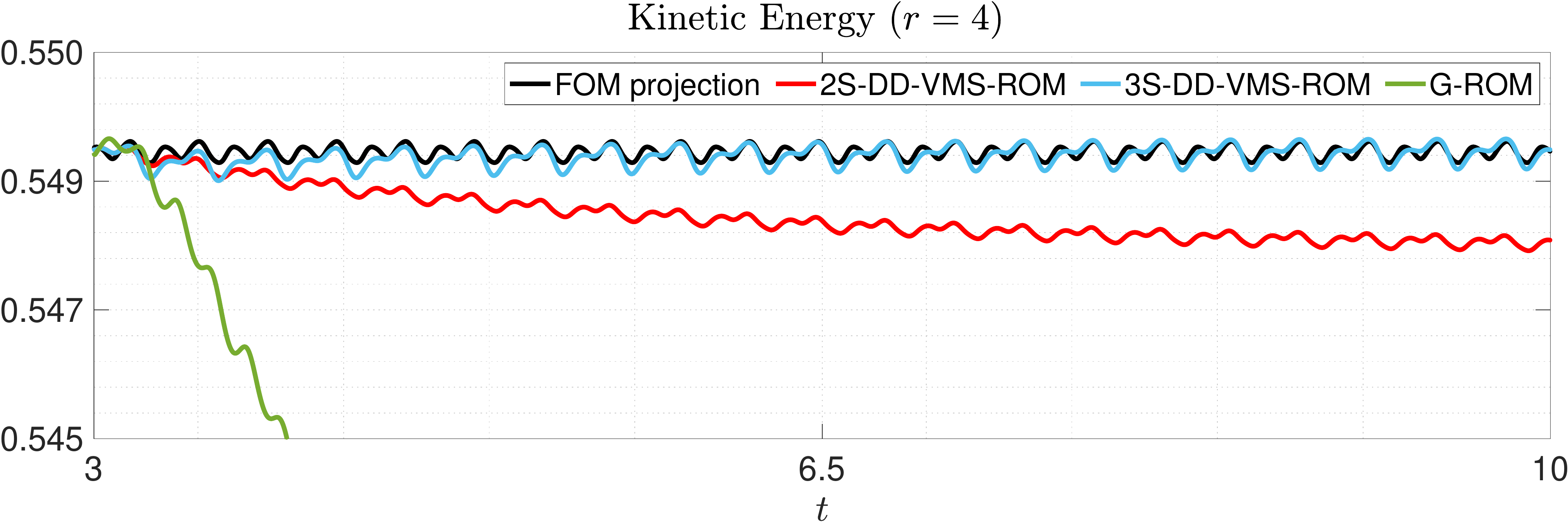}
		\includegraphics[width=0.9\textwidth,height=0.22\textwidth]{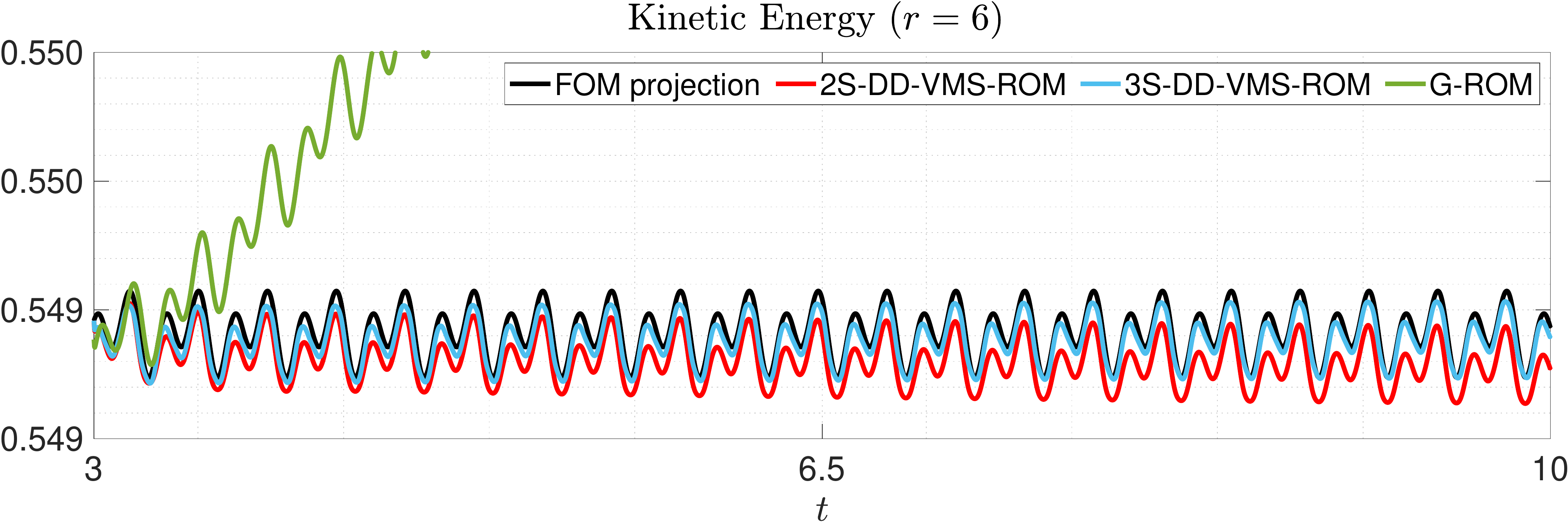}
		\includegraphics[width=0.9\textwidth,height=0.22\textwidth]{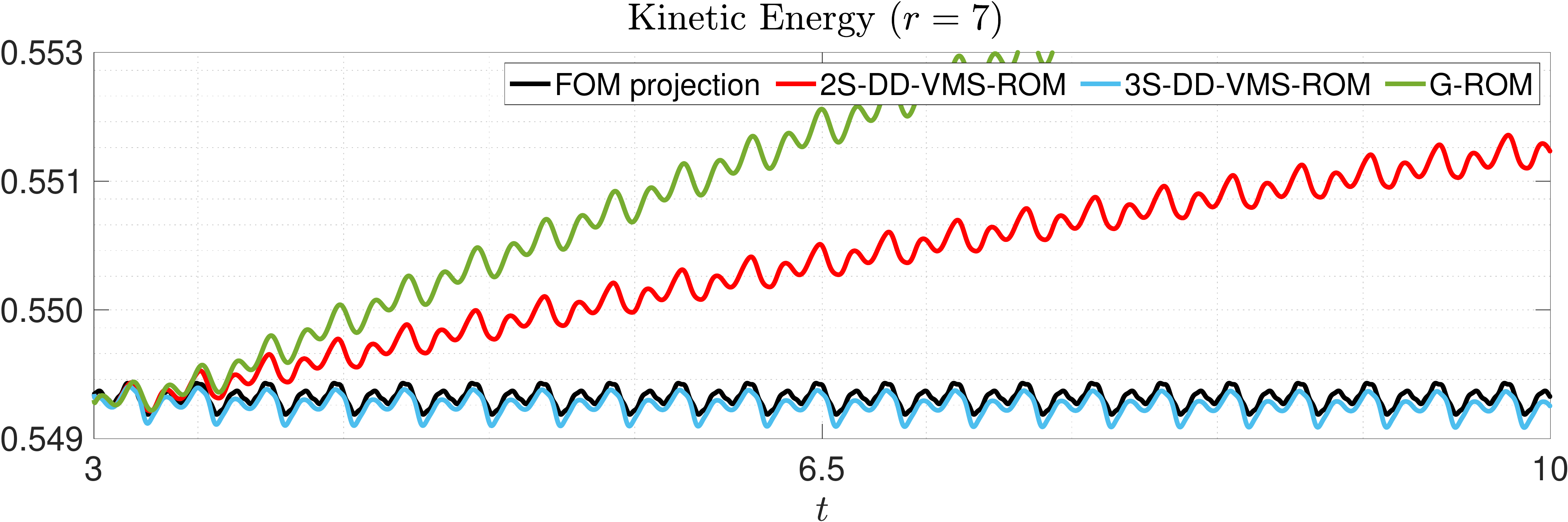}
	\end{center}
	\caption {
			Flow past a cylinder, $Re=100$, predictive regime.
			Time evolution of the kinetic energy for FOM projection, G-ROM, 2S-DD-VMS-ROM, and 3S-DD-VMS-ROM for different $r$ values. 
	\label{fig:nse-predictive-re100-ke}
	} 
\end{figure}

In Figure~\ref{fig:nse-predictive-re100-ke}, for $r=4, 6, 7$, we plot the time evolution of the kinetic energy of the 
FOM,
the G-ROM, the 2S-DD-VMS-ROM, and the new 3S-DD-VMS-ROM in the predictive regime. 
For all $r$ values, the G-ROM kinetic energy approximation is very inaccurate. 
In contrast, the new 3S-DD-VMS-ROM accurately approximates the exact 
FOM
kinetic energy for $r=4,6,7$. 
The 2S-DD-VMS-ROM kinetic energy approximation is accurate for $r=6$, but not for $r=4$ and, especially, for $r=7$.

\bigskip

The errors listed in Tables~\ref{table:nse-reconstructive-re100}--\ref{table:nse-predictive-re100} and the plots in Figures~\ref{fig:nse-reconstructive-re100-ke}--\ref{fig:nse-predictive-re100-ke} show that, in the reconstructive, cross-validation, and predictive regimes, the 3S-DD-VMS-ROM is consistently the most accurate ROM.
Furthermore, the 3S-DD-VMS-ROM is more accurate than the 2S-DD-VMS-ROM, especially in the predictive regime.



\subsubsection{Numerical Results for  $Re=500$} 
	\label{sec:numerical-results-nse-re500}

In this section, we present numerical results for the flow past a cylinder at $Re=500$.

In Table~\ref{table:nse-reconstructive-re500}, for different $r$ values, we list the average $L^{2}$ error~\eqref{eqn:l2-error} 
for the G-ROM, the 2S-DD-VMS-ROM, and the new 3S-DD-VMS-ROM in the reconstructive regime.
We also list 
the $r_1$ values for the 3S-DD-VMS-ROM.
These results show that, for all $r$ values, the 2S-DD-VMS-ROM  and 3S-DD-VMS-ROM are several times (sometimes more than one order of magnitude) more accurate than the standard G-ROM.
Furthermore, the 3S-DD-VMS-ROM is generally more accurate than the 2S-DD-VMS-ROM. 
For example, for $r=2$, the 3S-DD-VMS-ROM is almost twice more accurate as the 2S-DD-VMS-ROM.
\begin{table}[h!]
\small
	\begin{center}
		\smallskip
		\begin{tabular}{|c|c|c|c|c|c|c|c|}
			\hline
			\multicolumn{1}{|c|}{ $r$} &  
			\multicolumn{1}{c|}{ G-ROM } &  
			\multicolumn{1}{c|}{ 2S-DD-VMS-ROM } &  
			\multicolumn{2}{c|}{3S-DD-VMS-ROM}
						 \\  \cline{2-5}&{$\mathcal{E}(L^2)$ }  
						 &{$\mathcal{E}(L^2)$ }  
						 &{$r_1$} & { $\mathcal{E}(L^2)$ } 
									 \\
			 \hline  
$2$
&2.892e-01 
&7.029e-03 
&1&3.937e-03 
 	 \\
			 \hline    
$3$
& 3.344e-01 
& 8.138e-02 
&2& 7.517e-02  
 	 \\
			 \hline    
$4$
&3.478e-01 
&4.195e-03  
&3&4.145e-03 
 	 \\
			 \hline    
$5$
&3.795e-01 
&6.811e-02  
&2& 5.915e-02 
	 \\
			 \hline  
$6$
& 6.338e-02  
&3.864e-03   
&2&3.294e-03 
	 \\
			 \hline   
$7$
&5.738e-02  
&1.789e-02  
&2&1.563e-02 
	 \\
			 \hline  
$8$
&5.339e-02 
&5.734e-03 
&6&4.809e-03 
\\	\hline	 
		\end{tabular}
	\end{center}
	\caption {
			Flow past a cylinder, $Re=500$, reconstructive regime.
			Average $L^{2}$ errors for G-ROM, 2S-DD-VMS-ROM, and 3S-DD-VMS-ROM for different $r$ values. 
	\label{table:nse-reconstructive-re500}
	} 
\end{table}

In Table~\ref{table:nse-cross-re500}, for different $r$ values, we list the average $L^{2}$ error~\eqref{eqn:l2-error} 
for the G-ROM, the 2S-DD-VMS-ROM, and the new 3S-DD-VMS-ROM in the cross-validation regime.
We also list the $r_1$ values for the 3S-DD-VMS-ROM.
These results show that, for all $r$ values, the 2S-DD-VMS-ROM  and 3S-DD-VMS-ROM are several times (sometimes even two orders of magnitude) more accurate than the standard G-ROM.
Furthermore, the 3S-DD-VMS-ROM is generally more accurate than the 2S-DD-VMS-ROM. 
Specifically, for $r = 4, 5, 8$, the 3S-DD-VMS-ROM is almost twice more accurate than the 2S-DD-VMS-ROM.

\begin{table}[h!]
\small
	\begin{center}
		\smallskip
		\begin{tabular}{|c|c|c|c|c|c|c|c|}
						\hline
			\multicolumn{1}{|c|}{ $r$} &  
			\multicolumn{1}{c|}{ G-ROM } &  
			\multicolumn{1}{c|}{ 2S-DD-VMS-ROM } &  
			\multicolumn{2}{c|}{3S-DD-VMS-ROM}
						 \\  \cline{2-5}&{$\mathcal{E}(L^2)$ }  
						 &{$\mathcal{E}(L^2)$ }  
						 &{$r_1$} & { $\mathcal{E}(L^2)$ } 
									 \\
			 \hline  
$2$&
1.071e+00  
& 2.015e-02
&  1&1.501e-02
 	 \\
			 \hline    
$3$&
8.280e-01 
&   1.101e-01
& 2 &   8.428e-02
 	 \\
			 \hline    
$4$&
6.258e-01  
&   1.218e-02
& 3 &   4.648e-03  
 	 \\
			 \hline    
$5$&
6.440e-01
&    1.557e-01
& 3 &   7.329e-02 
	 \\
			 \hline  
$6$&
1.898e-01
&    5.733e-03
&  3&   4.056e-03
	 \\
			 \hline   
$7$& 
1.531e-01 
&     3.550e-02
& 2 &   2.033e-02 
	 \\
			 \hline  
$8$& 
1.678e-01
&     9.050e-03
& 1 &   5.480e-03
\\	\hline	 
		\end{tabular}
	\end{center}
	\caption{
			Flow past a cylinder, $Re=500$, cross-validation regime.
			Average $L^{2}$ errors for G-ROM, 2S-DD-VMS-ROM, and 3S-DD-VMS-ROM for different $r$ values. 
			 \label{table:nse-cross-re500}
			}
\end{table}

In  Table~\ref{table:nse-predictive-re500}, for different $r$ values, we list the average $L^{2}$ error~\eqref{eqn:l2-error} for the G-ROM, the 2S-DD-VMS-ROM, and the new 3S-DD-VMS-ROM in the predictive regime.
We also list 
the $r_1$ values for the 3S-DD-VMS-ROM.
These results show that, for all $r$ values, the 2S-DD-VMS-ROM  and 3S-DD-VMS-ROM are several times (sometimes even more than one order of magnitude) more accurate than the standard G-ROM.
More  importantly,  for all $r$ values (but  especially for large $r$ values), the 3S-DD-VMS-ROM is significantly more accurate than the 2S-DD-VMS-ROM:
For example, for  $r=5, 6, 7$, and $8$, the 3S-DD-VMS-ROM is {\it more than twice more accurate} than the 2S-DD-VMS-ROM.

\begin{table}[h!]
\small
	\begin{center}
		\smallskip
		\begin{tabular}{|c|c|c|c|c|c|c|c|}
			\hline
			\multicolumn{1}{|c|}{ $r$} &  
			\multicolumn{1}{c|}{ G-ROM } &  
			\multicolumn{1}{c|}{ 2S-DD-VMS-ROM } &  
			\multicolumn{2}{c|}{3S-DD-VMS-ROM}
						 \\  \cline{2-5}&{$\mathcal{E}(L^2)$ }  
						 &{$\mathcal{E}(L^2)$ }  
						 &{$r_1$} & { $\mathcal{E}(L^2)$ } 
									 \\
			 \hline  
$2$
&7.351e-01 
&1.004e-01 
&1& 1.004e-01 
 	 \\
			 \hline    
$3$
&7.088e-01  
&8.838e-02  
& 2&  8.497e-02 
 	 \\
			 \hline    
$4$
&5.871e-01 
&8.785e-03 
&1&   8.785e-03 
 	 \\
			 \hline    
$5$
&6.231e-01 
&9.735e-02  
&2&   3.640e-02 
	 \\
			 \hline  
$6$
&1.293e-01  
&2.288e-02  
&  4&   9.051e-03  
	 \\
			 \hline   
$7$
&1.069e-01  
&2.816e-02  
&6&  1.480e-02 
	 \\
			 \hline  
$8$
&1.130e-01 
&1.402e-02 
&6&   5.544e-03 
\\	\hline	 
		\end{tabular}
	\end{center}
	\caption {
			Flow past a cylinder, $Re=500$, predictive regime.
			Average $L^{2}$ errors for G-ROM, 2S-DD-VMS-ROM, and 3S-DD-VMS-ROM for different $r$ values. 
	\label{table:nse-predictive-re500}
	} 
	\end{table}

\bigskip

\begin{figure}[H]
	\begin{center}
		\includegraphics[width=0.9\textwidth,height=0.22\textwidth]{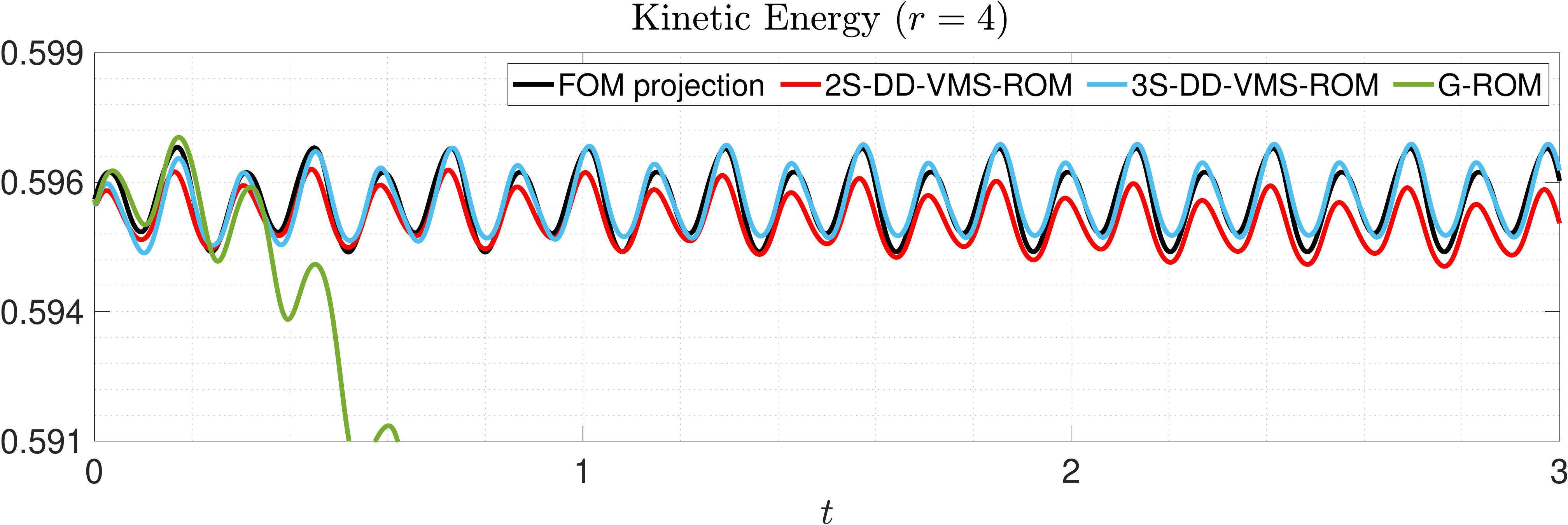}
		\includegraphics[width=0.9\textwidth,height=0.22\textwidth]{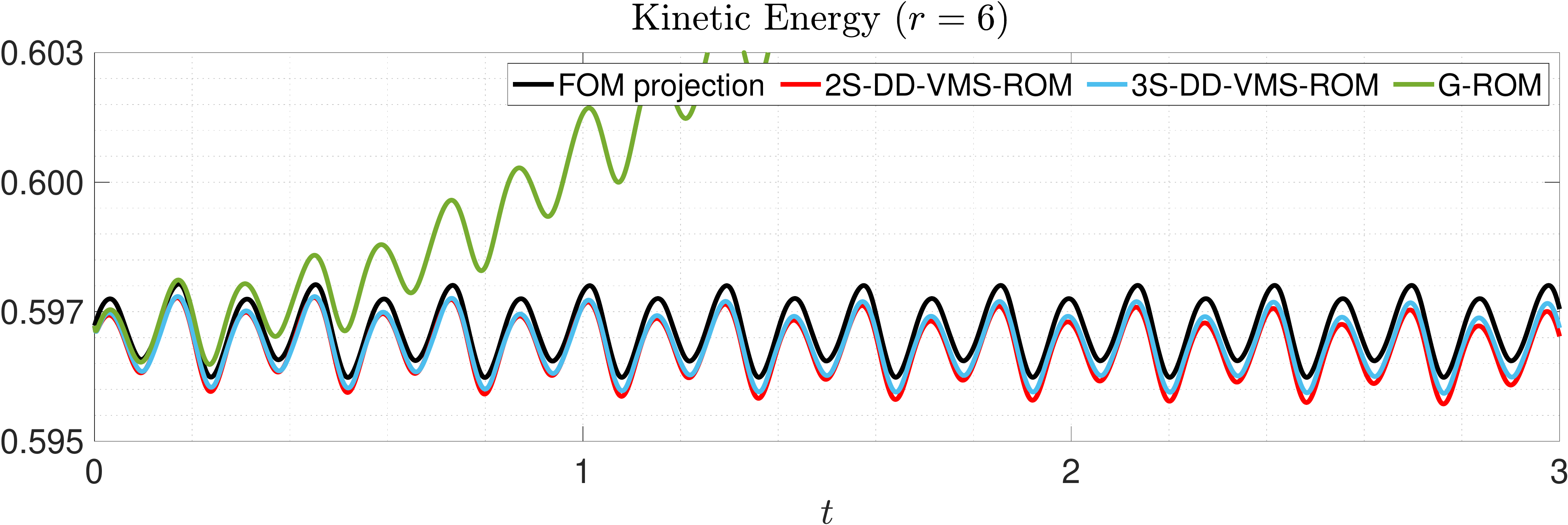}
		\includegraphics[width=0.9\textwidth,height=0.22\textwidth]{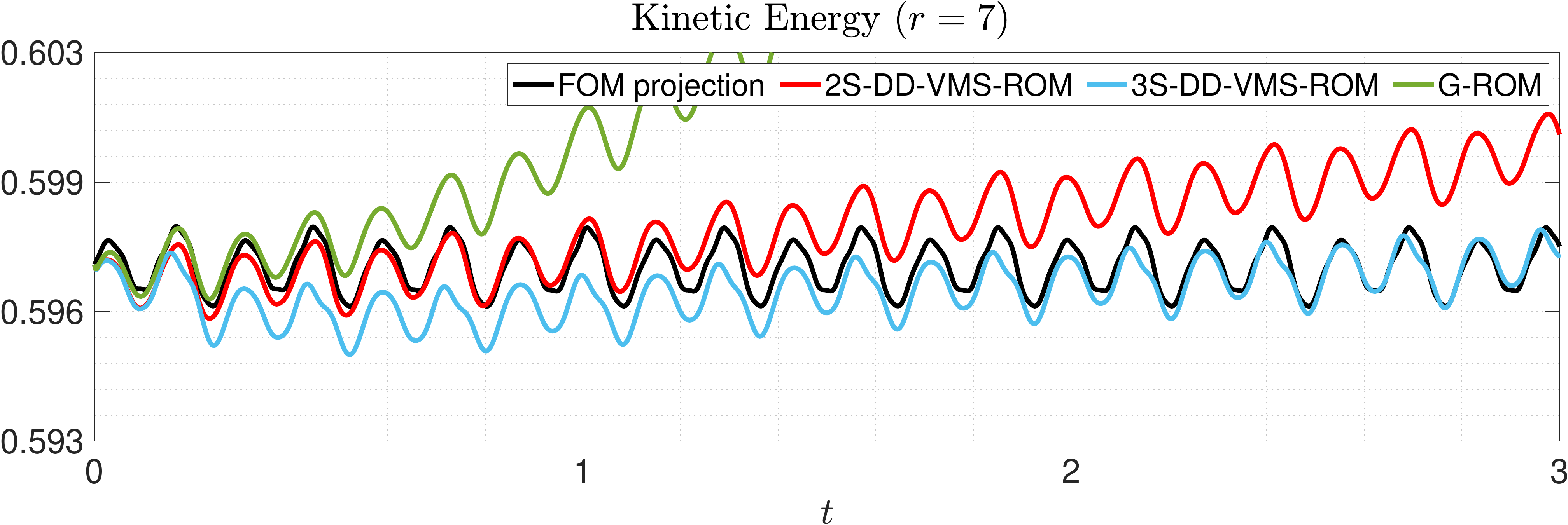}
	\end{center}
	\caption {
			Flow past a cylinder, $Re=500$, reconstructive regime.
			Time evolution of the kinetic energy for FOM projection, G-ROM, 2S-DD-VMS-ROM, and 3S-DD-VMS-ROM for different $r$ values. 
	\label{fig:nse-reconstructive-re500-ke}
	} 
\end{figure}

In Figure~\ref{fig:nse-reconstructive-re500-ke}, for $r=4, 6, 7$, we plot the time evolution of the kinetic energy of the 
FOM,
the G-ROM, the 2S-DD-VMS-ROM, and the new 3S-DD-VMS-ROM in the reconstructive regime. 
For all $r$ values, the G-ROM kinetic energy approximation is very inaccurate. 
In contrast, the new 3S-DD-VMS-ROM accurately approximates the exact
FOM
kinetic energy for $r=4,6,7$. 
The 2S-DD-VMS-ROM kinetic energy approximation is accurate for $r=4$ and $r=6$, but not for $r=7$.

\begin{figure}[H]
	\begin{center}
		\includegraphics[width=0.9\textwidth,height=0.22\textwidth]{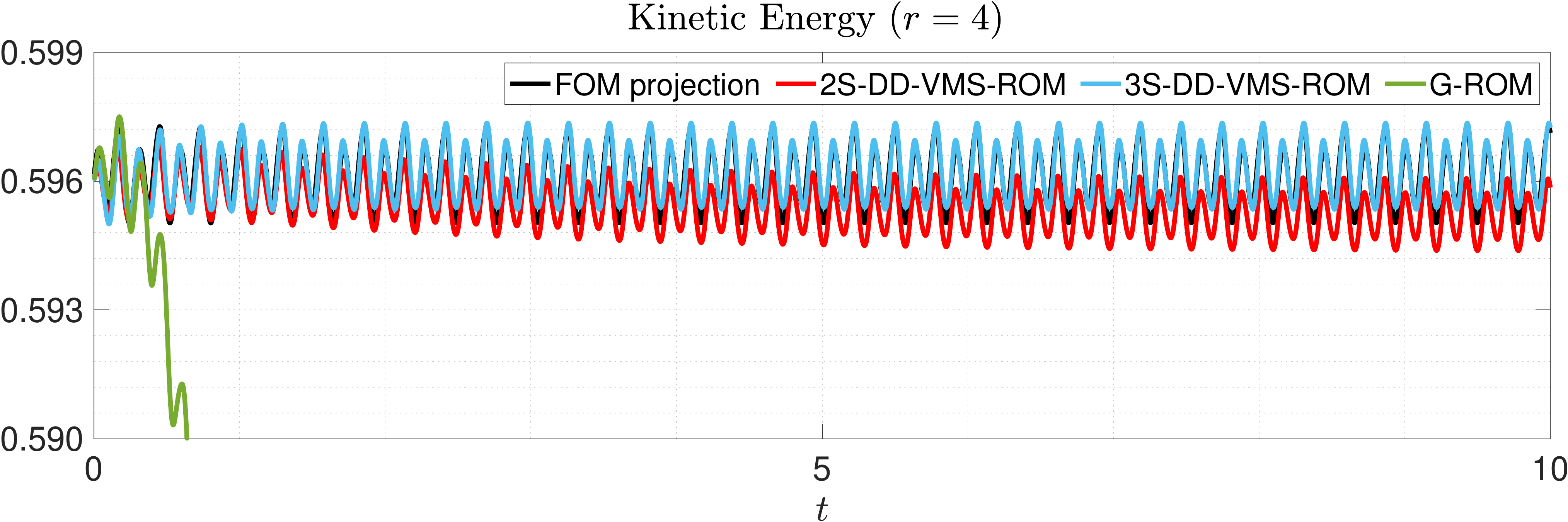}
		\includegraphics[width=0.9\textwidth,height=0.22\textwidth]{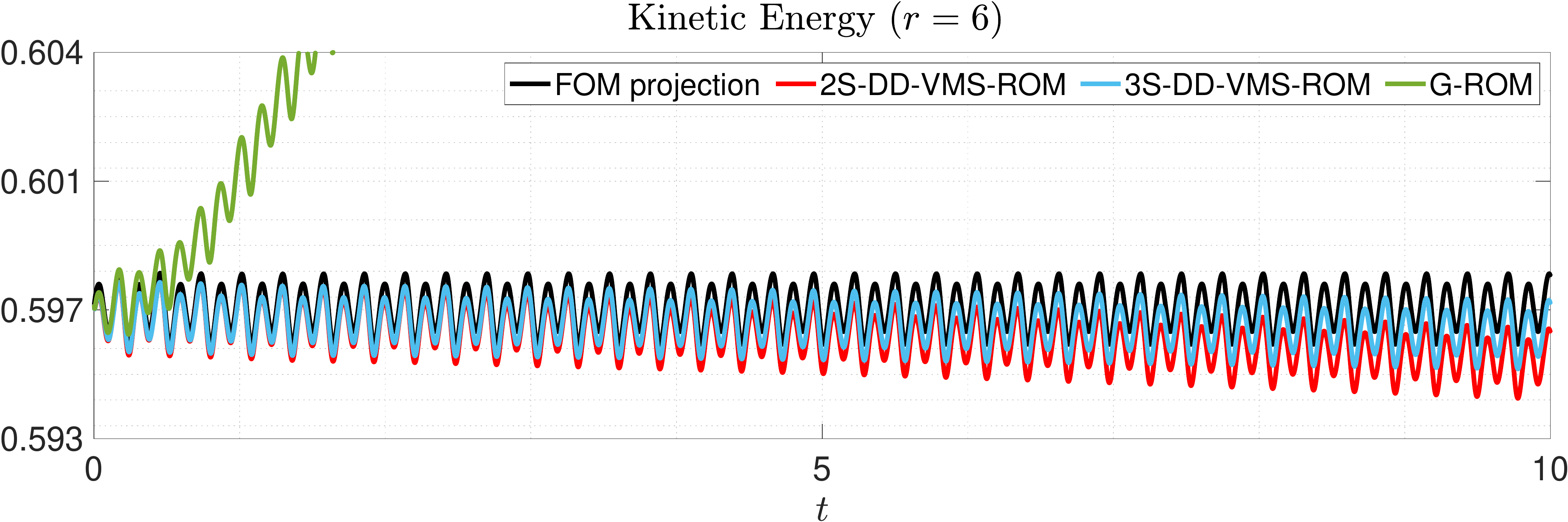}
		\includegraphics[width=0.9\textwidth,height=0.22\textwidth]{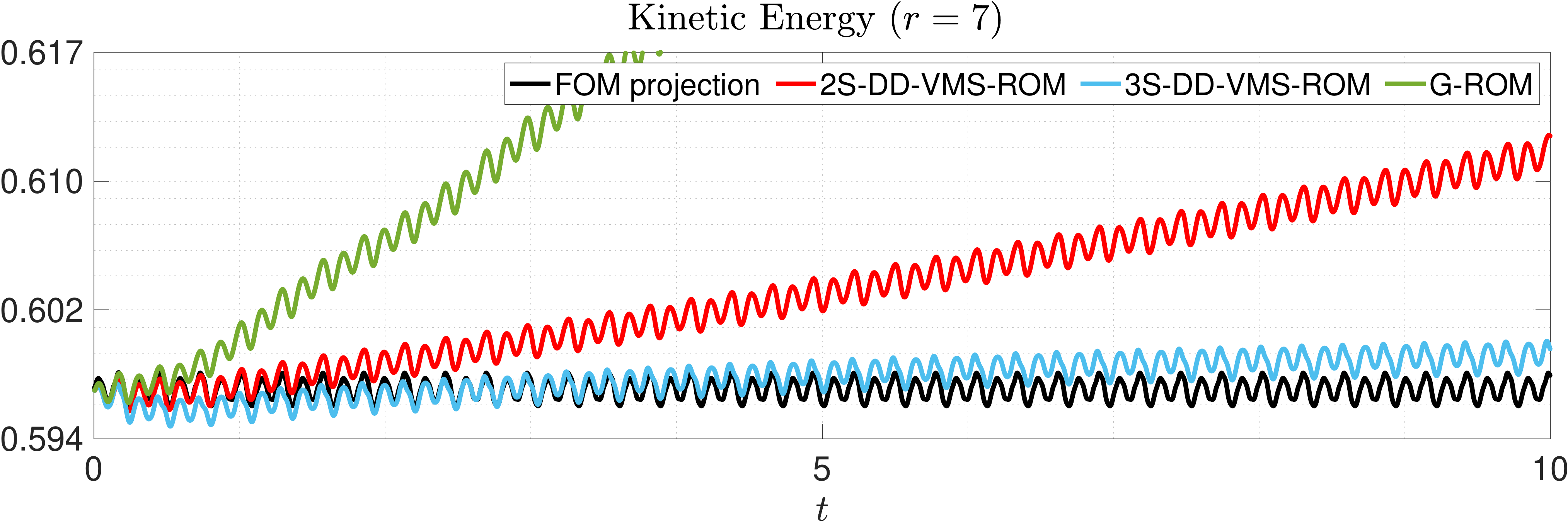}
	\end{center}
	\caption {
			Flow past a cylinder, $Re=500$, cross-validation regime.
			Time evolution of the kinetic energy for FOM projection, G-ROM, 2S-DD-VMS-ROM, and 3S-DD-VMS-ROM for different $r$ values. 
	\label{fig:nse-cross-re500-ke}
	} 
\end{figure}

In Figure~\ref{fig:nse-cross-re500-ke}, for $r=4, 6, 7$, we plot the time evolution of the kinetic energy of the 
FOM,
the G-ROM, the 2S-DD-VMS-ROM, and the new 3S-DD-VMS-ROM in the cross-validation regime. 
For all $r$ values, the G-ROM kinetic energy approximation is very inaccurate. 
In contrast, the new 3S-DD-VMS-ROM accurately approximates the exact
FOM
kinetic energy for $r=4,6,7$. 
The 2S-DD-VMS-ROM kinetic energy approximation is accurate for $r=4$ and $r=6$, but not for $r=7$.

\begin{figure}[H]

	\begin{center}
		\includegraphics[width=0.9\textwidth,height=0.22\textwidth]{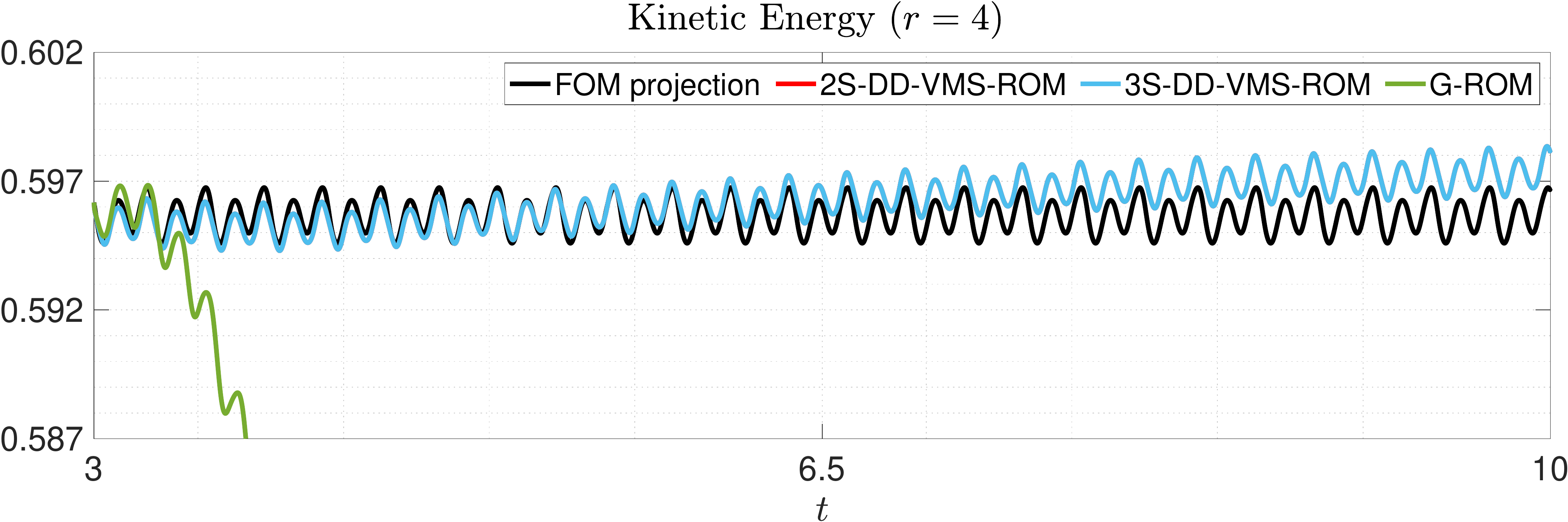}
		\includegraphics[width=0.9\textwidth,height=0.22\textwidth]{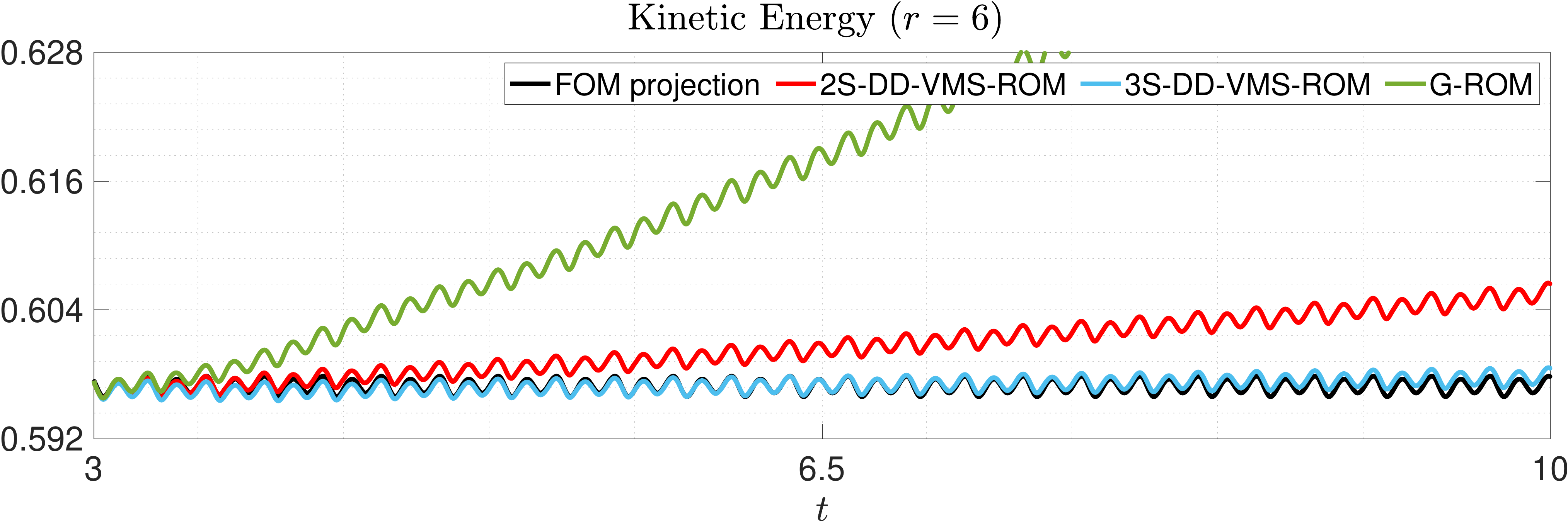}
		\includegraphics[width=0.9\textwidth,height=0.22\textwidth]{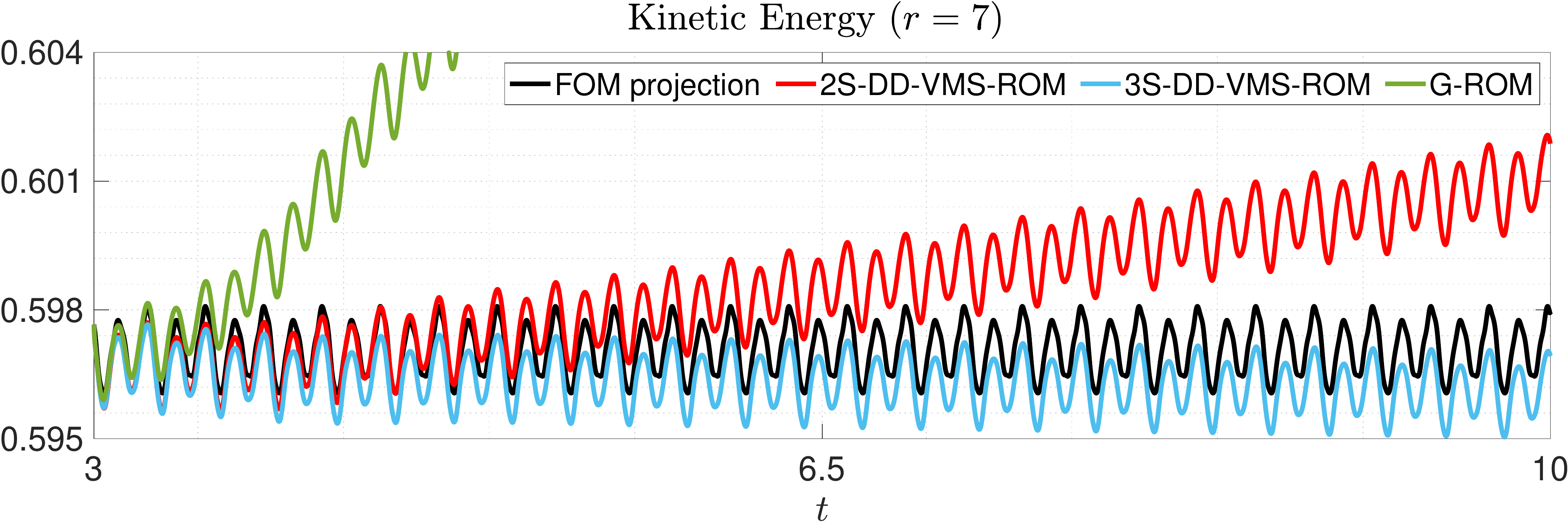}

	\end{center}
	\caption {
			Flow past a cylinder, $Re=500$, predictive regime.
			Time evolution of the kinetic energy for FOM projection, G-ROM, 2S-DD-VMS-ROM, and 3S-DD-VMS-ROM for different $r$ values. 
	\label{fig:nse-predictive-re500-ke}
	}
\end{figure}

In Figure~\ref{fig:nse-predictive-re500-ke}, for $r=4, 6, 7$, we plot the time evolution of the kinetic energy of the 
FOM,
the G-ROM, the 2S-DD-VMS-ROM, and the new 3S-DD-VMS-ROM in the predictive regime. 
For all $r$ values, the G-ROM kinetic energy approximation is very inaccurate. 
In contrast, the new 3S-DD-VMS-ROM accurately approximates the 
FOM
kinetic energy for $r=6$ and $r=7$. 
For $r=6$ and $r=7$, the 2S-DD-VMS-ROM kinetic energy approximation is less accurate than the 3S-DD-VMS-ROM kinetic energy approximation but more accurate than the G-ROM kinetic energy approximation.
For $r=4$, both the 3S-DD-VMS-ROM and the 2S-DD-VMS-ROM kinetic energy approximations are accurate. 

\bigskip


The errors listed in Tables~\ref{table:nse-reconstructive-re500}--\ref{table:nse-predictive-re500} and the plots in Figures~\ref{fig:nse-reconstructive-re500-ke}--\ref{fig:nse-predictive-re500-ke} show that, in the reconstructive,  cross-validation, and predictive regimes, the 3S-DD-VMS-ROM is consistently the most accurate ROM.
Furthermore, the 3S-DD-VMS-ROM is more accurate than the 2S-DD-VMS-ROM, especially in the predictive regime.



\subsubsection{Numerical Results for  $Re=1000$} 
	\label{sec:numerical-results-nse-re1000}

In this section, we present numerical results for the flow past a cylinder at $Re=1000$.

In Table~\ref{table:nse-reconstructive-re1000}, for different $r$ values, we list the average $L^{2}$ error~\eqref{eqn:l2-error} 
for the G-ROM, the 2S-DD-VMS-ROM, and the new 3S-DD-VMS-ROM in the reconstructive regime.
We also list 
the $r_1$ values for the 3S-DD-VMS-ROM.
These results show that, for all $r$ values, the 2S-DD-VMS-ROM  and 3S-DD-VMS-ROM are several times (sometimes even more than one order of magnitude) more accurate than the standard G-ROM.
Furthermore, the 3S-DD-VMS-ROM is generally more accurate than the 2S-DD-VMS-ROM.
For example, for $r=5$ and $r=8$, the 3S-DD-VMS-ROM is {\it almost twice more accurate} than the 2S-DD-VMS-ROM.

\begin{table}[h!]
\small
	\begin{center}
		\smallskip
		\begin{tabular}{|c|c|c|c|c|c|c|c|}
			\hline
			\multicolumn{1}{|c|}{ $r$} &  
			\multicolumn{1}{c|}{ G-ROM } &  
			\multicolumn{1}{c|}{ 2S-DD-VMS-ROM } &  
			\multicolumn{2}{c|}{3S-DD-VMS-ROM}
						 \\  \cline{2-5}&{$\mathcal{E}(L^2)$ }  
						 &{$\mathcal{E}(L^2)$ }  
						 &{$r_1$} & { $\mathcal{E}(L^2)$ } 
									 \\
			 \hline   
$2$
&4.937e-01 
&6.704e-03 
&1&6.692e-03 
 	 \\
			 \hline    
$3$
&5.112e-01 
&6.804e-02  
&1& 6.794e-02 
 	 \\
			 \hline    
$4$
&5.980e-01  
&1.287e-02  
&2&   9.869e-03 
 	 \\
			 \hline    
$5$
&6.579e-01  
&1.794e-01  
&3&   9.184e-02  
	 \\
			 \hline  
$6$
&1.503e-01  
& 1.086e-02  
&4&    8.210e-03  
	 \\
			 \hline   
$7$
&1.365e-01   
&  2.848e-02 
&5&   2.235e-02 
	 \\
			 \hline  
$8$
& 7.076e-02  
& 7.550e-03  
&4&    4.836e-03 
\\	\hline	 
		\end{tabular}
	\end{center}
	\caption {
			Flow past a cylinder, $Re=1000$, reconstructive regime.
			Average $L^{2}$ errors for G-ROM, 2S-DD-VMS-ROM, and 3S-DD-VMS-ROM for different $r$ values. 
	\label{table:nse-reconstructive-re1000}
	} 
\end{table}

In Table~\ref{table:cross-nse-re1000}, for different $r$ values, we list the average $L^{2}$ error~\eqref{eqn:l2-error} 
for the G-ROM, the 2S-DD-VMS-ROM, and the new 3S-DD-VMS-ROM in the cross-validation regime. These results 
show
that, for all $r$ values, the 2S-DD-VMS-ROM and 3S-DD-VMS-ROM are several times (sometimes even two orders of magnitude) more accurate than the standard G-ROM. Furthermore, the 3S-DD-VMS-ROM is generally more accurate than the 2S-DD-VMS-ROM, especially for large $r$ values. In particular, for $r= 5$, the 3S-DD-VMS-ROM is {\it almost five times more accurate} than the 2S-DD-VMS-ROM.  
\begin{table}[h!]
\small
	\begin{center}
		\smallskip
		\begin{tabular}{|c|c|c|c|c|c|c|c|}
			\hline
			\multicolumn{1}{|c|}{ $r$} &  
			\multicolumn{1}{c|}{ G-ROM } &  
			\multicolumn{1}{c|}{ 2S-DD-VMS-ROM } &  
			\multicolumn{2}{c|}{3S-DD-VMS-ROM}
						 \\  \cline{2-5}&{$\mathcal{E}(L^2)$ }  
						 &{$\mathcal{E}(L^2)$ }  
						 &{$r_1$} & { $\mathcal{E}(L^2)$ } 
									 \\
			 \hline   
$2$
&1.509e+00
&1.504e-02 
&1&1.503e-02 
 	 \\
			 \hline    
$3$
& 8.595e-01
&8.024e-02  
&1&8.024e-02  
 	 \\
			 \hline    
$4$
&6.583e-01
&2.538e-02  
&2&1.503e-02
 	 \\
			 \hline    
$5$
& 7.095e-01
&5.156e-01 
&3&1.026e-01
	 \\
			 \hline  
$6$
& 5.562e-01
& 3.132e-02  
&4&1.018e-02 
	 \\
			 \hline   
$7$
&4.760e-01
&6.482e-02 
&6&3.505e-02 
	 \\
			 \hline  
$8$
&2.692e-01
&1.691e-02 
&5&5.791e-03  
\\	\hline	 
		\end{tabular}
	\end{center}
	\caption{
			Flow past a cylinder, $Re=1000$, cross-validation regime.
			Average $L^{2}$ errors for G-ROM, 2S-DD-VMS-ROM, and 3S-DD-VMS-ROM for different $r$ values. 
	\label{table:cross-nse-re1000}
	} 
	\end{table}
	
In Table~\ref{table:nse-predictive-re1000}, for different $r$ values, we list the average $L^{2}$ error~\eqref{eqn:l2-error} for the G-ROM, the 2S-DD-VMS-ROM, and the new 3S-DD-VMS-ROM in the predictive regime.
We also list 
the $r_1$ values for the 3S-DD-VMS-ROM.
These results show that, for all $r$ values, the 2S-DD-VMS-ROM  and 3S-DD-VMS-ROM are significantly (sometimes several times) more accurate than the standard G-ROM.
More  importantly,  for all $r$ values (but  especially for large $r$ values), the 3S-DD-VMS-ROM is significantly more accurate than the 2S-DD-VMS-ROM:
For example, for  $r=6$, the 3S-DD-VMS-ROM is {\it more than five times more accurate} than the 2S-DD-VMS-ROM.

\begin{table}[h!]
\small
	\begin{center}
		\smallskip
		\begin{tabular}{|c|c|c|c|c|c|c|c|}
			\hline
			\multicolumn{1}{|c|}{ $r$} &  
			\multicolumn{1}{c|}{ G-ROM } &  
			\multicolumn{1}{c|}{ 2S-DD-VMS-ROM } &  
			\multicolumn{2}{c|}{3S-DD-VMS-ROM}
						 \\  \cline{2-5}&{$\mathcal{E}(L^2)$ }  
						 &{$\mathcal{E}(L^2)$ }  
						 &{$r_1$} & { $\mathcal{E}(L^2)$ } 
									 \\
			 \hline  
$2$
&1.146e+00 
& 3.857e-01 
&1&  2.860e-01 
 	 \\
			 \hline    
$3$
&9.217e-01 
&4.522e-01 
&1&1.357e-01 
 	 \\
			 \hline    
$4$
&7.207e-01 
&1.679e-01 
&2&7.070e-02  
 	 \\
			 \hline    
$5$
&7.281e-01 
&5.620e-01 
&3&2.331e-01  
	 \\
			 \hline  
$6$
&3.545e-01 
&2.279e-01 
&2&3.733e-02  
	 \\
			 \hline   
$7$
&3.027e-01  
&2.273e-01  
&3&7.922e-02 
	 \\
			 \hline  
$8$
&1.587e-01  
&5.849e-02  
&3&4.394e-02 
\\	\hline	 
		\end{tabular}
	\end{center}
	\caption {
			Flow past a cylinder, $Re=1000$, predictive regime.
			Average $L^{2}$ errors for G-ROM, 2S-DD-VMS-ROM, and 3S-DD-VMS-ROM for different $r$ values. 
	\label{table:nse-predictive-re1000}
	} 
	\end{table}

\bigskip

\begin{figure}[H]
	\begin{center}
		\includegraphics[width=0.9\textwidth,height=0.22\textwidth]{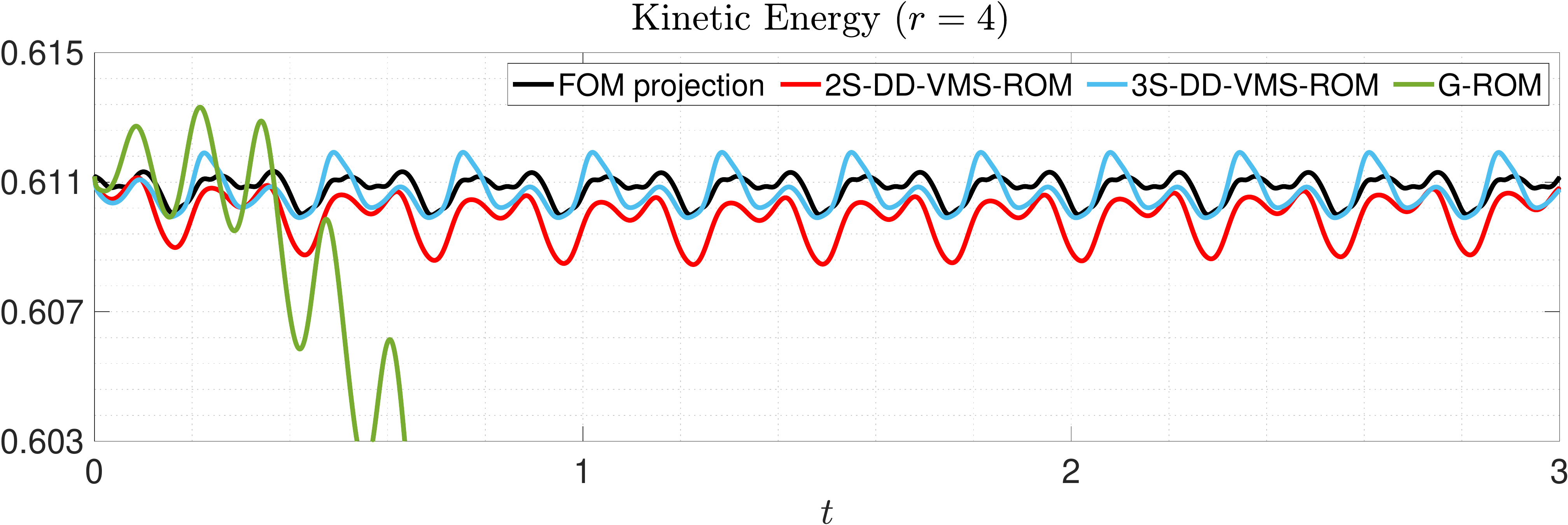}
		\includegraphics[width=0.9\textwidth,height=0.22\textwidth]{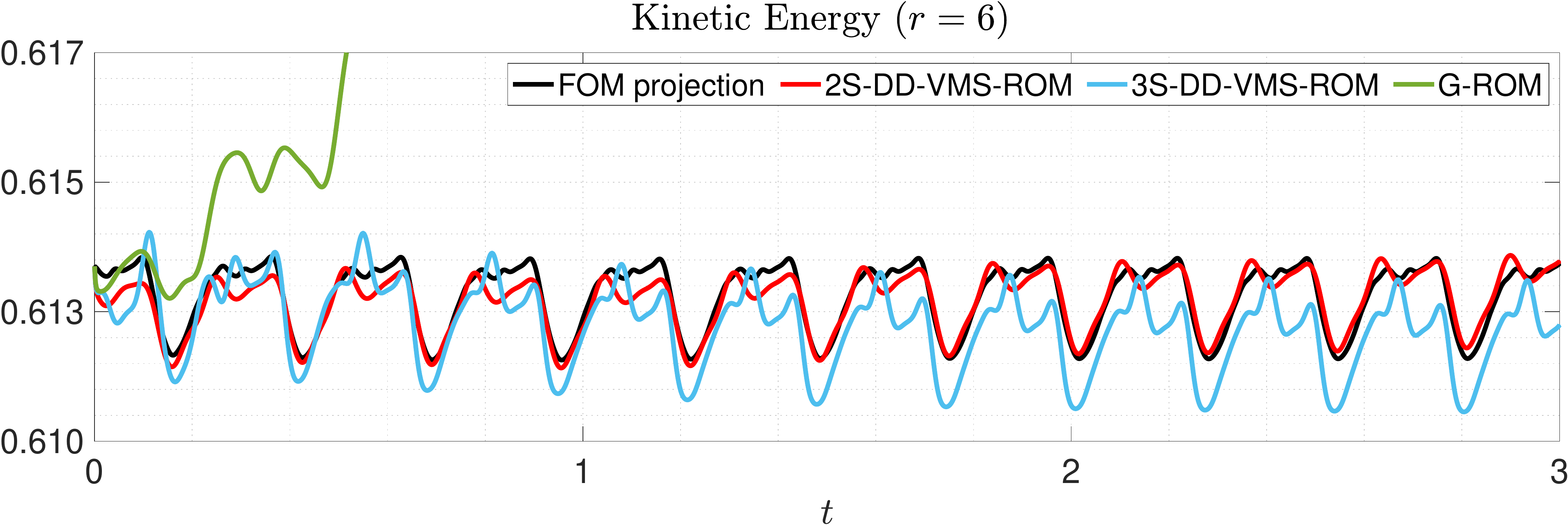}
		\includegraphics[width=0.9\textwidth,height=0.22\textwidth]{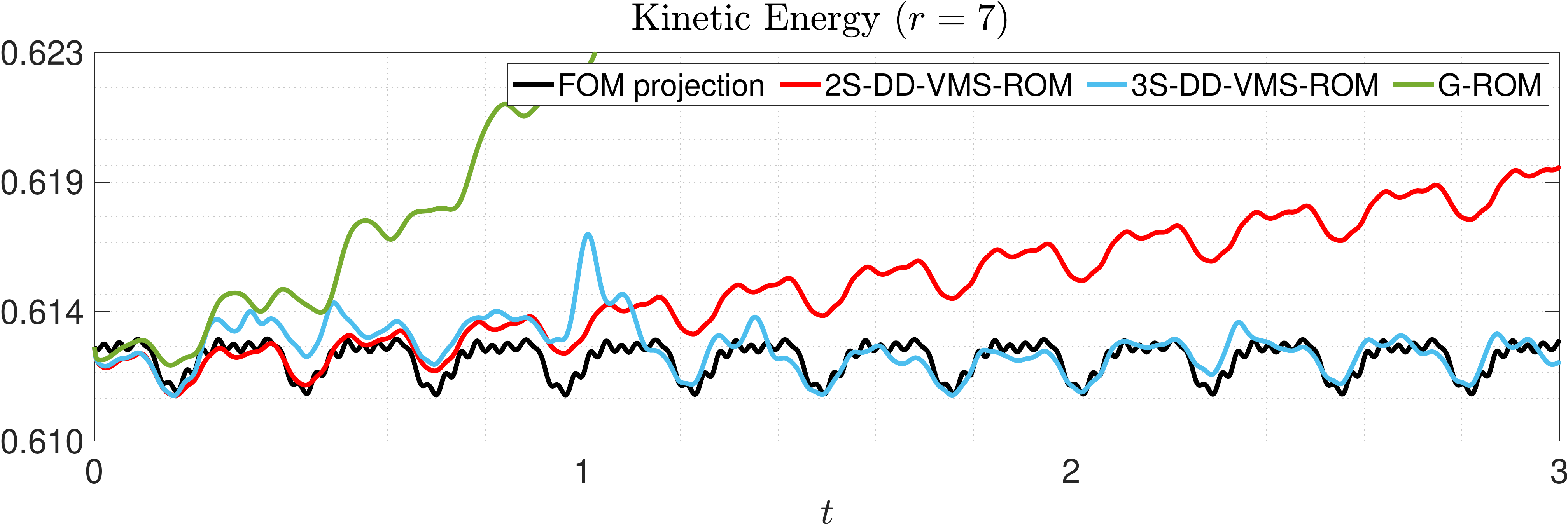}
	\end{center}
	\caption {
			Flow past a cylinder, $Re=1000$, reconstructive regime.
			Time evolution of the kinetic energy for FOM projection, G-ROM, 2S-DD-VMS-ROM, and 3S-DD-VMS-ROM for different $r$ values. 
	\label{fig:nse-reconstructive-re1000-ke}
	} 
\end{figure}

 In Figure~\ref{fig:nse-reconstructive-re1000-ke}, for $r=4, 6, 7$, we plot the time evolution of the kinetic energy of the 
FOM,
 the G-ROM, the 2S-DD-VMS-ROM, and the new 3S-DD-VMS-ROM in the reconstructive regime. 
 For all cases, the evolution of the G-ROM kinetic energy is very inaccurate.
 In contrast, for $r=4$ and $r=6$, both the 3S-DD-VMS-ROM and the 2S-DD-VMS-ROM successfully reproduce the 
FOM
 kinetic energy. 
 For $r=7$,  3S-DD-VMS-ROM kinetic energy yields small oscillations for $0 \leq t \leq 1$, but it quickly converges to the 
FOM
 kinetic energy after $t>1$. 
 On the other hand, the 2S-DD-VMS-ROM kinetic energy approximation is not accurate.

\begin{figure}[H]
	\begin{center}
		\includegraphics[width=0.9\textwidth,height=0.22\textwidth]{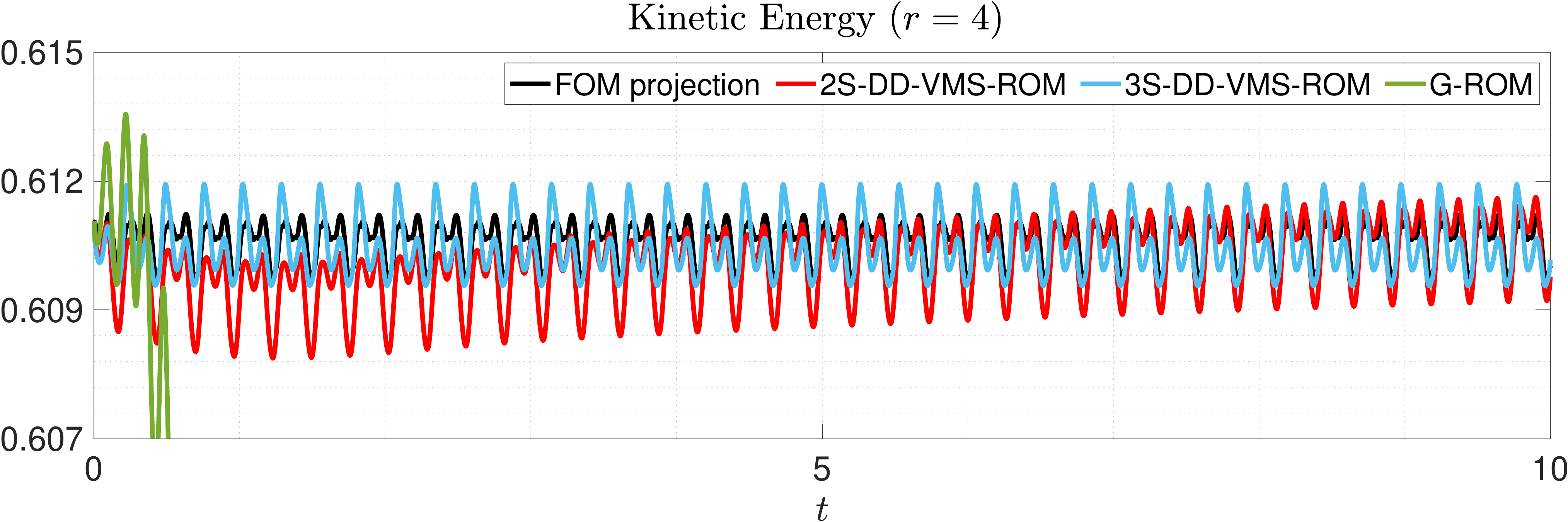}
		\includegraphics[width=0.9\textwidth,height=0.22\textwidth]{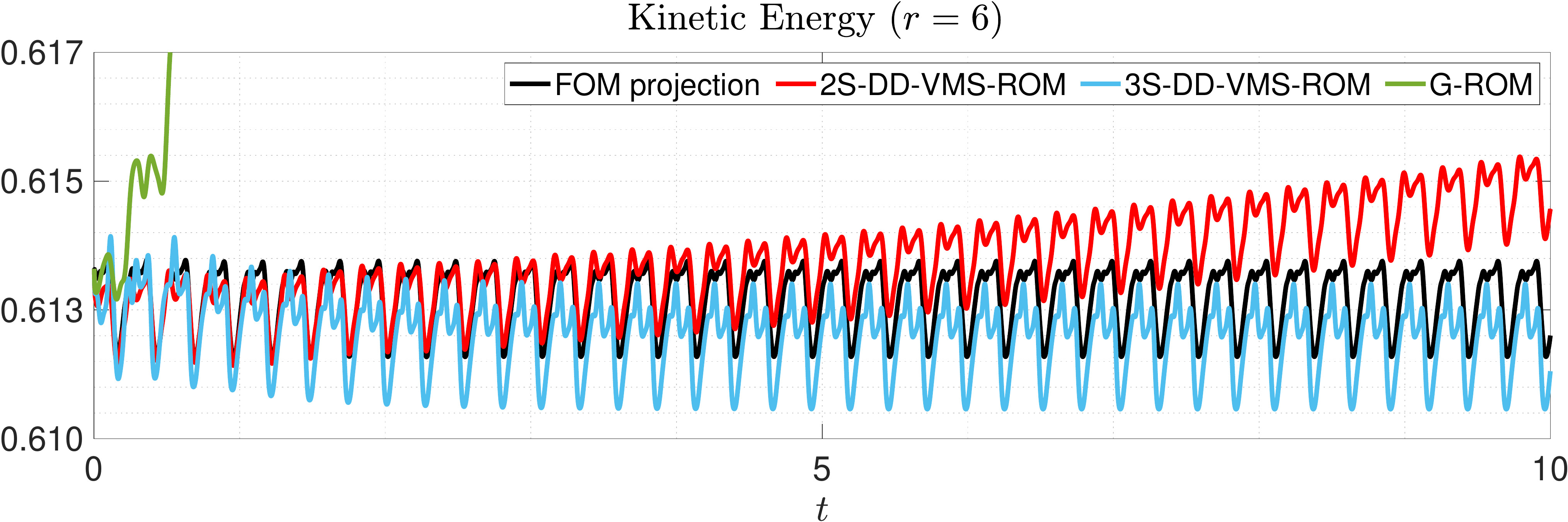}
		\includegraphics[width=0.9\textwidth,height=0.22\textwidth]{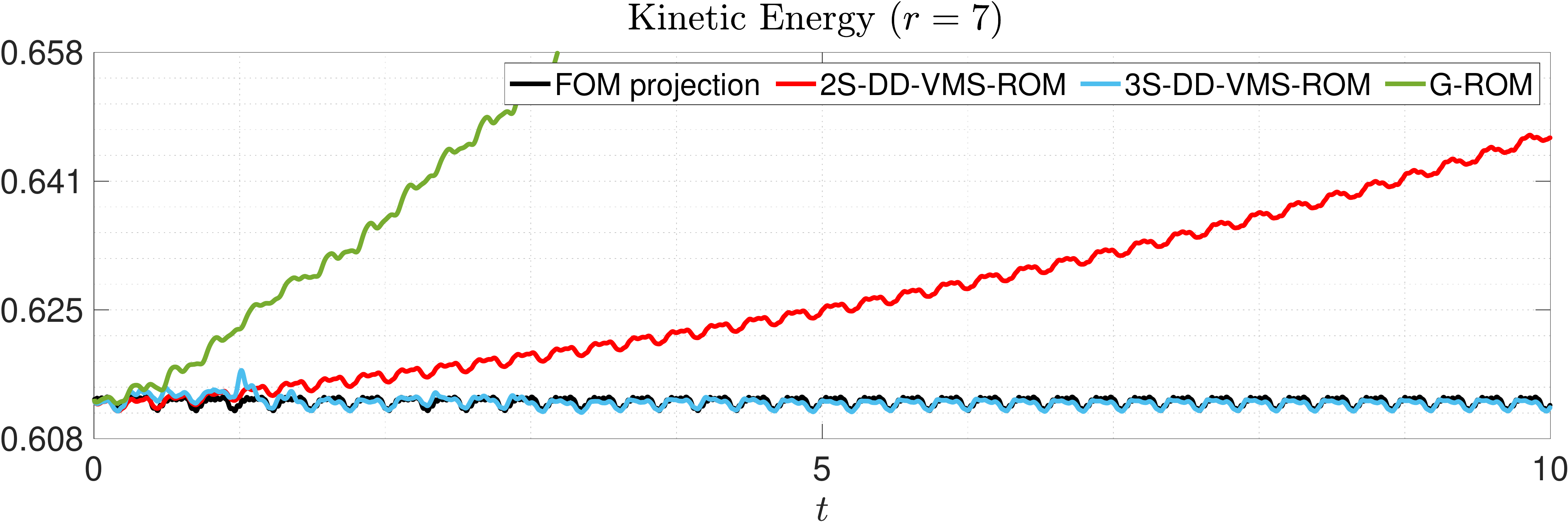}
	\end{center}
	\caption {
			Flow past a cylinder, $Re=1000$, cross-validation regime.
			Time evolution of the kinetic energy for FOM projection, G-ROM, 2S-DD-VMS-ROM, and 3S-DD-VMS-ROM for different $r$ values. 
	\label{fig:nse-cross-re1000-ke}
	} 
\end{figure}

 In Figure~\ref{fig:nse-cross-re1000-ke}, for $r=4, 6, 7$, we plot the time evolution of the kinetic energy of the 
FOM,
 the G-ROM, the 2S-DD-VMS-ROM, and the new 3S-DD-VMS-ROM in the cross-validation regime. 
 For all cases, the 
 evolution
 of the G-ROM kinetic energy is very inaccurate. 
 In contrast, for all cases, the 3S-DD-VMS-ROM successfully reproduces the exact 
FOM
 kinetic energy. 
 The 2S-DD-VMS-ROM kinetic energy is accurate for $r=4$, but not for $r=6$ and, especially, for $r=7$.

\begin{figure}[H]
	\begin{center}
		\includegraphics[width=0.9\textwidth,height=0.22\textwidth]{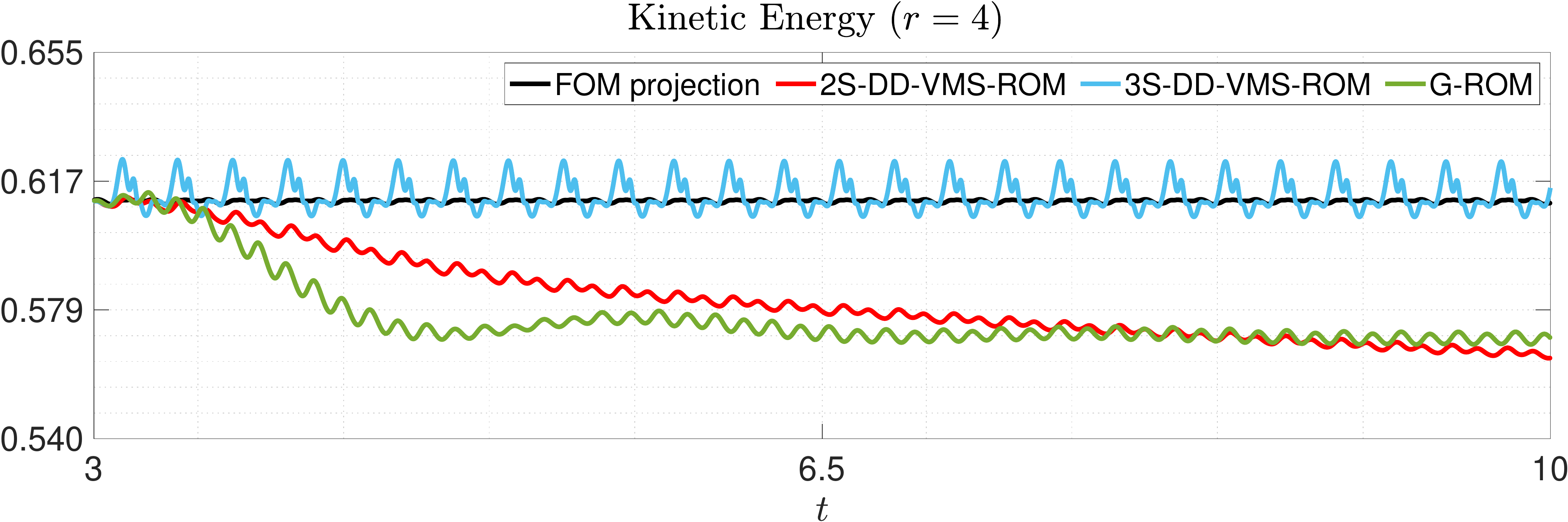}
		\includegraphics[width=0.9\textwidth,height=0.22\textwidth]{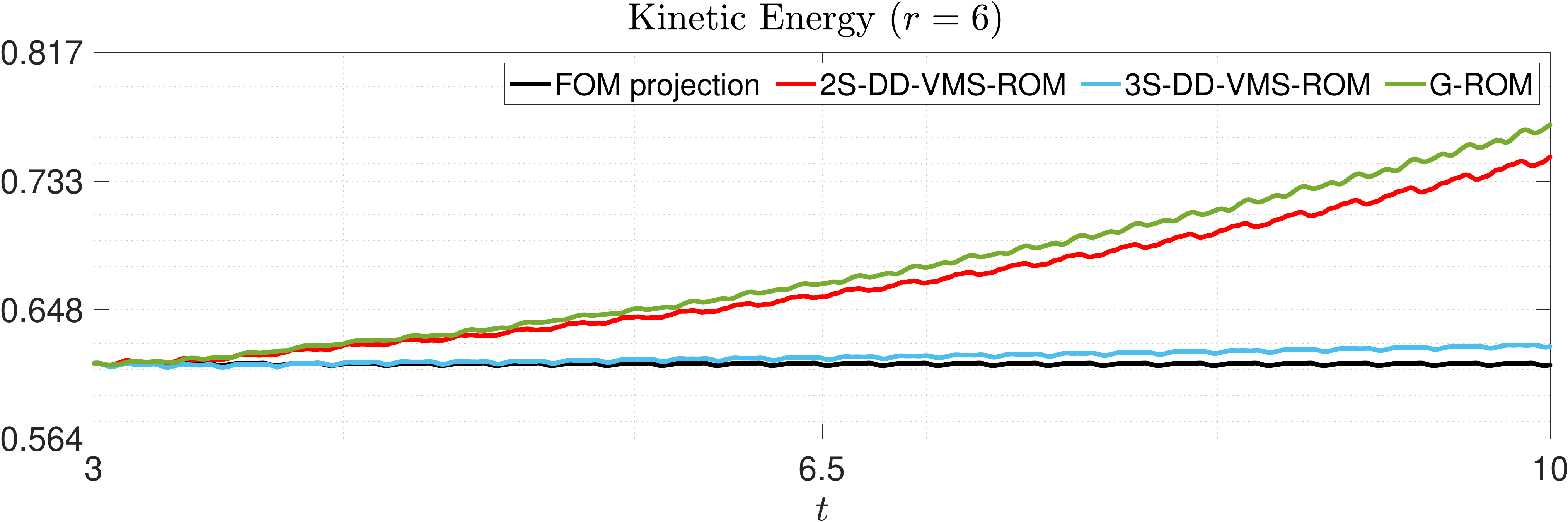}
		\includegraphics[width=0.9\textwidth,height=0.22\textwidth]{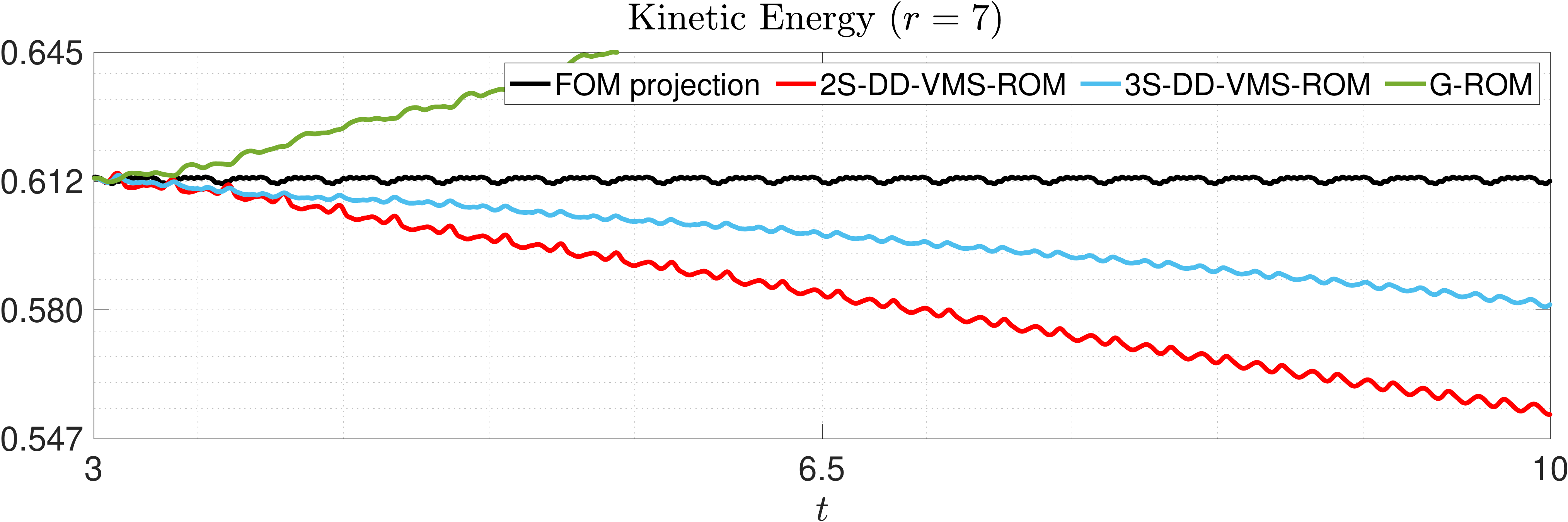}
	\end{center}
	\caption {
			Flow past a cylinder, $Re=1000$, predictive regime.
			Time evolution of the kinetic energy for FOM projection, G-ROM, 2S-DD-VMS-ROM, and 3S-DD-VMS-ROM for different $r$ values. 
	\label{fig:nse-predictive-re1000-ke}
	}
\end{figure}

In Figure~\ref{fig:nse-predictive-re1000-ke}, for three $r$ values, we plot the time evolution of the kinetic energy of the  G-ROM, the 2S-DD-VMS-ROM, and the new 3S-DD-VMS-ROM in the predictive regime.
For all cases, the evolution of the G-ROM kinetic energy is very inaccurate. 
For $r=4$, the 3S-DD-VMS-ROM kinetic energy approximation is accurate, whereas the 2S-DD-VMS-ROM and the G-ROM kinetic energies are inaccurate. 
For $r=6$ and $r=7$, although the 3S-DD-VMS-ROM kinetic energy approximations are not as accurate, they are still much more accurate than the 2S-DD-VMS-ROM and, especially, the G-ROM kinetic energy approximations.

\bigskip


The errors listed in Tables~\ref{table:nse-reconstructive-re1000}--\ref{table:nse-predictive-re1000} and the plots in Figures~\ref{fig:nse-reconstructive-re1000-ke}--\ref{fig:nse-predictive-re1000-ke} show that, in the reconstructive,  cross-validation, and predictive regimes, the 3S-DD-VMS-ROM is consistently the most accurate ROM.
Furthermore, the 3S-DD-VMS-ROM is more accurate than the 2S-DD-VMS-ROM, especially in the predictive regime.

\clearpage
{
\subsection{Quasi-Geostrophic Equations (QGE)}
    \label{sec:qge}

In this section, we investigate the 2S-DD-VMS-ROM~\eqref{eqn:vms-rom-two-scales-3} and the new 3S-DD-VMS-ROM~\eqref{eqn:vms-rom-three-scales-7}     
in the numerical simulation of 
the quasi-geostrophic equations (QGE)
\begin{align}
\frac{\partial \omega}{\partial t}+J(\omega,\psi)-Ro^{-1}\frac{\partial \psi}{\partial x} &= Re^{-1} \Delta \omega+Ro^{-1}F,\label{eq:qge1}\\
\omega &=-\Delta \psi, \label{eq:qge2}
\end{align}
which are used to model the large scale ocean circulation~\cite{MW06,vallis2006atmospheric}.
In~\eqref{eq:qge1}--\eqref{eq:qge2}, $\omega$ is the vorticity, $\psi$ is the streamfunction, $Re$ is the Reynolds number, and $Ro$ is the Rossby number. 

\paragraph{Computational Setting}
We follow~\cite{greatbatch2000four,mou2019data,san2015stabilized,san2011approximate} and consider 
a symmetric double-gyre wind forcing given by 
\begin{align}
	F = \sin(\pi(y-1)) ,
	\label{eq:qge:forcing}
\end{align}
the computational domain 
$\Omega = [0,1]\times [0,2]$, the time domain  $[0,80]$, and the parameters $Re=450$ and $Ro = 0.0036$.
We also assume that $\psi$ and $\omega$ satisfy homogeneous Dirichlet  boundary conditions:
\begin{align} \label{eq:qge:bdry_cond}
\psi(t, x,y) =0,\qquad \omega(t, x,y)=0 \qquad \text{for}\quad (x,y)\in\partial \Omega \; \text{ and } \;  t \ge 0.
\end{align}
\paragraph{Snapshot Generation}
For the FOM discretization, we use a spectral method with a $257 \times 513$ spatial resolution and an explicit Runge-Kutta method. 
We follow~\cite{mou2019data,san2015stabilized,san2011approximate} and run the FOM on the time interval $[0, 80]$. 
The flow displays a transient behavior on the time interval $[0,10]$, and then converges to a statistically steady state on the time interval $[10,80]$.
We record the FOM solutions on the time interval $[10,80]$ every \(10^{-2}\) simulation time units, 
which ensures that the snapshots used in the construction of the ROM basis are equally spaced. 

\paragraph{ROM Construction}

To construct the ROM basis, we follow the procedure described in Section 3.2 in~\cite{mou2019data} (see also~\cite{san2015stabilized,san2011approximate}).
First, we collect $701$ equally spaced FOM vorticity snapshots in the time interval 
$[10,80]$ at equidistant time intervals.
Next, for computational efficiency, we  interpolate the FOM vorticity onto a uniform mesh with the resolution $257\times 513$ over the spatial domain $\Omega = [0,1]\times [0,2]$, i.e., with a mesh size $\Delta x =\Delta y = 1/256$.
Finally, we use the interpolated snapshots and solve the corresponding eigenvalue problem to generate the ROM basis.

To train $\tA,\tB$ (for the 2S-DD-VMS-ROM) and $\tA_L,\tB_L$ and $\tA_S,\tB_S$ (for the 3S-DD-VMS-ROM), we use the same FOM data that was used to generate the ROM basis.
Furthermore, to increase the computational efficiency of the 2S-DD-VMS-ROM and 3S-DD-VMS-ROM, we replace the $R$-dimensional FOM data with its $d$-dimensional approximation, where the parameter $d$ satisfies $1 \leq d \leq R$ (for details, see Section 5.3 in~\cite{xie2018data}, Section 4.3 in~\cite{mohebujjaman2019physically}, and Section 3.2 in~ \cite{mou2019data}).
Specifically, we replace $\btau^{FOM}$ with $ \btau^{d}$ (for the 2S-DD-VMS-ROM) and  $\btau_L^{FOM}$ and $\btau_S^{FOM}$ with   $\btau_L^{d}$ and $\btau_S^{d}$, respectively (for the 3S-DD-VMS-ROM).
In our QGE numerical simulations, we choose $d=3r$ to maintain a good balance between numerical accuracy and computational efficiency.

\subsubsection{Numerical Results}
Next, we present results for the 2S-DD-VMS-ROM~\eqref{eqn:vms-rom-two-scales-3} and the new 3S-DD-VMS-ROM~\eqref{eqn:vms-rom-three-scales-7} in the numerical simulation of the QGE~\eqref{eq:qge1}--\eqref{eq:qge2}.
For clarity of presentation, we consider only the reconstructive regime.

To assess the ROM performance, we follow~\cite{mou2019data} and use the $L^2$ error of the time-averaged ROM streamfunction over the time interval $[10, 80]$: 
\begin{equation} \label{Eq_err_streamfunc}
\overline{\mathcal{E}}(L^2) 
=
\bigl \|\overline{\psi^{FOM}({\bf x}, \cdot)}-\overline{\psi^{ROM}({\bf x, \cdot})} \bigr \|^2_{L^2} \, ,
\end{equation}
where $\overline{(\cdot)}$ denotes the time average over the time interval $[10,80]$, and ${\bf x} = (x,y)$.
In Table~\ref{table:qge-reconstructive}, for different $r$ values, we list $\overline{\mathcal{E}}(L^2)$ for the G-ROM, the 2S-DD-VMS-ROM, and the new 3S-DD-VMS-ROM.
We also list the $r_1$ values used for the 3S-DD-VMS-ROM.
These results show that, for all r values, the 2S-DD-VMS-ROM and 3S-DD-VMS-ROM are orders of magnitude (sometimes two and even three orders of magnitude) more accurate than the standard G-ROM.
Furthermore, the 3S-DD-VMS-ROM is generally more accurate than the 2S-DD-VMS-ROM: For example, for $r=10$, $r=15$, $r=20$, and $r=25$, the 3S-DD-VMS-ROM is about {\it three times more accurate} than the 2S-DD-VMS-ROM.

\begin{table}[H]
\small
	\begin{center}
		\smallskip
		\begin{tabular}{|c|c|c|c|c|c|c|c|}
			\hline
			\multicolumn{1}{|c|}{ $r$} &  
			\multicolumn{1}{c|}{ G-ROM } &  
			\multicolumn{1}{c|}{ 2S-DD-VMS-ROM } &  
			\multicolumn{2}{c|}{3S-DD-VMS-ROM}
						 \\  \cline{2-5}&{$\overline{\mathcal{E}}(L^2)$ }  
						 &{$\overline{\mathcal{E}}(L^2)$ }  
						 &{$r_1$} & { $\overline{\mathcal{E}}(L^2)$ } 
	\\
			 \hline    
$10$
&   3.734e+02
&   5.174e-01
&   5
&   1.996e-01
 	 \\
			 \hline    
			
$15$
&   1.035e+02
&   3.853e-01
&   8
&   1.260e-01
 	 \\
			 \hline    			
			
$20$
&   1.371e+01
&   1.653e-01
&   9
&   5.175e-02
 	 \\
			 \hline   
$25$ 
&3.491e+00
&3.434e-01
&10
&5.640e-02
 	 \\
			 \hline    

		\end{tabular}
	\end{center}
	\caption {
		\label{table:qge-reconstructive}
	{
        QGE, $Re = 450$, $Ro=0.0036$, reconstructive regime. $L^2$ errors of the time-averaged streamfunction for G-ROM, 2S-DD-VMS-ROM, and 3S-DD-VMS-ROM for different r values.
        }
    	} 
\end{table}

In Figure~\ref{fig:qge-ke}, for $r = 10$ and $r = 20$, we plot the time evolution of the kinetic energy of the FOM, the G-ROM, the 2S-DD-VMS-ROM, and the new 3S-DD-VMS-ROM in the reconstructive regime. 
These plots support the conclusions in Table~\ref{table:qge-reconstructive}: 
For $r = 10$, the G-ROM kinetic energy takes off very quickly and stabilizes at a
level which is roughly $200$ times higher than the FOM kinetic energy on average. 
In contrast, both the 2S-DD-VMS-ROM and the  3S-DD-VMS-ROM produce kinetic energies of the same order of magnitude as the FOM kinetic energy. 
Furthermore, the 3S-DD-VMS-ROM performs better than the 2S-DD-VMS-ROM in reproducing the peaks and the peak frequencies.
As expected, for larger $r$ values, the G-ROM's performance improves.
For example, for $r = 20$, the G-ROM, the 2S-DD-VMS-ROM, and the 3S-DD-VMS-ROM kinetic energies perform similarly.
We note, however, that for later times (e.g., on the time interval $[60, 80]$), the G-ROM kinetic energy is somewhat higher than the FOM kinetic energy, while the 2S-DD-VMS-ROM and 3S-DD-VMS-ROM kinetic
energies are closer to the FOM kinetic energy.  

\begin{figure}[H]
\centering
	\includegraphics[width=\textwidth,height=0.3\textwidth]{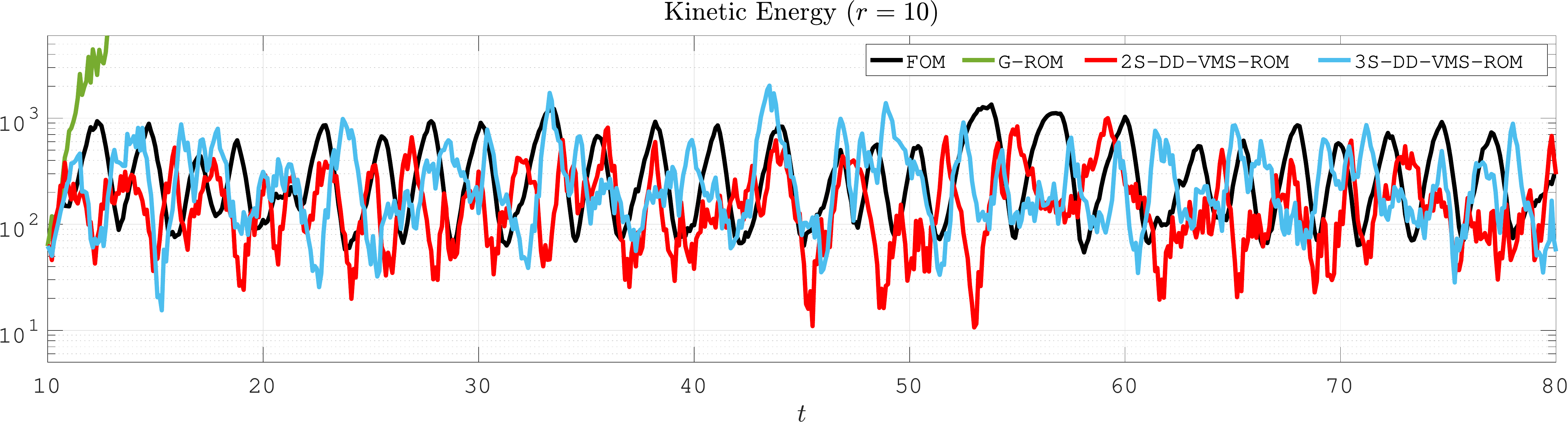}
	\includegraphics[width=\textwidth,height=0.3\textwidth]{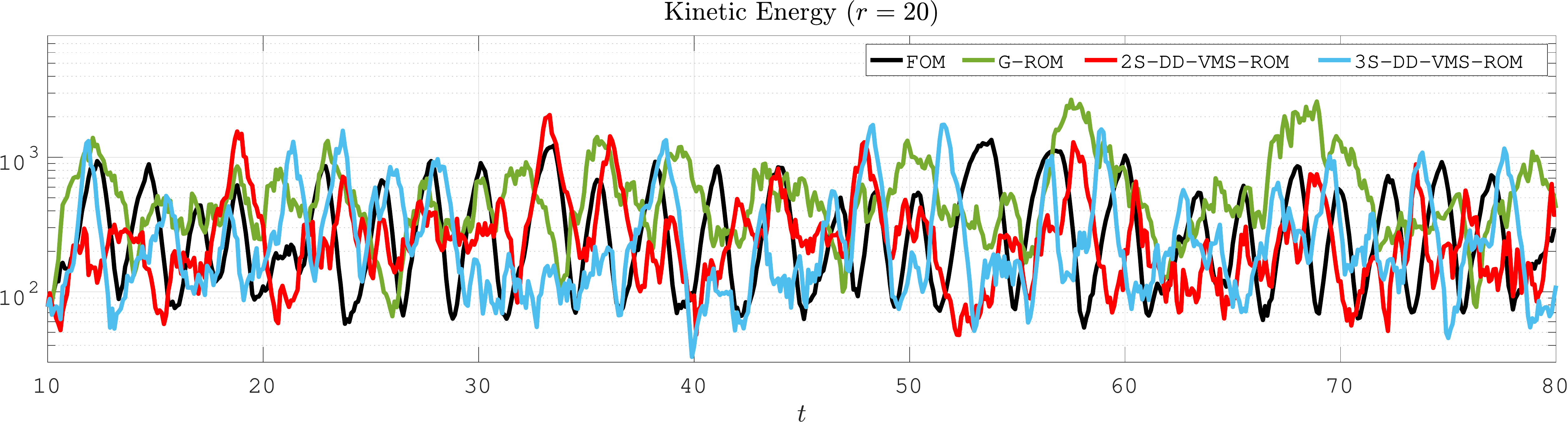}
	\caption{
	\label{fig:qge-ke}
	{
	    QGE, $Re = 450$, $Ro=0.0036$, reconstructive regime. Time evolution of the kinetic energy for FOM, G-ROM, 2S-DD-VMS-ROM, and 3S-DD-VMS-ROM for different $r$ values.
        }
        }
\end{figure}

We follow~\cite{mou2019data,san2015stabilized} and, in Figure~\ref{fig:qge-mean-psi-reconstructive}, for $r = 10$ and $r = 20$, we plot the time-average of the streamfunction $\psi$ over the time interval $[10,80]$ for the FOM,  G-ROM, 2S-DD-VMS-ROM, and 3S-DD-VMS-ROM. 
We emphasize that we use the same scale for the FOM, 2S-DD-VMS-ROM, and 3S-DD-VMS-ROM plots. 
The plots in Figure~\ref{fig:qge-mean-psi-reconstructive}  support the conclusions in Table~\ref{table:qge-reconstructive}: 
For both $r=10$ and $r=20$, the 2S-DD-VMS-ROM and 3S-DD-VMS-ROM successfully reproduce the four gyre structure in the time-averaged streamfunction, whereas the G-ROM fails. Furthermore, the 3S-DD-VMS-ROM is more accurate than the 2S-DD-VMS-ROM. 

\begin{figure}[H]
\centering
		\includegraphics[width=0.24\textwidth]{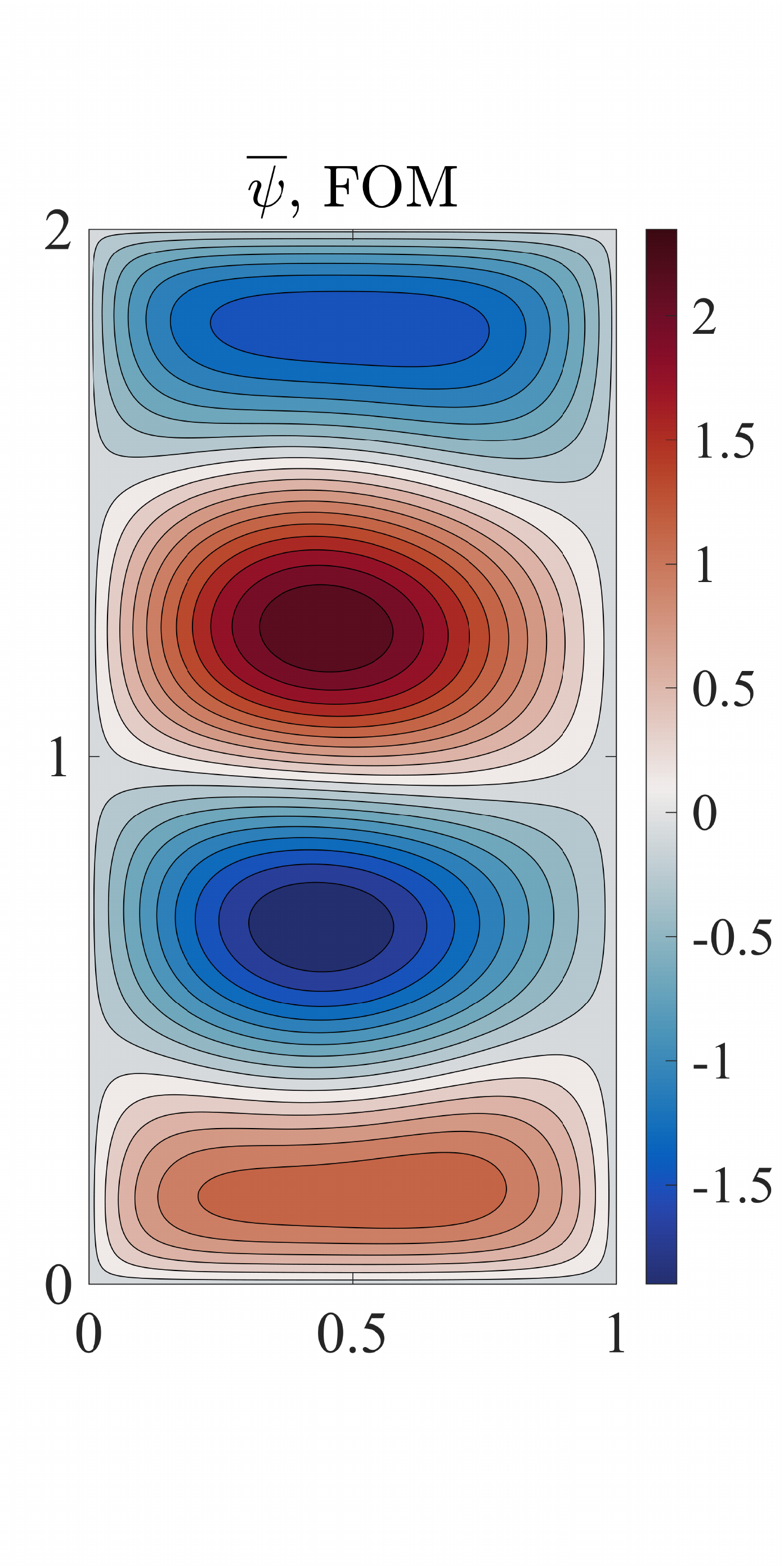}
		\includegraphics[width=0.24\textwidth]{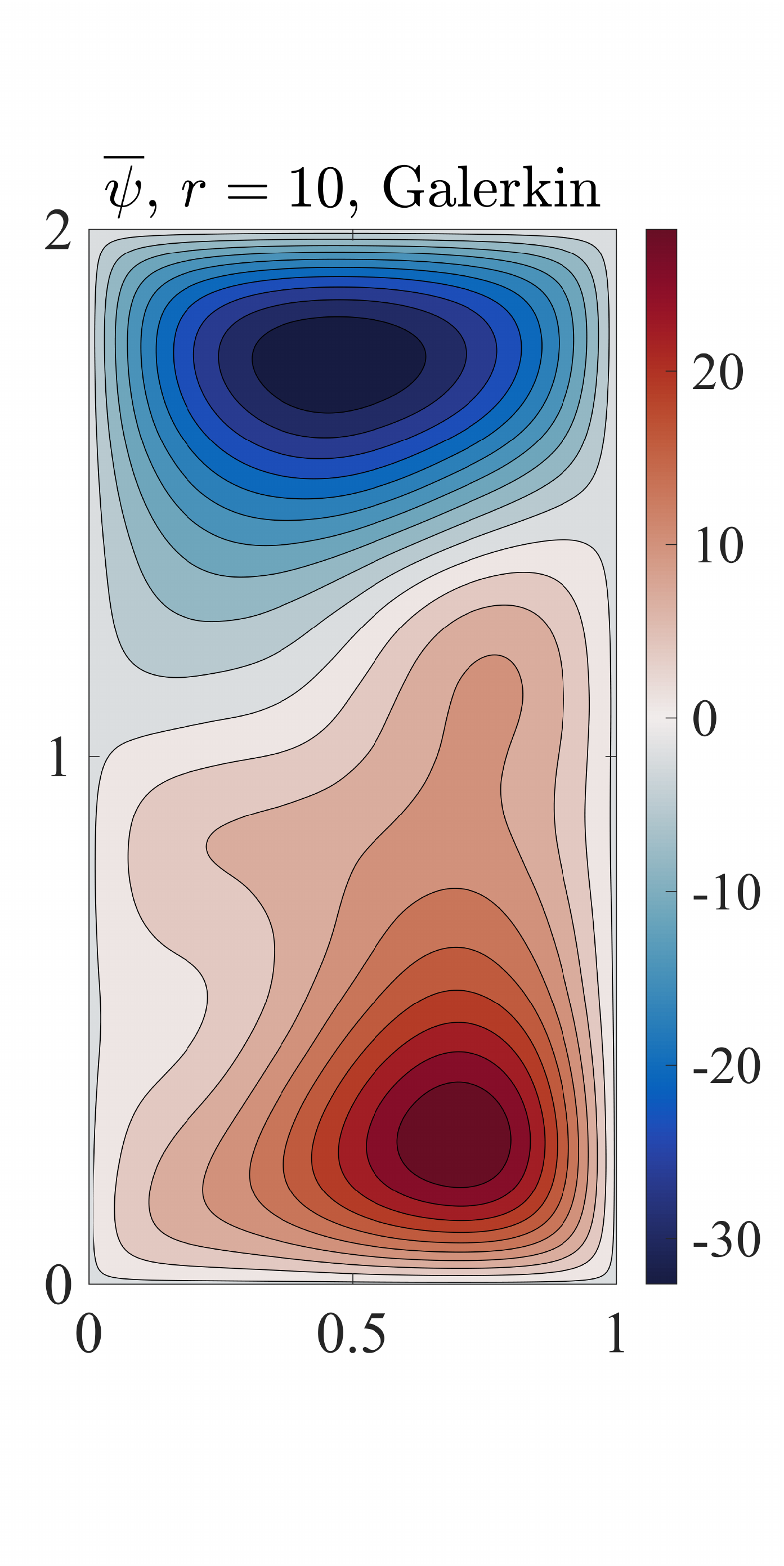}
		\includegraphics[width=0.24\textwidth]{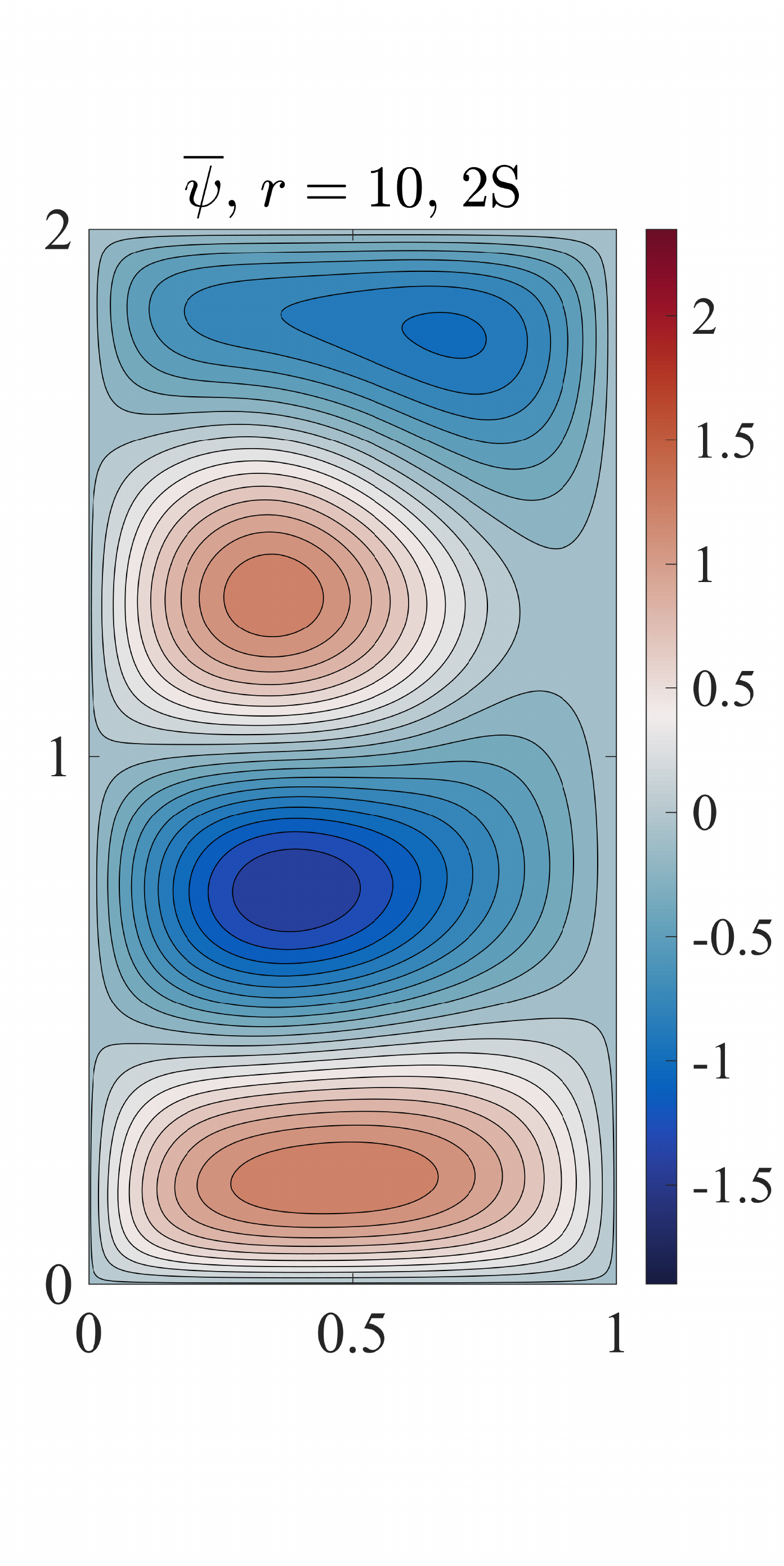}
		\includegraphics[width=0.24\textwidth]{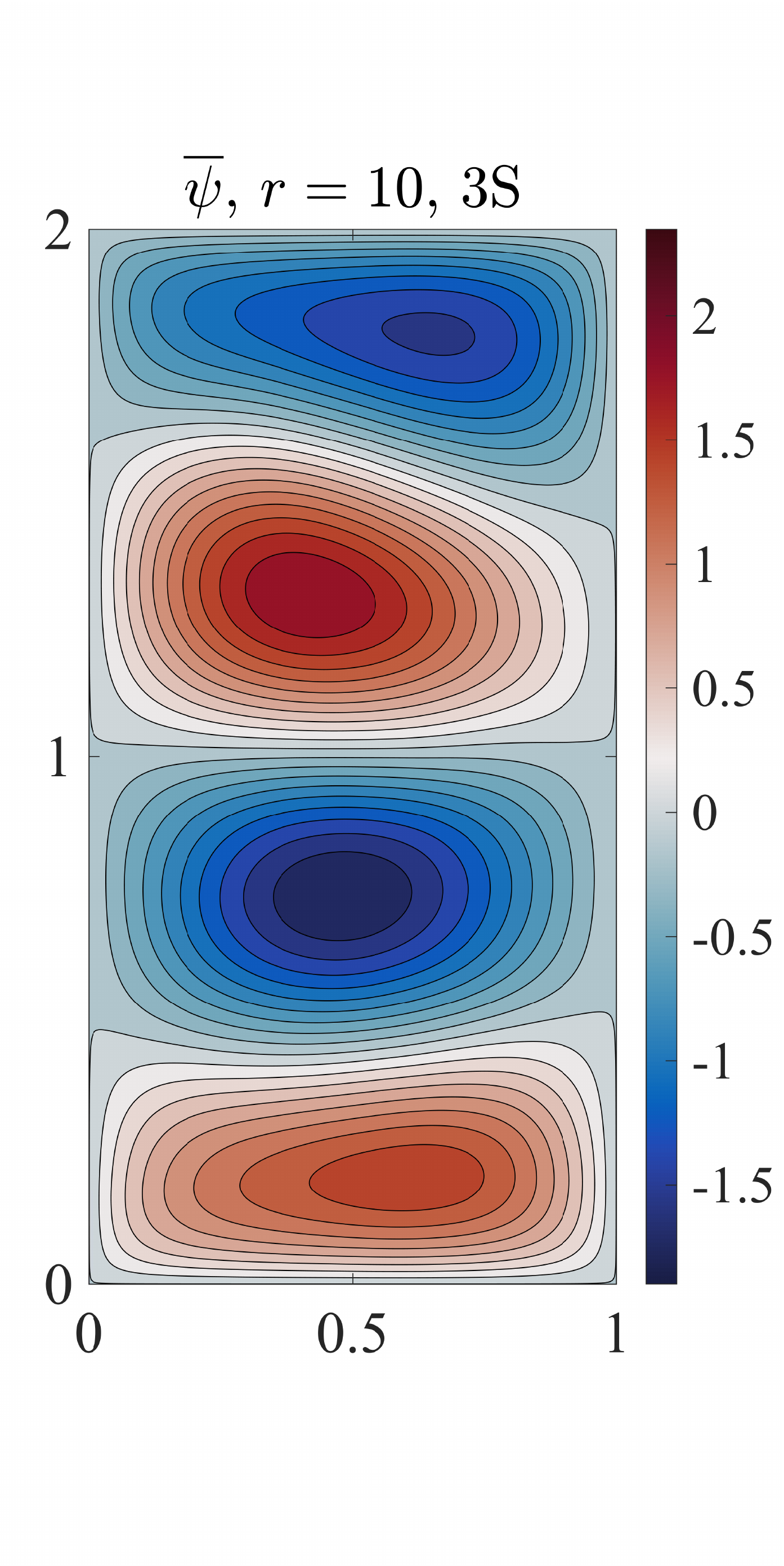}
		\\
		\vspace{-0.8cm}
		\includegraphics[width=0.24\textwidth]{figures/figures_qge/mean_psi/psi_mean_DNS.pdf}
		\includegraphics[width=0.24\textwidth]{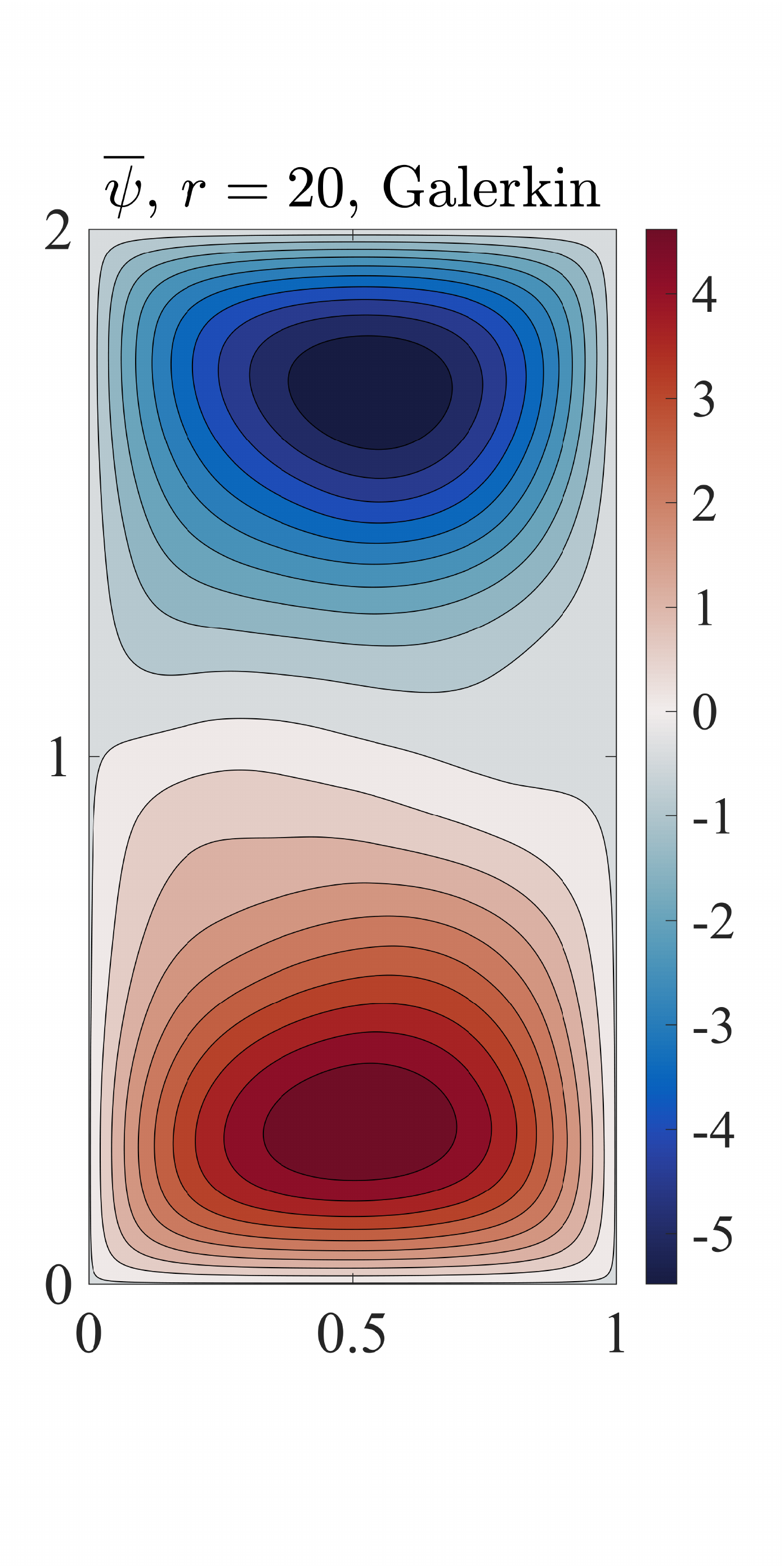}
		\includegraphics[width=0.24\textwidth]{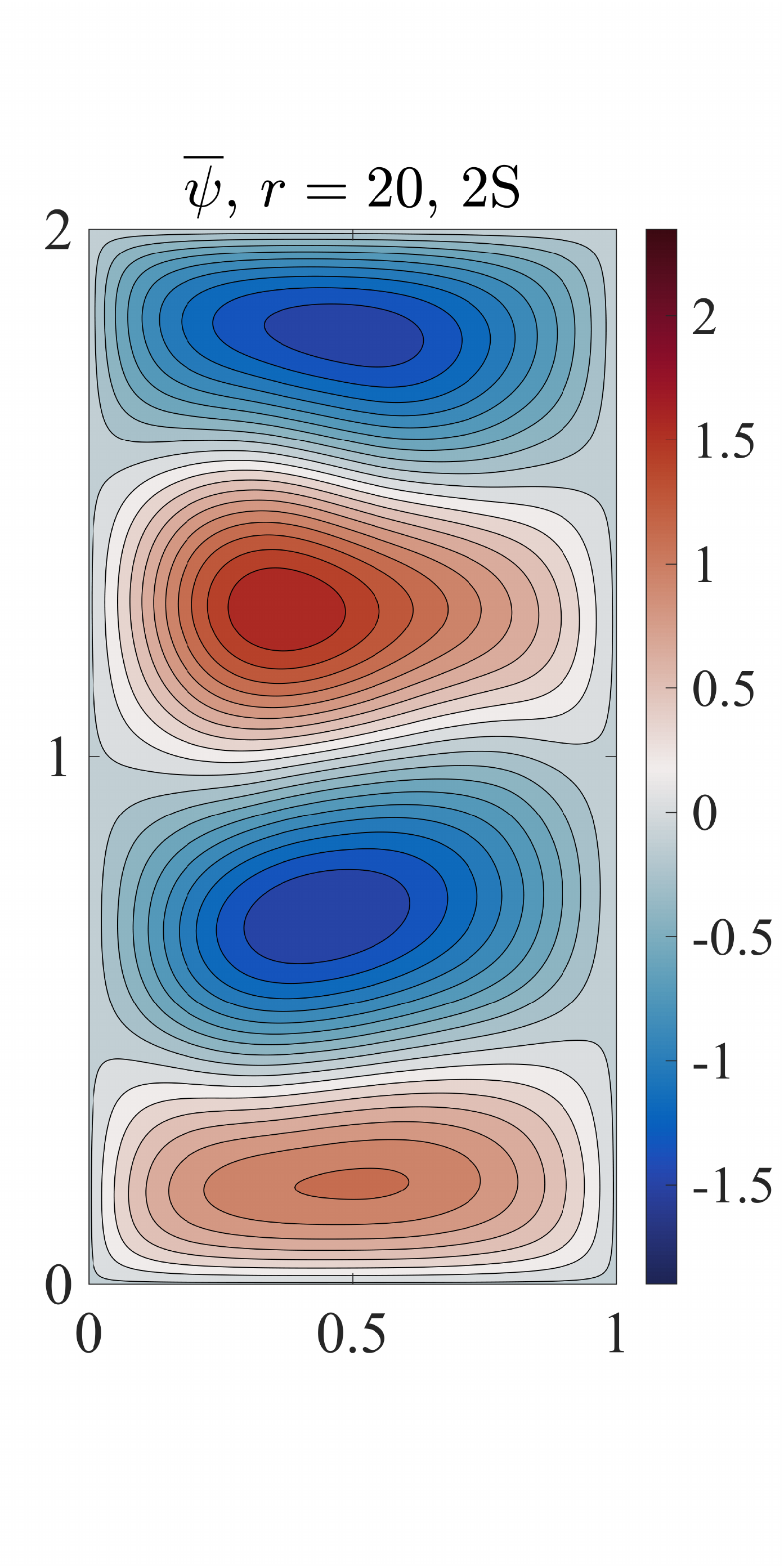}
		\includegraphics[width=0.24\textwidth]{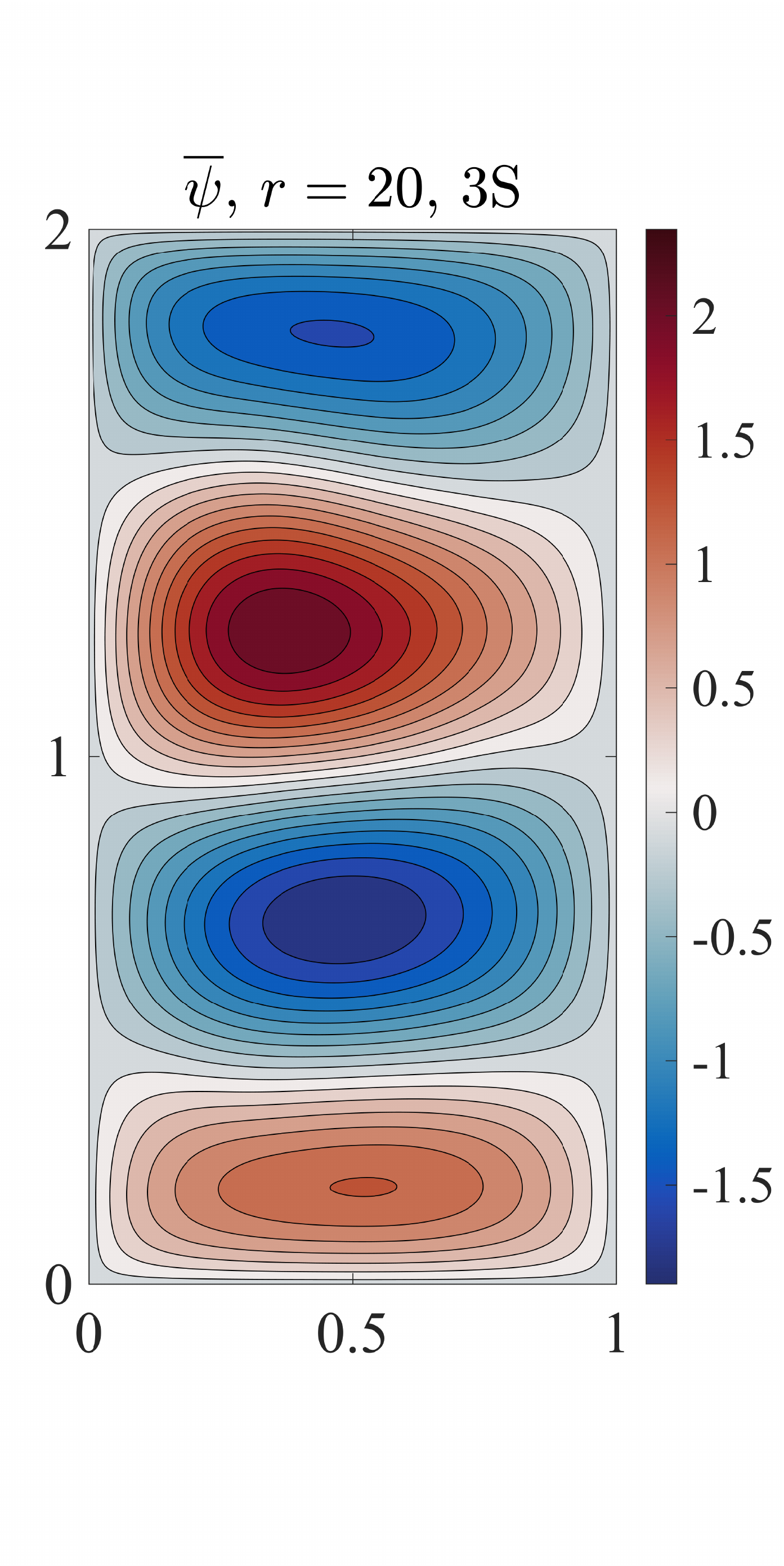}
		\\
		\vspace{-0.8cm}
	\caption{
	\label{fig:qge-mean-psi-reconstructive}
	{
        QGE, $Re = 450$, $Ro=0.0036$, reconstructive regime. Time-averaged streamfunction $\psi$ over the interval $[10, 80]$ for FOM, G-ROM, 2S-DD-VMS-ROM, and 3S-DD-VMS-ROM for different $r$ values.
        }
        }
\end{figure}

The errors listed in Table~\ref{table:qge-reconstructive} and the plots in Figures~\ref{fig:qge-ke}--\ref{fig:qge-mean-psi-reconstructive} show that, in the reconstructive regime, the 3S-DD-VMS-ROM is consistently the most accurate ROM.

}

\clearpage


{
\subsection{Backward Facing Step}
    \label{sec:bfs}
    
In this section, we investigate the 2S-DD-VMS-ROM~\eqref{eqn:vms-rom-two-scales-3} and the new 3S-DD-VMS-ROM~\eqref{eqn:vms-rom-three-scales-7}  in the numerical simulation of a two-dimensional flow over a backward facing step at $Re=1000$. 

\paragraph{Computational Setting}
As a mathematical model, we use the NSE~\eqref{eqn:nse-1}--\eqref{eqn:nse-2}.
We use the same computational domain as that used in Section 4.4~\cite{baiges2015reduced} and Section 8.2.2 in~\cite{reyes2020projection}, i.e., a $44\times 9$ rectangle with a unit height step placed at $(4, 0)$ (see the top plot in Figure~\ref{fig:bfs-geometry}).

\begin{figure}[H]
	\begin{center}
		\includegraphics[width=1\textwidth]{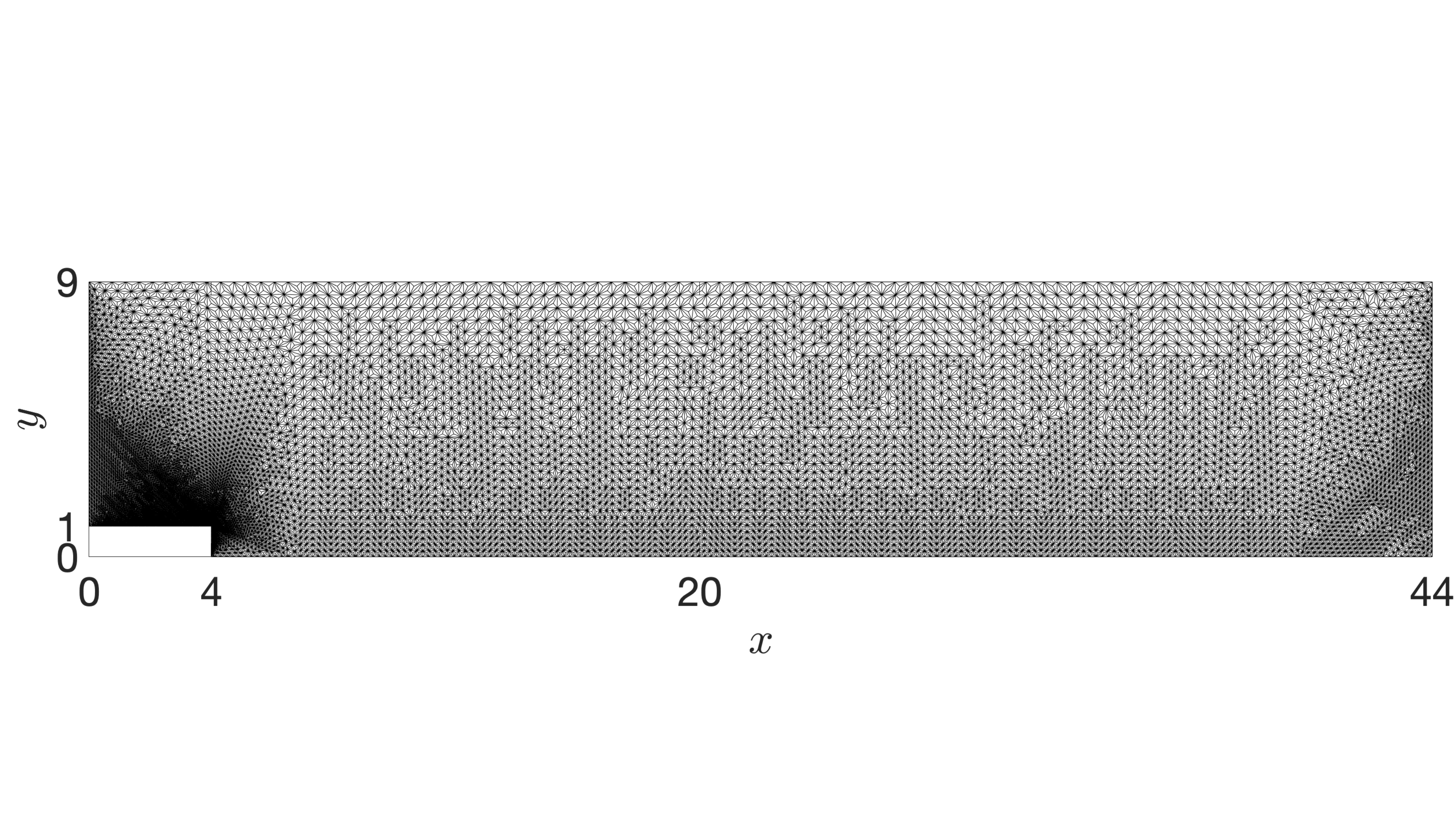} \\
		\includegraphics[width=1\textwidth]{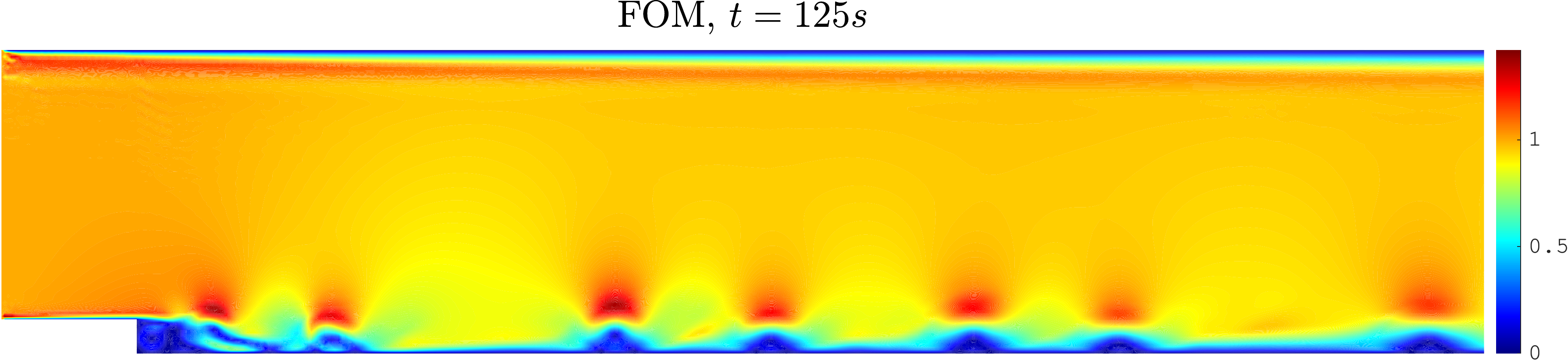}
	\end{center}
	\caption{
	        \label{fig:bfs-geometry}
	        {
            Backward facing step, $Re=1000$.
            Geometry and finite element mesh (top).
            Magnitude of FOM velocity field at $t = 125$  (bottom).
            }
            }
\end{figure}

\paragraph{Snapshot Generation} 
For the spatial discretization, we use a barycenter refinement mesh of a Delaunay generated triangulation, which allows for $(P_2,P_1^{disc})$ Scott-Vogelius
elements to be LBB stable (for details, see \cite{john2017divergence}). 
The mesh (see the top plot in~Figure~\ref{fig:bfs-geometry}) has $209508$ velocity and $156285$ pressure degrees of freedom.  
We use the linearized BDF2 method and a time step size $\Delta t=0.05$ for both FOM and ROM time discretizations. 
On the first time step, we use the backward Euler method so that we have two initial time step solutions required for the BDF2 scheme. 
For illustration purposes, in Figure~\ref{fig:bfs-geometry} (the bottom plot), we display the magnitude of the FOM velocity field at $t = 125$.

In Figure~\ref{fig:bfs-dns-ke}, we plot the time evolution of the FOM kinetic energy on the time interval $[100,150]$. 
This plot shows that the flow over a backward facing step that we consider is not periodic or periodic-like.
The numerical results in the remainder of this section will show that this setting is more challenging for reduced order modeling than the other three test problems considered in Sections~\ref{sec:numerical-results-burgers}--\ref{sec:qge}.

\begin{figure}[H]
	\begin{center}
		\includegraphics[width=0.9\textwidth,height=0.25\textwidth]{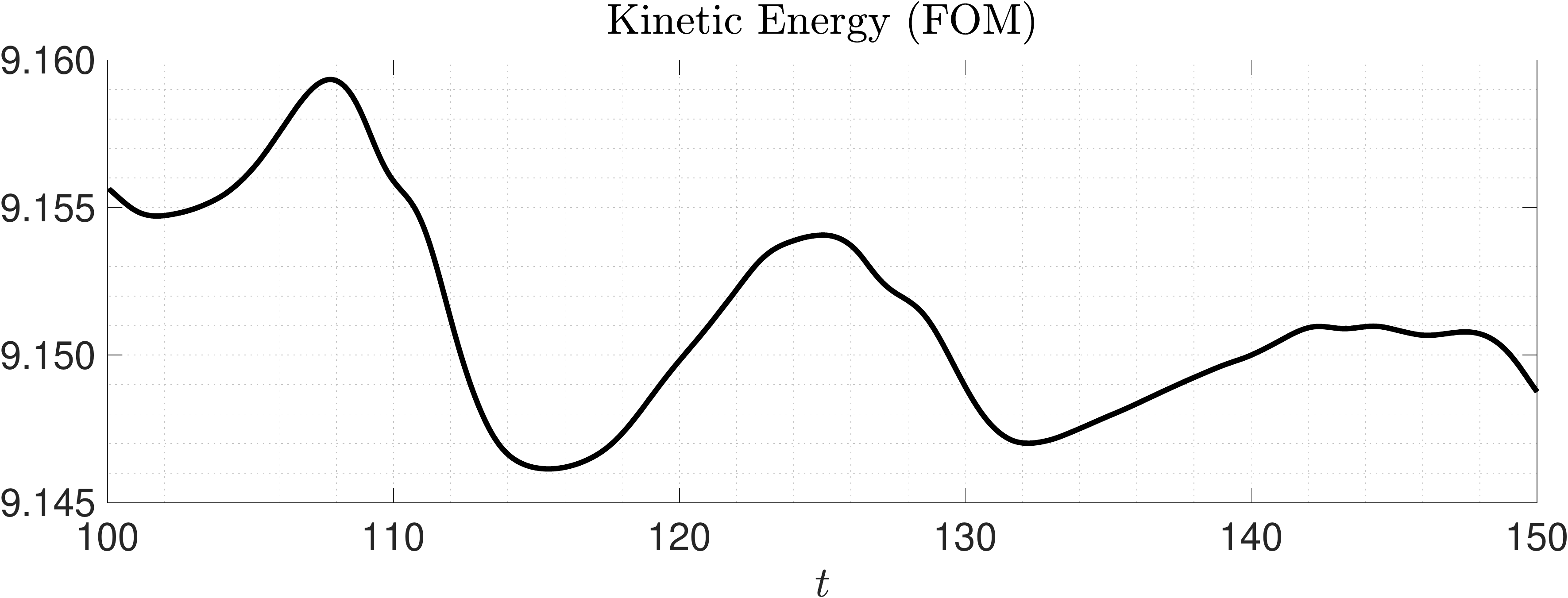}
	\end{center}
	\caption{
	\label{fig:bfs-dns-ke}
	{
Backward facing step, $Re=1000$. 
Time evolution of the FOM kinetic energy.
		}
		}
\end{figure}

\paragraph{ROM Construction}
To build the ROM basis functions, we follow~\cite{reyes2020projection} and collect $1000$) equally spaced FOM snapshots on the time interval $[100.05,150]$.

To train $\tA,\tB$ (for the 2S-DD-VMS-ROM) and $\tA_L,\tB_L$ and $\tA_S,\tB_S$ (for the 3S-DD-VMS-ROM), we use the same FOM data that was used to generate the ROM basis.
Furthermore, to increase the computational efficiency of the 2S-DD-VMS-ROM and 3S-DD-VMS-ROM, we use the  approach described in Section~\ref{sec:qge} and replace $\btau^{FOM}$ with $ \btau^{3r}$ (for the 2S-DD-VMS-ROM) and  $\btau_L^{FOM}$ and $\btau_S^{FOM}$ with   $\btau_L^{3r}$
and $\btau_S^{3r}$, respectively (for the 3S-DD-VMS-ROM).
To further reduce the computational cost of the 3S-DD-VMS-ROM, we adopt a generic way in choosing $r_1$ for large $r$ values (i.e., $r\ge 30$) and let $r_1 = \lfloor r/2 \rfloor$.

\subsubsection{Numerical Results}

Next, we present results for the 2S-DD-VMS-ROM~\eqref{eqn:vms-rom-two-scales-3} and the new 3S-DD-VMS-ROM~\eqref{eqn:vms-rom-three-scales-7} in the numerical simulation of the flow over a backward facing step at $Re=1000$.
For clarity of presentation, we consider only the reconstructive regime.

In Table~\ref{table:bfs-re1000-reconstructive}, for different $r$ values, we list the average $L^2$ error~\eqref{eqn:l2-error} for the G-ROM, the 2S-DD-VMS-ROM, and the new 3S-DD-VMS-ROM. 
We also list the $r_1$ values for the 3S-DD-VMS-ROM. 
These results show that, for all $r$ values, both the 2S-DD-VMS-ROM and the 3S-DD-VMS-ROM are about $30\%$ more accurate than the standard G-ROM.
Furthermore, the 3S-DD-VMS-ROM is consistently more accurate than the 2S-DD-VMS-ROM.
This improvement is significant for low $r$ values (i.e., $2 \leq r \leq 15$), and modest for large $r$ values (i.e., $20 \leq r \leq 60$).

\begin{table}[h!]
\small
	\begin{center}
		\smallskip
		\begin{tabular}{|c|c|c|c|c|c|c|c|}
			\hline
			\multicolumn{1}{|c|}{ $r$} &  
			\multicolumn{1}{c|}{ G-ROM } &  
			\multicolumn{1}{c|}{ 2S-DD-VMS-ROM } &  
			\multicolumn{2}{c|}{3S-DD-VMS-ROM}
						 \\  \cline{2-5}&{$\mathcal{E}(L^2)$ }
						&{$\mathcal{E}(L^2)$} 
						 &{$r_1$} & { $\mathcal{E}(L^2)$ } 
									 \\
			 \hline 
  
$2$
&    1.0270e+00
&    9.6593e-01
& 1
&    8.6129e-01
 	 \\
			 \hline    
$5$
&   1.4864e+00
&   1.1671e+00
& 1
& 1.1070e+00
 	 \\
			 \hline    
			
$10$
&   1.8401e+00
&   1.5064e+00
& 2
& 1.2932e+00 
 	 \\
			 \hline    
			
$15$
&   1.4733e+00
&   1.0909e+00
& 9
&   7.4297e-01
 	 \\
			 \hline    			
			
$20$
&   1.0392e+00
&   7.5813e-01
& 3
&  7.0704e-01
 	 \\
			 \hline    

			
$30$
& 9.1723e-01  
& 7.7908e-01 
& 15
& 7.5835e-01
 	 \\
			 \hline    			
			
			
$40$
&  4.4118e-01 
& 2.9694e-01  
& 20
& 2.7753e-01
 	 \\
			 \hline    			
			
			
$50$
&  2.5578e-01 
&  1.6002e-01  
& 25
& 1.5586e-01
 	 \\
			 \hline    			
			
			
$60$
& 1.6772e-01 
& 1.1679e-01   
& 30
& 1.1276e-01
\\	
            \hline
		\end{tabular}
	\end{center}
	\caption {
	{ 
		Backward facing step, $Re = 1000$, reconstructive regime. 
		Average $L^2$ errors for G-ROM, 2S-DD-VMS-ROM, and 3S-DD-VMS-ROM for different $r$ values. 
	\label{table:bfs-re1000-reconstructive}
	}
	} 
\end{table}

We follow~\cite{baiges2015reduced} (see also~\cite{reyes2020projection}) and, in Figure~\ref{fig:bfs-pt-vy-re1000}, for $r=15$, we plot a pointwise quantity, i.e., the time evolution of the $y$-component of the velocity, $v$, for the FOM, 2S-DD-VMS-ROM, and 3S-DD-VMS-ROM at the  point with coordinates $(19,1)$, which is physically located behind the step. 

This plot shows that both the 2S-DD-VMS-ROM and the 3S-DD-VMS-ROM are significantly more accurate than the G-ROM.
Furthermore, the 3S-DD-VMS-ROM is more accurate than the 2S-DD-VMS-ROM.

\begin{figure}[H]
	\begin{center}
	\includegraphics[width=\textwidth]{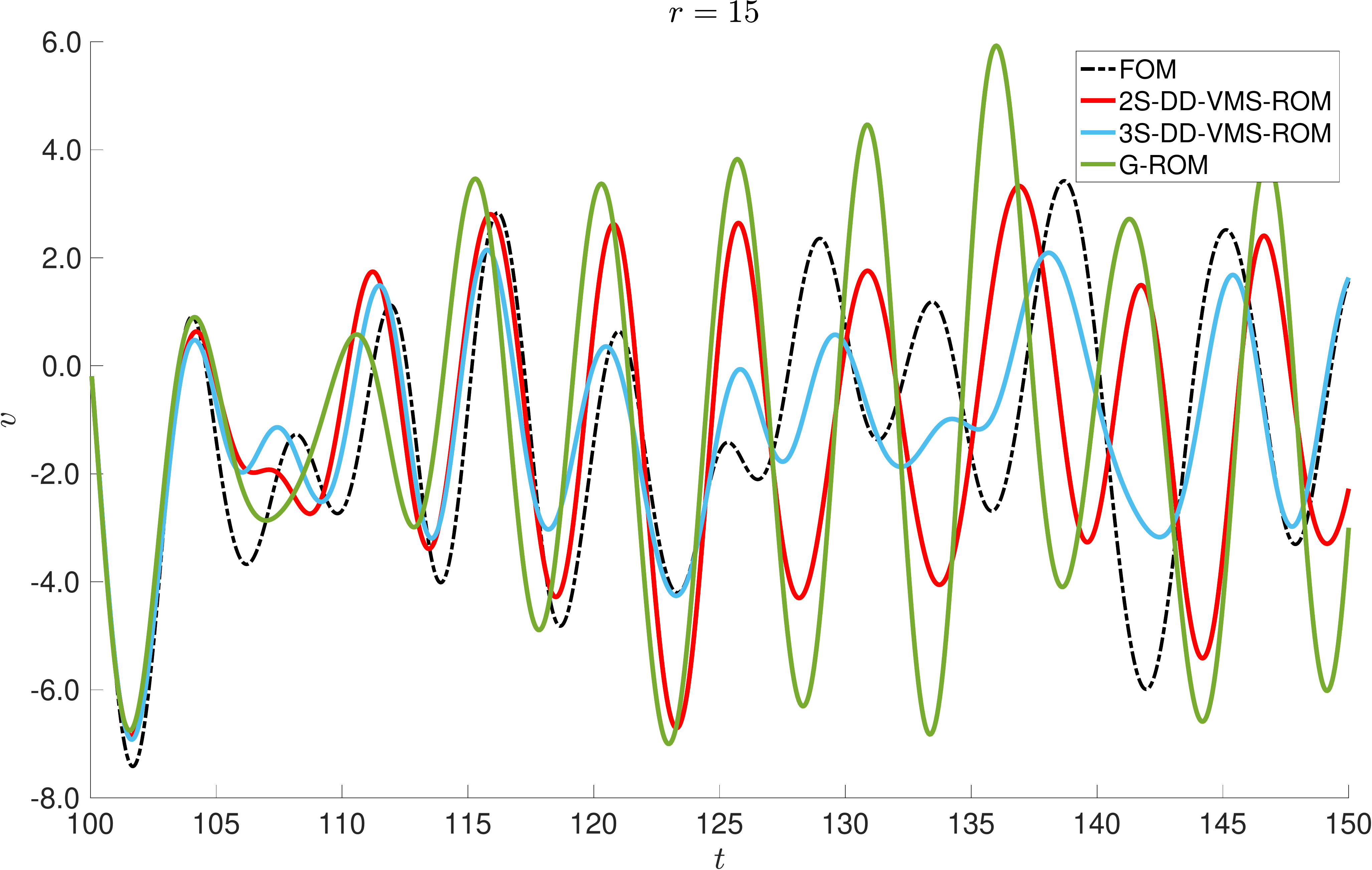}

	\end{center}
	\caption{		
	\label{fig:bfs-pt-vy-re1000}
	{
        	Backward facing step, $Re=1000$, reconstructive regime.
			Time evolution of the $y$-component of the velocity, $v$, of FOM, G-ROM, 2S-DD-VMS-ROM, and 3S-DD-VMS-ROM with $r=15$ at the  point with coordinates $(19,1)$.  
		}
		}
\end{figure}

In Figure~\ref{fig:bfs-ke-reconstructive}, for $r=30, 40$, and $60$,
we plot the time evolution of the kinetic energy of the FOM, G-ROM, 2S-DD-VMS-ROM, and 3S-DD-VMS-ROM. 
These plots support the conclusions in Table~\ref{table:bfs-re1000-reconstructive}.
Specifically, for low $r$ values (i.e., $r=30$), the G-ROM, 2S-DD-VMS-ROM, and 3S-DD-VMS-ROM results are relatively inaccurate.
However, for medium $r$ values (i.e., $r=40$), the 2S-DD-VMS-ROM and 3S-DD-VMS-ROM results are significantly more accurate than the G-ROM results.
As expected, for high $r$ values  (i.e., $r=60$),  the 2S-DD-VMS-ROM, 3S-DD-VMS-ROM, and G-ROM results are all accurate.
Furthermore, for $r=40$ the 3S-DD-VMS is more accurate than the 2S-DD-VMS-ROM.
For $r=30$ and $r=60$, the 2S-DD-VMS-ROM and the 3S-DD-VMS-ROM perform similarly.

\begin{figure}[H]
	\begin{center}
		\includegraphics[width=0.9\textwidth,height=0.25\textwidth]{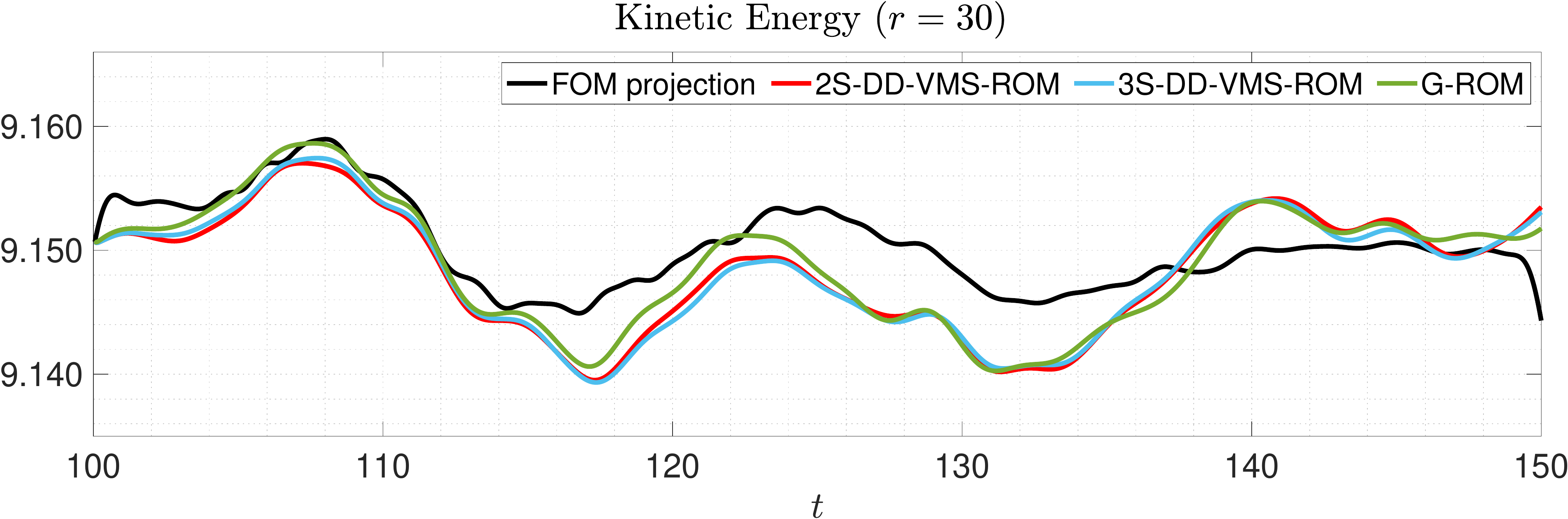}
		\includegraphics[width=0.9\textwidth,height=0.25\textwidth]{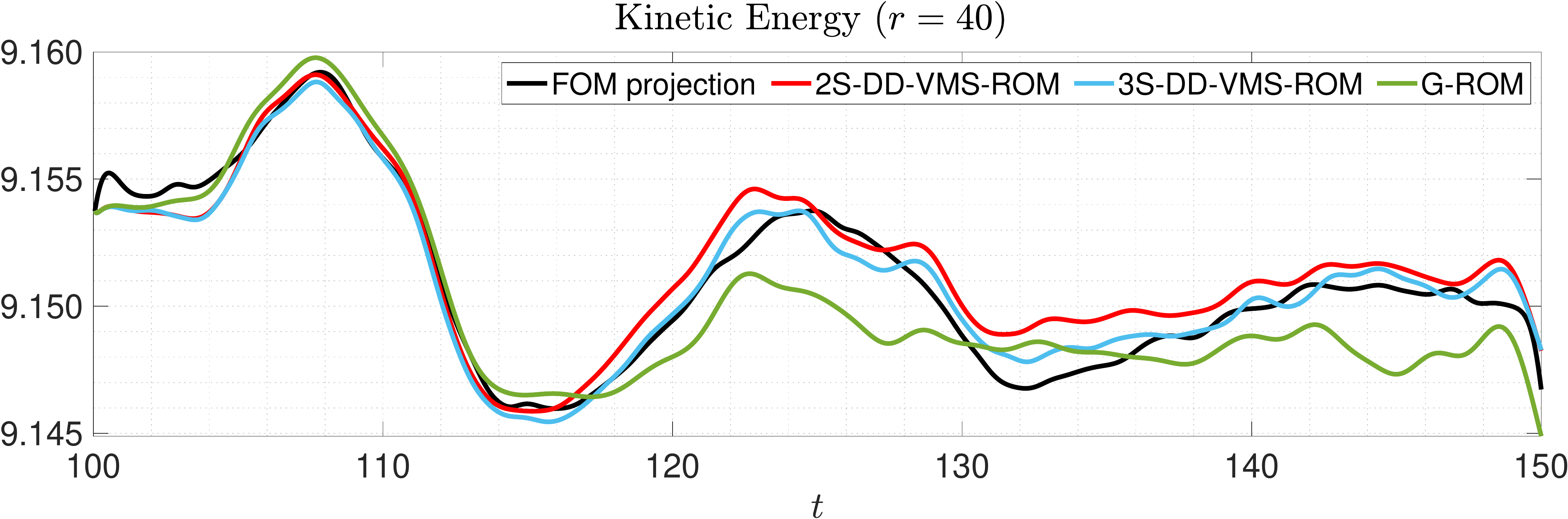}
		\includegraphics[width=0.9\textwidth,height=0.25\textwidth]{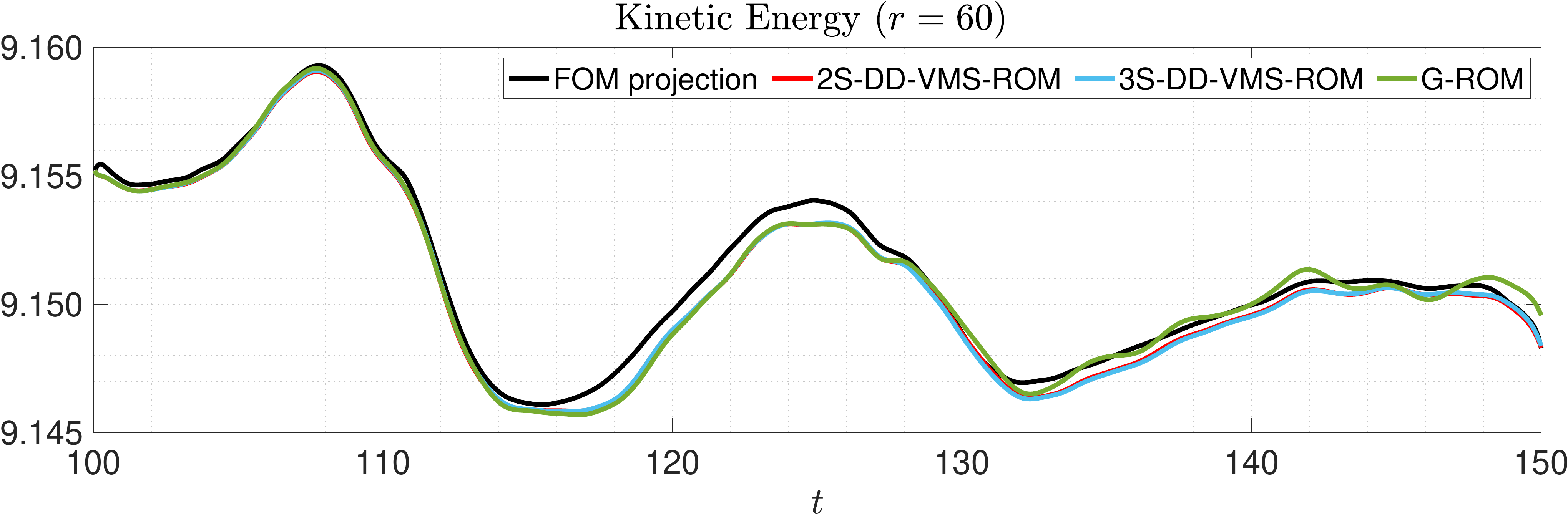}
	\end{center}
	\caption{
	\label{fig:bfs-ke-reconstructive}
	{
Backward facing step, $Re=1000$, reconstructive regime. 
Time evolution of the kinetic energy for FOM projection, G-ROM, 2S-DD-VMS-ROM and 3S-DD-VMS-ROM for different $r$ values. 
		}
			}
\end{figure}

We follow~\cite{baiges2015reduced,reyes2020projection} and, in Figure~\ref{fig:bfs-velocity-spectrum}, for $r=15$, we plot a pointwise quantity, i.e., the spectrum of the $y-$component of the velocity, $v$, for the FOM, 2S-DD-VMS-ROM, and 3S-DD-VMS-ROM at the point with coordinates $(19,1)$. 
This plot shows that the 2S-DD-VMS-ROM spectrum is more accurate than the G-ROM spectrum.
Furthermore, the 3S-DD-VMS-ROM spectrum is more accurate than the 2S-DD-VMS-ROM spectrum.

\begin{figure}[H]
	\begin{center}
	\includegraphics[width=0.8\textwidth]{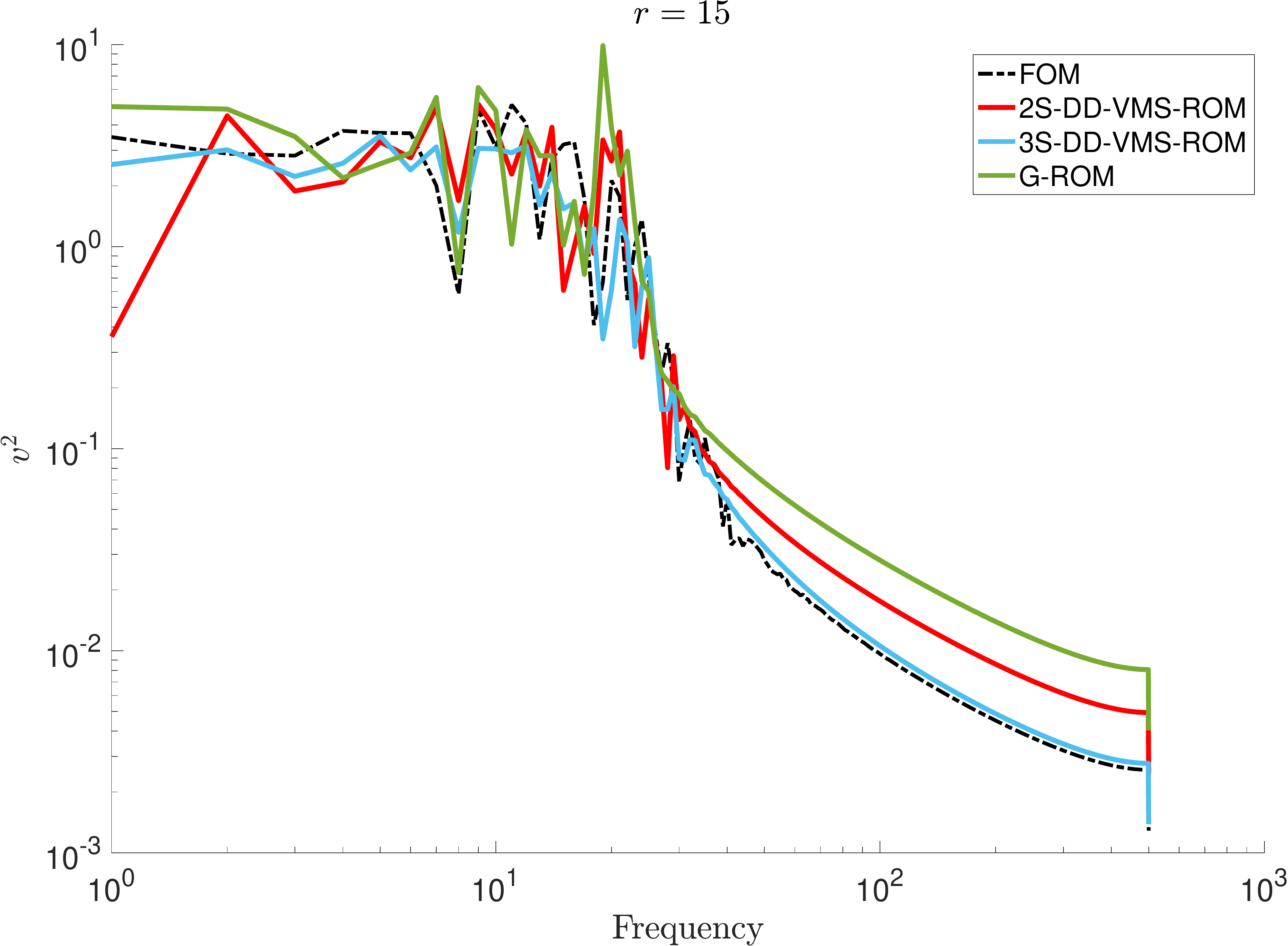}

	\end{center}
	\caption{
	\label{fig:bfs-velocity-spectrum}
	{
        	Backward facing step, $Re=1000$, reconstructive regime. 
        	The spectrum of the $y$-component of the velocity for  FOM, G-ROM, 2S-DD-VMS-ROM, and 3S-DD-VMS-ROM with  $r=15$ at the  point with coordinates $(19,1)$.
		}
		}
\end{figure}

The errors listed in Table~\ref{table:bfs-re1000-reconstructive} and the plots in Figures~\ref{fig:bfs-pt-vy-re1000}--\ref{fig:bfs-velocity-spectrum} show that, in the reconstructive regime, both the 2S-DD-VMS-ROM and the 3S-DD-VMS-ROM are more accurate than the G-ROM. 
Furthermore, the 3S-DD-VMS-ROM is more accurate than the 2S-DD-VMS-ROM.
However, for the backward facing step test problem, this improvement is not as significant as for the other three test cases investigated in Sections~\ref{sec:numerical-results-burgers}--\ref{sec:qge}.

}

\clearpage


{
\subsection{Qualitative Comparison of 2S-DD-VMS-ROM and 3S-DD-VMS-ROM}

In the previous sections, we performed a quantitative comparison of the 2S-DD-VMS-ROM and the 3S-DD-VMS-ROM in the numerical simulation of the Burgers equation (Section~\ref{sec:numerical-results-burgers}), the flow past a cylinder (Section~\ref{sec:numerical-results-nse}),
the QGE (Section~\ref{sec:qge}), and the flow over a backward facing step (Section~\ref{sec:bfs}).
In all our numerical simulations, the 3S-DD-VMS-ROM was more accurate than the 2S-DD-VMS-ROM, although this improvement was less significant for the flow over a backward facing step.
In this section, we present a qualitative comparison of the 2S-DD-VMS-ROM and 3S-DD-VMS-ROM.

We believe that the 3S-DD-VMS-ROM is more accurate than the 2S-DD-VMS-ROM in our numerical tests because the 3S-DD-VMS-ROM is more {\it flexible} than the 2S-DD-VMS-ROM.
Specifically, the 2S-DD-VMS-ROM has only one control parameter in the truncated SVD used in Algorithm~\ref{alg:2s-dd-vms-rom}, i.e., the tolerance $tol$.
The 3S-DD-VMS-ROM, on the other hand, has two control parameters in the truncated SVD used in Algorithm~\ref{alg:3s-dd-vms-rom}: (i) the tolerance $tol_L$ used in the truncated SVD for the  least squares problem for the large resolved scales, and (ii) the tolerance $tol_S$ used in the truncated SVD for the least squares problem for the small resolved scales.
Thus, in principle, by choosing optimal values for the two modeling parameters in the 3S-DD-VMS-ROM (i.e., $tol_L$ and $tol_S$), we can obtain more accurate results than those obtained with the 2S-DD-VMS-ROM, which has only one modeling parameter (i.e., $tol$).
The truncated SVD components of the 2S-DD-VMS-ROM and 3S-DD-VMS-ROM algorithms aim at alleviating the ill-conditioning that is common in data-driven least squares problems (see, e.g.,~\cite{peherstorfer2016data,mohebujjaman2019physically,xie2018data}).
Our numerical investigation shows that the tolerances used in the truncated SVD have a significant effect on the 2S-DD-VMS-ROM and 3S-DD-VMS-ROM results.
Furthermore, our numerical results confirm that having more flexibility in choosing the two tolerances in the 3S-DD-VMS-ROM yields more accurate results.

For example, for the Burgers equation, the results in Table~\ref{table:burgers-reconstructive-3a} show that, for $r=3$, choosing two different tolerances in the 3S-DD-VMS-ROM (i.e., $tol_{L}=10^{0}$ and $tol_{S}=10^{-2}$) yields more accurate results than the 2S-DD-VMS-ROM, which uses only one tolerance (i.e., $tol=10^{0}$).
Indeed, the 3S-DD-VMS-ROM average $L^2$ error is more than an order of magnitude lower than the 2S-DD-VMS-ROM average $L^2$ error.

The flow past a circular cylinder test case yields similar conclusions. 
We follow~\cite{baiges2015reduced} and, in Figure~\ref{fig:nse-pt-vy-re1000}, for $r=5$, we plot the time evolution of the $y$-component of the velocity, $v$, of the FOM, 2S-DD-VMS-ROM, and 3S-DD-VMS-ROM at the  point with coordinates $(0.43,0.2)$, which is physically located behind the circular cylinder. 
The plot in Figure~\ref{fig:nse-pt-vy-re1000} clearly shows that choosing two different tolerances in the 3S-DD-VMS-ROM algorithm yields more accurate results than the 2S-DD-VMS-ROM, which uses only one tolerance.
\begin{figure}[H]
	\begin{center}
			\includegraphics[width=0.9\textwidth]{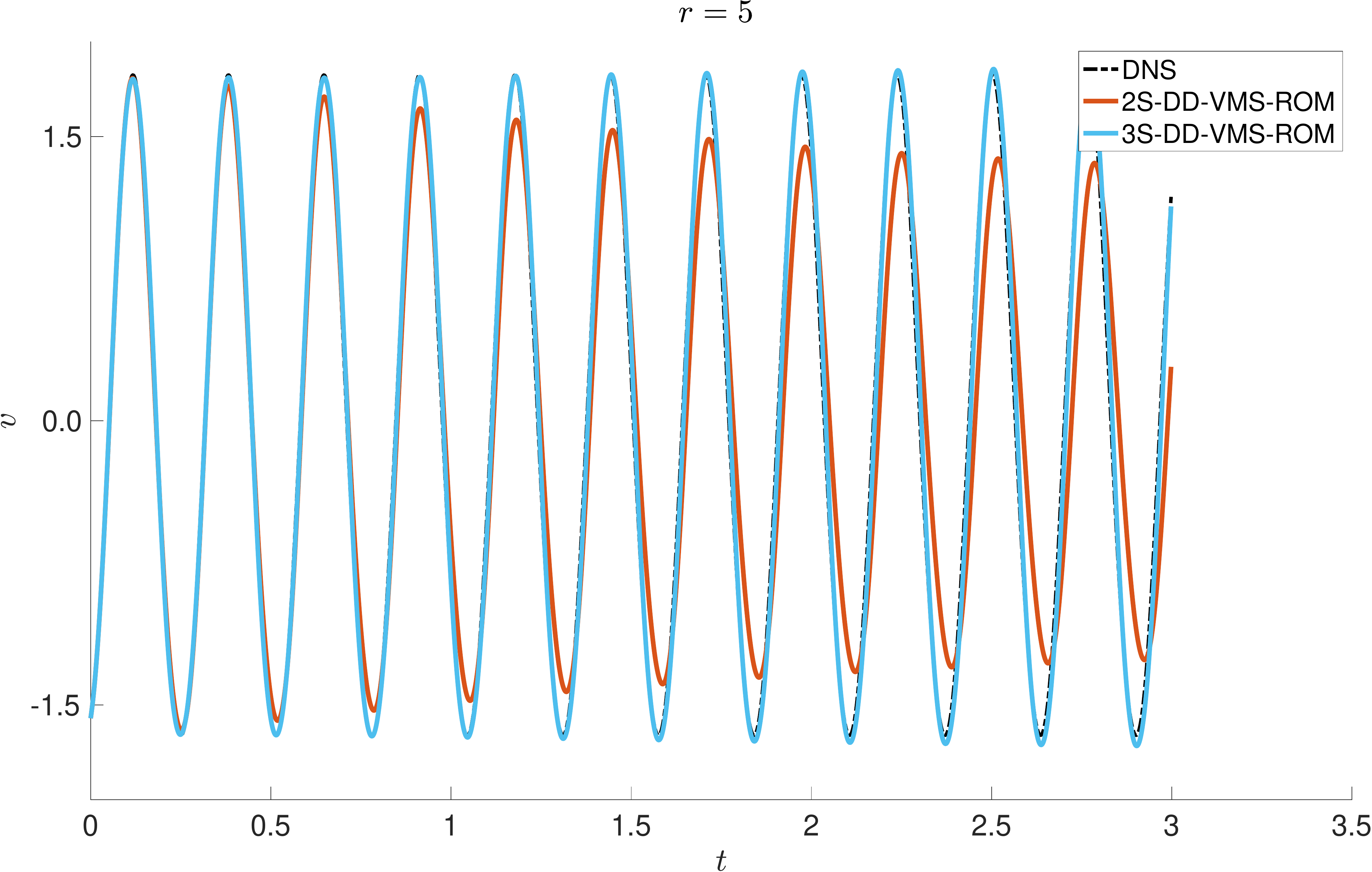}
	\end{center}
	\caption {
	        {
			Flow past a cylinder, $Re=1000$, reconstructive regime.
			Time evolution of the $y$-component of the velocity, $v$, of the FOM, 2S-DD-VMS-ROM, and 3S-DD-VMS-ROM with  $r=5$ at the  point with coordinates $(0.43,0.2)$.
			}
	\label{fig:nse-pt-vy-re1000}
	}
\end{figure}

Furthermore, we follow~\cite{baiges2015reduced} and, for the 2S-DD-VMS-ROM and  3S-DD-VMS-ROM, in Figure~\ref{fig:nse-subscale-re1000} we plot the first component of the vectors $\btau^{FOM}$ and $\btau^{ROM}$ with the FOM and ROM representations of the VMS-ROM closure terms, which are defined in~\eqref{eqn:ansatz-nse} for the 2S-DD-VMS-ROM and in~\eqref{eqn:vms-rom-three-scales-3}--\eqref{eqn:vms-rom-three-scales-4} for the 3S-DD-VMS-ROM.
Specifically, at each time step $t_j, \, j=1,\ldots,M$,
\begin{eqnarray}
    \btau^{FOM}(t_j) 
    =  &-& \bigl[ 
			\bigl( \bigl( {\bu_R^{FOM}(t_j)} \cdot \nabla \bigr) \, {\bu_R^{FOM}(t_j)} \, , \bphi_{i} \bigr)
			\nonumber \\[0.3cm]
			&& - \bigl( \bigl( {\bu_r^{FOM}(t_j)} \cdot \nabla \bigr) \, {\bu_r^{FOM}(t_j)} \, , \bphi_{i} \bigr)
		\bigr]\,,
\end{eqnarray}
where $\bu_R^{FOM}(t_j)$ and $\bu_r^{FOM}(t_j)$ are defined in~\eqref{eqn:u-r-fom}, and
\begin{eqnarray}
    \btau^{ROM}(t_j)  
    = \tA \, \ba^{ROM}(t_j) 
    + \ba^{ROM}(t_j)^\top\,\tB\,\ba^{ROM}(t_j)\,,
\label{eq:subscale-rom}
\end{eqnarray}
where $\tA$ and $\tB$ are the DD-VMS-ROM operators, and $\ba^{ROM}(t_j)$ is the ROM solution at time step $t_j$.
The plot in Figure~\ref{fig:nse-subscale-re1000} shows that the first component of the 2S-DD-VMS-ROM and 3S-DD-VMS-ROM closure terms are different. 
Thus, we conclude that the tolerance used in the truncated SVD has a significant effect on the ROM closure model and on the corresponding ROM results (as illustrated in Figure~\ref{fig:nse-pt-vy-re1000}). 

\begin{figure}[H]
	\begin{center}
     \includegraphics[width=0.7\textwidth]{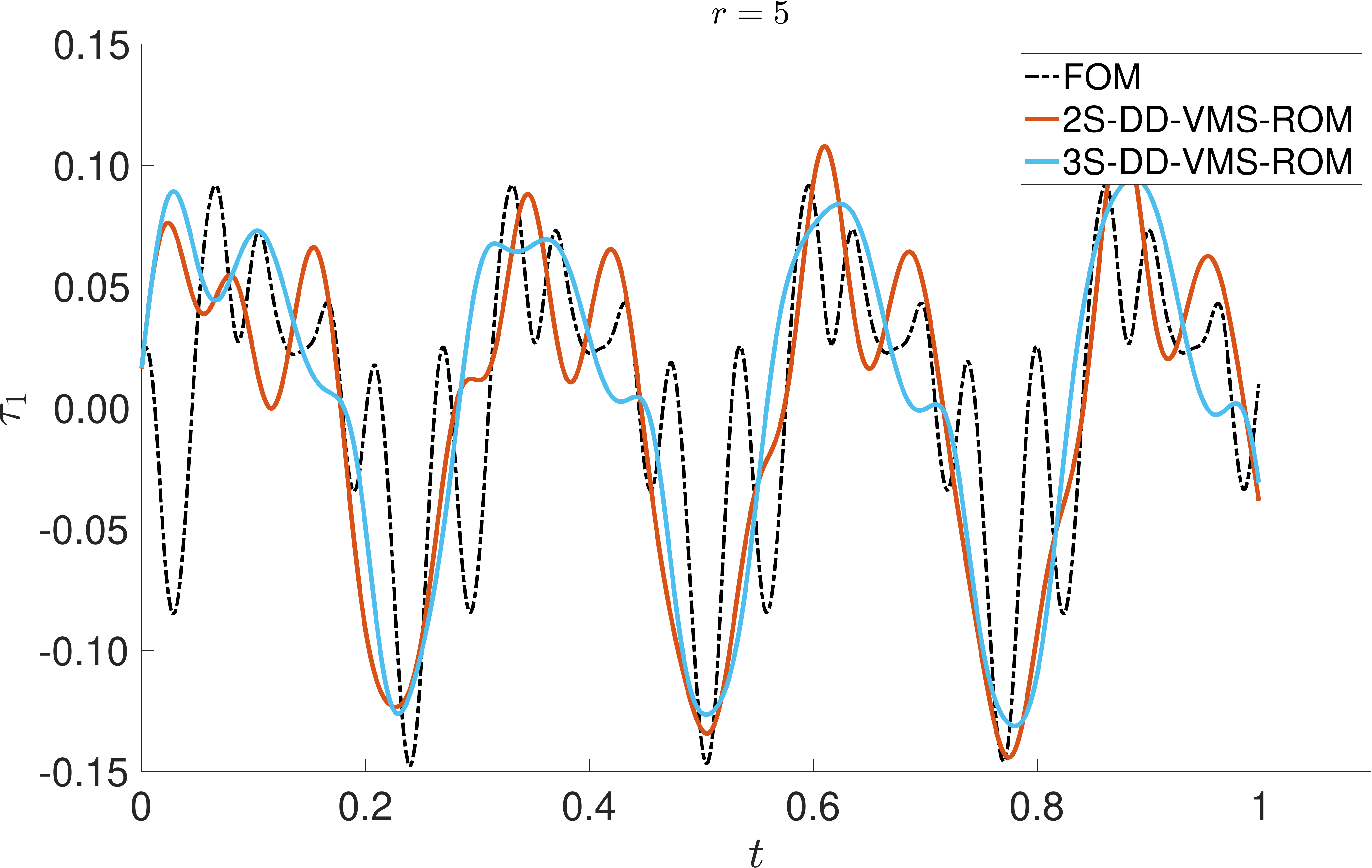}
	\end{center}
	\caption {{
	        {
			Flow past a cylinder, $Re=1000$, reconstructive regime.
			Time evolution of the first component  of the subscales for the FOM, 2S-DD-VMS-ROM, and 3S-DD-VMS-ROM with  $r=5$.
			}
	\label{fig:nse-subscale-re1000}
	}
	}
\end{figure}
{
For the QGE test case, the results in Table~\ref{table:qge-reconstructive} show that choosing two different tolerances in the 3S-DD-VMS-ROM yields more accurate results than the 2S-DD-VMS-ROM, which uses only one tolerance. 
For example, for $r=25$, the 3S-DD-VMS-ROM  $L^2$ error
is more than six times lower than the 2S-DD-VMS-ROM $L^2$ error.

For the backward facing step test case, the results in Table~\ref{table:bfs-re1000-reconstructive} (see also 
Figures~\ref{fig:bfs-pt-vy-re1000}--\ref{fig:bfs-velocity-spectrum}) support the same conclusion. 
Indeed, the 3S-DD-VMS-ROM (which uses two different tolerances) is 
more accurate than the 2S-DD-VMS-ROM (which uses only one tolerance).
This improvement is significant for low $r$ values (i.e., $2 \leq r \leq 15$), and modest for large $r$ values (i.e., $20 \leq r \leq 60$).

}

\bigskip

We emphasize that both the quantitative comparisons (in Sections~\ref{sec:numerical-results-burgers}--\ref{sec:bfs}) and the qualitative comparison in this section are only valid for the 2S-DD-VMS-ROM and the 3S-DD-VMS-ROM.  
Thus, our conclusions do not carry over to other types of VMS-ROMs, e.g.,~\cite{baiges2015reduced,bergmann2009enablers,eroglu2017modular,guler2019decoupled,iliescu2013variational,iliescu2014variational,reyes2020projection,reyes2018reduced,roop2013proper,stabile2019reduced,tello2019fluid,wang2012proper}.
In particular, we do not perform a general comparison of two-scale VMS-ROMs and three-scale VMS-ROMs.
Instead, we take a more modest step and compare two specific examples from the two classes, i.e., the 2S-DD-VMS-ROM and the 3S-DD-VMS-ROM, respectively.
We believe that extending to the ROM setting two-scale and three-scale VMS models developed for classical numerical discretizations (see, e.g., the surveys in~\cite{ahmed2017review,codina2018variational,john2016finite,rasthofer2018recent}), and comparing the resulting two-scale and three-scale VMS-ROMs is a worthy research endeavor that could yield conclusions that are different from the conclusions drawn from our numerical investigation (see, e.g.,~\cite{ahmed2020assessment} for the finite element setting).
This, however, is beyond the scope of this paper.

}


\section{Conclusions and Outlook}
	\label{sec:conclusions}

In this paper, we propose a new data-driven variational multiscale reduced order model (DD-VMS-ROM) framework.
We construct the new DD-VMS-ROM framework in two steps:
In the first step, we leverage the VMS methodology and the hierarchical structure of the ROM basis to provide explicit mathematical formulas for the interaction among the ROM spatial scales.
In the second step, we use the available full order model (FOM) data to construct structural VMS-ROM closure models for the interactions among scales. 
We investigate two DD-VMS-ROMs:
(i) The two-scale DD-VMS-ROM (2S-DD-VMS-ROM) considers two scales: resolved scales and unresolved scales.
For the 2S-DD-VMS-ROM, we construct one ROM closure model for the interaction between the resolved and unresolved scales.
(ii) The three-scale DD-VMS-ROM (3S-DD-VMS-ROM) considers three scales: resolved large scales, resolved small scales, and unresolved scales.
For the 3S-DD-VMS-ROM, we construct one ROM closure model for the interaction between the resolved large and resolved small scales, and another ROM closure model for the interaction between resolved small scales and unresolved scales.
We test the 2S-DD-VMS-ROM and 3S-DD-VMS-ROM in the numerical simulation of
{four test cases:
(i) the 1D Burgers equation with  viscosity  coefficient $\nu = 10^{-3}$; 
(ii) a 2D flow past a circular cylinder at Reynolds numbers $Re=100$, $Re=500$, and $Re=1000$;
(iii) the quasi-geostrophic equations at Reynolds number $Re=450$ and Rossby number $Ro=0.0036$; and 
(iv) a 2D flow over a backward facing step at Reynolds number $Re=1000$.
We consider the reconstructive regime for all the test cases, and the cross-validation and predictive regimes for the Burgers equation and  the 2D flow past a circular cylinder test cases.
}
The numerical results show that both the 2S-DD-VMS-ROM and the 3S-DD-VMS-ROM are more accurate than the standard Galerkin ROM (G-ROM). 
{
Furthermore, the 3S-DD-VMS-ROM is more accurate than the 2S-DD-VMS-ROM, although this improvement is less significant for the flow over a backward facing step.
}

We intend to pursue several research avenues in the development of the new DD-VMS-ROM  framework.
The first research direction that we plan to investigate is finding the optimal parameter $r_1$ and the optimal tolerances $tol_{L}$ and $tol_{S}$ in the new 3S-DD-VMS-ROM.
In this paper, we used a trial and error approach to find these parameters.
We intend to investigate whether providing rigorous error estimates~\cite{hesthaven2015certified,iliescu2014are,quarteroni2015reduced} or leveraging physical insight~\cite{HLB96} can provide parameters that yield more accurate results.
Another research direction that we plan to pursue is the development of new DD-VMS-ROM closure models by leveraging ideas from VMS methods for finite element discretizations (see, e.g., Section 8.8 in~\cite{john2016finite}), e.g., the time-dependent subscale-orthogonal methods~\cite{codina2002stabilized,reyes2020projection,reyes2018reduced}.
We also plan to explore different topological structures for the ROM closure term. 
In the present study, we assume that the structure of the ROM closure model function ${\boldsymbol{g}}$ is similar to the structure of the Galerkin model function $\bff$ and we utilize a least squares approach to determine the shape of ${\boldsymbol{g}}$. 
We emphasize that, without loss of generality, our DD-VMS-ROM framework can be formulated by utilizing a supervised machine learning approach~\cite{rahman2018hybrid,san2018extreme,san2018machine,san2018neural}, a topic that we would like to explore in the future.
Finally, we intend to explore the extension of the new DD-VMS-ROM to the Petrov-Galerkin framework~\cite{carlberg2017galerkin,carlberg2011efficient,grimberg2020stability,parish2019adjoint}.

\section*{Acknowledgements}

The work of the first, second, and fourth authors was supported by National Science Foundation grants DMS-1821145 and BMMB-1929731.
The work of the third author was supported by the U.S. Department of Energy, Office of Science, Office of Advanced Scientific Computing (ASCR), under Award number DE-SC0019290.

\bibliographystyle{plain}
\bibliography{reference,qge}

\end{document}